\documentclass[11pt,reqno]{amsart}
\usepackage{amsmath,amssymb,amsthm,mathtools}
\usepackage{tikz-cd}
\usepackage{tikz}
\usepackage{enumitem}
\usepackage{upgreek}
\usepackage[margin=1in]{geometry}
\definecolor{darkblue}{rgb}{0.0, 0.2, 0.8}
\usepackage[linkcolor=darkblue,colorlinks=true,citecolor=darkblue]{hyperref}
\theoremstyle{definition}
\newtheorem{defn}{Definition}[section]
\newtheorem{rem}[defn]{Remark}
\newtheorem{hypotheses}[defn]{Hypotheses}
\theoremstyle{plain}
\newtheorem{thm}[defn]{Theorem}
\newtheorem{prop}[defn]{Proposition}
\newtheorem{lem}[defn]{Lemma}
\newtheorem{lemma}[defn]{Lemma}
\newtheorem{cor}[defn]{Corollary}
\newtheorem{conjecture}[defn]{Conjecture}
\newcommand{\F}{\mathbb{F}_2}
\newcommand{\sr}{\textsc{r}}
\newcommand{\e}{\mathsf{e}}
\newcommand{\rrsc}{\textsc{rrsc}}
\newcommand{\RRSC}{\textsc{rrsc}}

\newcommand{\red}{\mathrm{red}}
\newcommand{\scone}{\textsc{i}}
\newcommand{\sctwo}{\textsc{ii}}
\newcommand{\HMR}{\text{\small \textit{HMR}}}
\newcommand{\hmr}{\text{\small \textit{HMR}}}
\newcommand{\hfr}{\text{\small \textit{HFR}}}
\newcommand{\HM}{\text{\small \textit{HM}}}
\newcommand{\cmr}{\text{\small \textit{CMR}}}
\newcommand{\swr}{\text{\small \textit{SWR}}}
\newcommand{\BF}{\text{\small \textit{BF}}}
\newcommand{\tcmr}{\widetilde{\cmr}}
\newcommand{\tH}{\widetilde{\HMR}}

\newcommand{\toH}{\widecheck{\HMR}}
\newcommand{\thmr}{\widetilde{\HMR}}
\newcommand{\tHMR}{\widetilde{\HMR}}
\newcommand{\grtwo}{\operatorname{gr}_{}^{\scriptscriptstyle (2)}}
\DeclareMathOperator{\dbcv}{\mathsf{\Sigma}}
\DeclareMathOperator{\id}{\mathrm{id}}
\DeclareMathOperator{\rot}{\mathsf{rot}}
\DeclareMathOperator{\fix}{\mathsf{Fix}}
\DeclareMathOperator{\im}{im}
\DeclareMathOperator{\nullity}{null}
\DeclareMathOperator{\grad}{grad}
\DeclareMathOperator{\coker}{coker}
\DeclareMathOperator{\Arf}{Arf}
\DeclareMathOperator{\Hess}{Hess}
\DeclareMathOperator{\ev}{\mathrm{ev}}

\DeclareMathOperator{\gr}{gr}
\DeclareMathOperator{\Gr}{Gr}
\newcommand\hcancel[2][black]{\setbox0=\hbox{$#2$}%
    \rlap{\raisebox{.45\ht0}{\textcolor{#1}{\rule{\wd0}{1pt}}}}#2}
\DeclareMathOperator{\pertL}{\hcancel{\mathcal{L}}}

\makeatletter
\newbox\WC@base
\newbox\WC@hat
\DeclareRobustCommand\widecheck[1]{{\mathpalette\@widecheck{#1}}}

\def\@widecheck#1#2{%
	\setbox\WC@base\hbox{\m@th$#1#2$}%
	\setbox\WC@hat\hbox{\m@th$#1\widehat{\vrule\@width\z@\@height\ht\WC@base\vrule\@height\z@\@width\wd\WC@base}$}%
	\dp\WC@hat=-\ht\WC@base
	\@tempdima=\ht\WC@base
	\advance\@tempdima 2\ht\WC@hat
	\divide\@tempdima\thr@@
	\setbox\WC@hat\hbox{\raise\@tempdima\hbox{\scalebox{1}[-1]{\lower\@tempdima\box\WC@hat}}}%
	{\ooalign{\box\WC@hat \cr \box\WC@base}}%
}
\makeatother
\newcommand{\figFgamma}{
\tikzset{every picture/.style={line width=0.75pt}} 
\begin{tikzpicture}[x=0.75pt,y=0.75pt,yscale=-1,xscale=1]
	\draw    (100,80) -- (100,200) ;
	\draw    (220,80) -- (220,110) ;
	\draw    (220,140) -- (220,200) ;
	\draw    (100,80) -- (220,80) ;
	\draw    (100,200) -- (220,200) ;
	\draw    (220,110) .. controls (150.66,109.92) and (150.16,139.92) .. (220,140) ;
	\draw    (220,80) .. controls (260.66,80.42) and (260.66,110.42) .. (220,110) ;
	\draw    (220,140) -- (300,140) ;
	\draw    (220,200) -- (300,200) ;
	\draw    (300,140) -- (300,160) ;
	\draw    (300,180) -- (300,200) ;
	\draw    (300,160) .. controls (249.66,160.42) and (249.66,179.92) .. (300,180) ;
	\draw    (300,140) .. controls (331.16,140.42) and (329.66,159.92) .. (300,160) ;
	\draw    (300,180) .. controls (331.16,180.42) and (329.66,199.92) .. (300,200) ;
	\draw [color={rgb, 255:red, 0; green, 0; blue, 255 }  ,draw opacity=1 ]   (220,160) .. controls (279.16,160.42) and (271.66,179.42) .. (220,180) ;
	\draw [color={rgb, 255:red, 0; green, 0; blue, 255 }  ,draw opacity=1 ]   (220,90) .. controls (239.16,89.42) and (239.16,99.92) .. (220,100) ;
	\draw [color={rgb, 255:red, 0; green, 0; blue, 255 }  ,draw opacity=1 ]   (220,100) .. controls (159.93,100.32) and (139.34,111.08) .. (140,130) .. controls (140.66,148.92) and (176.66,160.42) .. (220,160) ;
	\draw [color={rgb, 255:red, 0; green, 0; blue, 255 }  ,draw opacity=1 ]   (100,90) -- (220,90) ;
	\draw [color={rgb, 255:red, 0; green, 0; blue, 255 }  ,draw opacity=1 ]   (100,180) -- (220,180) ;
	
	\draw (101,149.4) node [anchor=north west][inner sep=0.75pt]    {$\overline{W}$};
	\draw (73,122.4) node [anchor=north west][inner sep=0.75pt]    {$Y_{00}$};
	\draw (214.5,59.9) node [anchor=north west][inner sep=0.75pt]    {$S^{3}$};
	\draw (209,202.9) node [anchor=north west][inner sep=0.75pt]    {$Y_{0}$};
	\draw (254.5,82.9) node [anchor=north west][inner sep=0.75pt]    {$B^{4}$};
	\draw (319,161.9) node [anchor=north west][inner sep=0.75pt]    {$W_{F,M}$};
	\draw (241,181.4) node [anchor=north west][inner sep=0.75pt]  [color={rgb, 255:red, 0; green, 0; blue, 255 }  ,opacity=1 ]  {$F_{0}$};
	\draw (221,111.4) node [anchor=north west][inner sep=0.75pt]  [color={rgb, 255:red, 0; green, 0; blue, 255 }  ,opacity=1 ]  {$D^{2}$};
	\draw (154.5,161) node [anchor=north west][inner sep=0.75pt]  [color={rgb, 255:red, 0; green, 0; blue, 255 }  ,opacity=1 ] [align=left] {2 bands};
\end{tikzpicture}}

\newcommand{\figseifertcap}{
\tikzset{every picture/.style={line width=0.75pt}} 
\begin{tikzpicture}[x=0.75pt,y=0.75pt,yscale=-1,xscale=1]
	\draw    (90,100) -- (90,220) ;
	\draw    (220,140) .. controls (140.75,141) and (147.25,178) .. (220,180) ;
	\draw [color={rgb, 255:red, 0; green, 0; blue, 255 }  ,draw opacity=1 ]   (90,140) .. controls (186.46,140.92) and (186.25,179.5) .. (90,180) ;
	\draw    (220,100) -- (220,140) ;
	\draw    (220,180) -- (220,220) ;
	\draw    (220,100) .. controls (289.96,99.92) and (289.96,139.92) .. (220,140) ;
	\draw    (220,180) .. controls (281.46,179.92) and (281.96,220.42) .. (220,220) ;
	\draw    (220,100) -- (148.46,100) -- (90,100) ;
	\draw    (220,220) -- (148.46,220) -- (90,220) ;
	\draw (71,152.4) node [anchor=north west][inner sep=0.75pt]    {$Y$};
	\draw (111,152.4) node [anchor=north west][inner sep=0.75pt]  [color={rgb, 255:red, 0; green, 0; blue, 255 }  ,opacity=1 ]  {$F$};
	\draw (130,195.4) node [anchor=north west][inner sep=0.75pt]    {$W_F$};
	\draw (232,112.4) node [anchor=north west][inner sep=0.75pt]    {$M$};
	\draw (232,192.4) node [anchor=north west][inner sep=0.75pt]    {$M$};
	\draw (201,102.4) node [anchor=north west][inner sep=0.75pt]    {$Y'$};
	\draw (202,202.4) node [anchor=north west][inner sep=0.75pt]    {$Y'$};
\end{tikzpicture}}

\newcommand{\figcobWo}{
\tikzset{every picture/.style={line width=0.75pt}} 
\begin{tikzpicture}[x=0.75pt,y=0.75pt,yscale=-1,xscale=1,scale=0.9]
\draw    (119.98,70) -- (119.98,236.57) ;
\draw    (119.98,70) -- (379.66,70) ;
\draw    (379.66,70) -- (379.66,236.57) ;
\draw [color={rgb, 255:red, 208; green, 2; blue, 27 }  ,draw opacity=1 ]   (119.98,125.52) .. controls (222.85,113.86) and (158.44,186.05) .. (119.98,147.73) ;
\draw  [color={rgb, 255:red, 208; green, 2; blue, 27 }  ,draw opacity=1 ] (140.01,196.71) .. controls (140.01,190.58) and (155.66,185.61) .. (174.97,185.61) .. controls (194.27,185.61) and (209.92,190.58) .. (209.92,196.71) .. controls (209.92,202.84) and (194.27,207.82) .. (174.97,207.82) .. controls (155.66,207.82) and (140.01,202.84) .. (140.01,196.71) -- cycle ;
\draw [color={rgb, 255:red, 208; green, 2; blue, 27 }  ,draw opacity=1 ]   (379.66,108.42) .. controls (288.27,91.21) and (249.49,115.23) .. (209.87,108.42) .. controls (170.25,101.61) and (204.23,57.02) .. (239.83,86.21) .. controls (275.44,115.4) and (307.32,60.03) .. (379.66,86.21) ;
\draw  [color={rgb, 255:red, 208; green, 2; blue, 27 }  ,draw opacity=1 ] (139.96,97.76) .. controls (139.96,88.56) and (146.67,81.1) .. (154.94,81.1) .. controls (163.21,81.1) and (169.92,88.56) .. (169.92,97.76) .. controls (169.92,106.96) and (163.21,114.42) .. (154.94,114.42) .. controls (146.67,114.42) and (139.96,106.96) .. (139.96,97.76) -- cycle ;
\draw    (119.98,236.57) -- (379.66,236.57) ;
\draw [color={rgb, 255:red, 208; green, 2; blue, 27 }  ,draw opacity=1 ]   (320,140) .. controls (272,127.5) and (145.95,153.84) .. (129.97,192.15) .. controls (113.99,230.47) and (209.37,223.25) .. (239.83,192.15) .. controls (270.29,161.06) and (308.34,183.85) .. (339.71,181.05) .. controls (371.07,178.25) and (373.5,156.5) .. (340,140) ;
\draw  [color={rgb, 255:red, 208; green, 2; blue, 27 }  ,draw opacity=1 ] (224.85,158.84) .. controls (224.85,154.24) and (240.5,150.51) .. (259.81,150.51) .. controls (279.11,150.51) and (294.76,154.24) .. (294.76,158.84) .. controls (294.76,163.44) and (279.11,167.17) .. (259.81,167.17) .. controls (240.5,167.17) and (224.85,163.44) .. (224.85,158.84) -- cycle ;
\draw [color={rgb, 255:red, 208; green, 2; blue, 27 }  ,draw opacity=1 ]   (379.66,192.15) .. controls (357.05,192.36) and (272.88,199.15) .. (294.81,209.53) .. controls (316.74,219.92) and (242.5,237.88) .. (379.66,214.36) ;
\draw  [color={rgb, 255:red, 208; green, 2; blue, 27 }  ,draw opacity=1 ] (259.81,208.81) .. controls (259.81,199.61) and (264.28,192.15) .. (269.8,192.15) .. controls (275.31,192.15) and (279.78,199.61) .. (279.78,208.81) .. controls (279.78,218.01) and (275.31,225.47) .. (269.8,225.47) .. controls (264.28,225.47) and (259.81,218.01) .. (259.81,208.81) -- cycle ;
\draw   (320,140) .. controls (320,134.48) and (324.48,130) .. (330,130) .. controls (335.52,130) and (340,134.48) .. (340,140) .. controls (340,145.52) and (335.52,150) .. (330,150) .. controls (324.48,150) and (320,145.52) .. (320,140) -- cycle ;
\draw  [color={rgb, 255:red, 208; green, 2; blue, 27 }  ,draw opacity=1 ][fill={rgb, 255:red, 208; green, 2; blue, 27 }  ,fill opacity=1 ] (337.75,140) .. controls (337.75,138.76) and (338.76,137.75) .. (340,137.75) .. controls (341.24,137.75) and (342.25,138.76) .. (342.25,140) .. controls (342.25,141.24) and (341.24,142.25) .. (340,142.25) .. controls (338.76,142.25) and (337.75,141.24) .. (337.75,140) -- cycle ;
\draw  [color={rgb, 255:red, 208; green, 2; blue, 27 }  ,draw opacity=1 ][fill={rgb, 255:red, 208; green, 2; blue, 27 }  ,fill opacity=1 ] (317.75,140) .. controls (317.75,138.76) and (318.76,137.75) .. (320,137.75) .. controls (321.24,137.75) and (322.25,138.76) .. (322.25,140) .. controls (322.25,141.24) and (321.24,142.25) .. (320,142.25) .. controls (318.76,142.25) and (317.75,141.24) .. (317.75,140) -- cycle ;
\draw (99,139.08) node [anchor=north west][inner sep=0.75pt]  [color={rgb, 255:red, 208; green, 2; blue, 27 }  ,opacity=1 ]  {$\mathbf{p}_{-}$};
\draw (380.63,140.19) node [anchor=north west][inner sep=0.75pt]  [color={rgb, 255:red, 208; green, 2; blue, 27 }  ,opacity=1 ]  {$\mathbf{p}_{+}$};
\draw (120.97,240.13) node [anchor=north west][inner sep=0.75pt]    {$Y_{-}$};
\draw (360.67,240.13) node [anchor=north west][inner sep=0.75pt]    {$Y_{+}$};
\draw (241,240.5) node [anchor=north west][inner sep=0.75pt]    {$W_{o}$};
\draw (181,112.4) node [anchor=north west][inner sep=0.75pt]  [color={rgb, 255:red, 208; green, 2; blue, 27 }  ,opacity=1 ]  {$\mathbf{a}$};
\draw (328,110.4) node [anchor=north west][inner sep=0.75pt]    {$B_{p}^{4}$};
\draw (301,122.4) node [anchor=north west][inner sep=0.75pt]  [color={rgb, 255:red, 208; green, 2; blue, 27 }  ,opacity=1 ]  {$q_{-}$};
\draw (339.75,145.4) node [anchor=north west][inner sep=0.75pt]  [color={rgb, 255:red, 208; green, 2; blue, 27 }  ,opacity=1 ]  {$q_{+}$};
\end{tikzpicture}}

\newcommand{\figcobW}{
\tikzset{every picture/.style={line width=0.75pt}} 
\begin{tikzpicture}[x=0.75pt,y=0.75pt,yscale=-1,xscale=1,scale=0.8]
\draw    (30,110) -- (30,260) ;
\draw    (30,110) -- (290,110) ;
\draw    (290,110) -- (290,260) ;
\draw [color={rgb, 255:red, 208; green, 2; blue, 27 }  ,draw opacity=1 ]   (30,160) .. controls (133,149.5) and (68.5,214.5) .. (30,180) ;
\draw  [color={rgb, 255:red, 208; green, 2; blue, 27 }  ,draw opacity=1 ] (50.05,224.1) .. controls (50.05,218.58) and (65.72,214.1) .. (85.05,214.1) .. controls (104.38,214.1) and (120.05,218.58) .. (120.05,224.1) .. controls (120.05,229.62) and (104.38,234.1) .. (85.05,234.1) .. controls (65.72,234.1) and (50.05,229.62) .. (50.05,224.1) -- cycle ;
\draw [color={rgb, 255:red, 208; green, 2; blue, 27 }  ,draw opacity=1 ]   (290,150) .. controls (198.5,134.5) and (159.67,156.13) .. (120,150) .. controls (80.33,143.87) and (114.35,103.71) .. (150,130) .. controls (185.65,156.29) and (217.57,106.42) .. (290,130) ;
\draw  [color={rgb, 255:red, 208; green, 2; blue, 27 }  ,draw opacity=1 ] (50,135) .. controls (50,126.72) and (56.72,120) .. (65,120) .. controls (73.28,120) and (80,126.72) .. (80,135) .. controls (80,143.28) and (73.28,150) .. (65,150) .. controls (56.72,150) and (50,143.28) .. (50,135) -- cycle ;
\draw    (30,260) -- (290,260) ;
\draw [color={rgb, 255:red, 208; green, 2; blue, 27 }  ,draw opacity=1 ]   (230,170) .. controls (186.5,154.5) and (56,185.5) .. (40,220) .. controls (24,254.5) and (119.5,248) .. (150,220) .. controls (180.5,192) and (218.6,212.52) .. (250,210) .. controls (281.4,207.48) and (277,187.5) .. (230,170) -- cycle ;
\draw  [color={rgb, 255:red, 208; green, 2; blue, 27 }  ,draw opacity=1 ] (135,190) .. controls (135,185.86) and (150.67,182.5) .. (170,182.5) .. controls (189.33,182.5) and (205,185.86) .. (205,190) .. controls (205,194.14) and (189.33,197.5) .. (170,197.5) .. controls (150.67,197.5) and (135,194.14) .. (135,190) -- cycle ;
\draw [color={rgb, 255:red, 208; green, 2; blue, 27 }  ,draw opacity=1 ]   (290,220) .. controls (267.37,220.19) and (183.09,226.3) .. (205.05,235.65) .. controls (227,245) and (152.68,261.17) .. (290,240) ;
\draw  [color={rgb, 255:red, 208; green, 2; blue, 27 }  ,draw opacity=1 ] (170,235) .. controls (170,226.72) and (174.48,220) .. (180,220) .. controls (185.52,220) and (190,226.72) .. (190,235) .. controls (190,243.28) and (185.52,250) .. (180,250) .. controls (174.48,250) and (170,243.28) .. (170,235) -- cycle ;
\draw (9,171.4) node [anchor=north west][inner sep=0.75pt]  [color={rgb, 255:red, 208; green, 2; blue, 27 }  ,opacity=1 ]  {$\mathbf{p}_{-}$};
\draw (298,172.4) node [anchor=north west][inner sep=0.75pt]  [color={rgb, 255:red, 208; green, 2; blue, 27 }  ,opacity=1 ]  {$\mathbf{p}_{+}$};
\draw (31,262.4) node [anchor=north west][inner sep=0.75pt]    {$Y_{-}$};
\draw (271,262.4) node [anchor=north west][inner sep=0.75pt]    {$Y_{+}$};
\draw (151,272.4) node [anchor=north west][inner sep=0.75pt]    {$W$};
\draw (141,153.4) node [anchor=north west][inner sep=0.75pt]  [color={rgb, 255:red, 208; green, 2; blue, 27 }  ,opacity=1 ]  {$\mathbf{a}$};
\end{tikzpicture}}

\newcommand{\figdiscandball}{
\tikzset{every picture/.style={line width=0.75pt}} 
\begin{tikzpicture}[x=0.75pt,y=0.75pt,yscale=-1,xscale=1,scale=0.8]
	\fill[gray!25] (107.645,134.21) ellipse (14.645 and 43.79);

	\fill[gray!25] (321.5,131) circle (20);
	\fill[gray!25] (391.5,131) circle (20);
	\draw    (320,191.38) -- (320,73.63) ;
	\draw [shift={(320,70.63)}, rotate = 90] [fill={rgb, 255:red, 0; green, 0; blue, 0 }  ][line width=0.08]  [draw opacity=0] (8.93,-4.29) -- (0,0) -- (8.93,4.29) -- cycle    ;
	\draw    (391.5,191) -- (391.5,74) ;
	\draw [shift={(391.5,71)}, rotate = 90] [fill={rgb, 255:red, 0; green, 0; blue, 0 }  ][line width=0.08]  [draw opacity=0] (8.93,-4.29) -- (0,0) -- (8.93,4.29) -- cycle    ;
	\draw  [color={rgb, 255:red, 0; green, 0; blue, 0 }  ,draw opacity=1 ][fill={rgb, 255:red, 0; green, 0; blue, 0 }  ,fill opacity=1 ] (318.5,131) .. controls (318.5,130.17) and (319.17,129.5) .. (320,129.5) .. controls (320.83,129.5) and (321.5,130.17) .. (321.5,131) .. controls (321.5,131.83) and (320.83,132.5) .. (320,132.5) .. controls (319.17,132.5) and (318.5,131.83) .. (318.5,131) -- cycle ;
	\draw  [color={rgb, 255:red, 0; green, 0; blue, 0 }  ,draw opacity=1 ][fill={rgb, 255:red, 0; green, 0; blue, 0 }  ,fill opacity=1 ] (390,132.5) .. controls (390,131.67) and (390.67,131) .. (391.5,131) .. controls (392.33,131) and (393,131.67) .. (393,132.5) .. controls (393,133.33) and (392.33,134) .. (391.5,134) .. controls (390.67,134) and (390,133.33) .. (390,132.5) -- cycle ;
	\draw  [dash pattern={on 0.84pt off 2.51pt}] (301.5,131) .. controls (301.5,119.95) and (310.45,111) .. (321.5,111) .. controls (332.55,111) and (341.5,119.95) .. (341.5,131) .. controls (341.5,142.05) and (332.55,151) .. (321.5,151) .. controls (310.45,151) and (301.5,142.05) .. (301.5,131) -- cycle ;
	\draw  [dash pattern={on 0.84pt off 2.51pt}] (371.5,131) .. controls (371.5,119.95) and (380.45,111) .. (391.5,111) .. controls (402.55,111) and (411.5,119.95) .. (411.5,131) .. controls (411.5,142.05) and (402.55,151) .. (391.5,151) .. controls (380.45,151) and (371.5,142.05) .. (371.5,131) -- cycle ;
	\draw (75,179.49) -- (197.59,88.30);
	\draw[shift={(200,86.51)},rotate=143.36,fill=black,draw=none]
	(8.93,-4.29) -- (0,0) -- (8.93,4.29) -- cycle;

	\draw (131.29,131.92) -- (67.48,88.64);
	\draw[shift={(65,86.96)},rotate=34.15,fill=black,draw=none]
	(8.93,-4.29) -- (0,0) -- (8.93,4.29) -- cycle;

	\draw (138.79,137.92) -- (200,177.37);
	\draw  [dash pattern={on 0.84pt off 2.51pt}] (65,130) .. controls (65,91.34) and (96.34,60) .. (135,60) .. controls (173.66,60) and (205,91.34) .. (205,130) .. controls (205,168.66) and (173.66,200) .. (135,200) .. controls (96.34,200) and (65,168.66) .. (65,130) -- cycle ;
	\draw   (93,134.21) .. controls (93,110.03) and (99.56,90.42) .. (107.64,90.42) .. controls (115.73,90.42) and (122.29,110.03) .. (122.29,134.21) .. controls (122.29,158.39) and (115.73,178) .. (107.64,178) .. controls (99.56,178) and (93,158.39) .. (93,134.21) -- cycle ;
	\draw  [color={rgb, 255:red, 0; green, 0; blue, 0 }  ,draw opacity=1 ][fill={rgb, 255:red, 0; green, 0; blue, 0 }  ,fill opacity=1 ] (106.5,154.5) .. controls (106.5,153.67) and (107.17,153) .. (108,153) .. controls (108.83,153) and (109.5,153.67) .. (109.5,154.5) .. controls (109.5,155.33) and (108.83,156) .. (108,156) .. controls (107.17,156) and (106.5,155.33) .. (106.5,154.5) -- cycle ;
	\draw  [color={rgb, 255:red, 0; green, 0; blue, 0 }  ,draw opacity=1 ][fill={rgb, 255:red, 0; green, 0; blue, 0 }  ,fill opacity=1 ] (105.5,115) .. controls (105.5,114.17) and (106.17,113.5) .. (107,113.5) .. controls (107.83,113.5) and (108.5,114.17) .. (108.5,115) .. controls (108.5,115.83) and (107.83,116.5) .. (107,116.5) .. controls (106.17,116.5) and (105.5,115.83) .. (105.5,115) -- cycle ;
	\draw (303.5,120.4) node [anchor=north west][inner sep=0.75pt]    {$r_{1}$};
	\draw (394.5,123.4) node [anchor=north west][inner sep=0.75pt]    {$r_{2}$};
	\draw (281.5,103.4) node [anchor=north west][inner sep=0.75pt]    {$B'_{1}$};
	\draw (210.5,120.4) node [anchor=north west][inner sep=0.75pt]    {$L_+$};
	\draw (410.5,103.4) node [anchor=north west][inner sep=0.75pt]    {$B'_{2}$};
	\draw (100,130.4) node [anchor=north west][inner sep=0.75pt]    {$\mathsf{D}$};
	\draw (116.5,90) node [anchor=north west][inner sep=0.75pt]    {$\gamma$};
\end{tikzpicture}}

\newcommand{\figWplusminus}{
	\tikzset{every picture/.style={line width=0.75pt}} 
	\begin{tikzpicture}[x=0.75pt,y=0.75pt,yscale=-1,xscale=1,scale=0.8]
		\draw    (110,120) -- (110,270) ;
		\draw    (110,270) -- (480,270) ;
		\draw    (150,120) -- (150,190) ;
		\draw   (340,220) .. controls (340,214.48) and (344.48,210) .. (350,210) .. controls (355.52,210) and (360,214.48) .. (360,220) .. controls (360,225.52) and (355.52,230) .. (350,230) .. controls (344.48,230) and (340,225.52) .. (340,220) -- cycle ;
		\draw    (150,190) -- (340,190) ;
		\draw    (150,250) -- (340,250) ;
		\draw    (150,250) -- (150,270) ;
		\draw    (480,120) -- (480,200) ;
		\draw    (480,240) -- (480,270) ;
		\draw    (340,190) .. controls (429.33,190.33) and (430.57,249.42) .. (340,250) ;
		\draw    (340,200) .. controls (300.57,199.92) and (300.57,239.42) .. (340,240) ;
		\draw    (340,200) -- (480,200) ;
		\draw    (340,240) -- (480,240) ;
		\draw    (110,120) -- (480,120) ;
		\draw    (520,120) -- (520,270) ;
		\draw    (480,270) -- (520,270) ;
		\draw    (480,120) -- (520,120) ;
		\draw [color={rgb, 255:red, 208; green, 2; blue, 27 }  ,draw opacity=1 ]   (110,220) -- (340,220) ;
		\draw [color={rgb, 255:red, 208; green, 2; blue, 27 }  ,draw opacity=1 ]   (360,220) -- (520,220) ;
		\draw [color={rgb, 255:red, 208; green, 2; blue, 27 }  ,draw opacity=1 ]   (110,130) .. controls (158.07,123.42) and (141.57,153.42) .. (110,150) ;
		\draw  [color={rgb, 255:red, 208; green, 2; blue, 27 }  ,draw opacity=1 ] (325,150) .. controls (325,144.48) and (340.67,140) .. (360,140) .. controls (379.33,140) and (395,144.48) .. (395,150) .. controls (395,155.52) and (379.33,160) .. (360,160) .. controls (340.67,160) and (325,155.52) .. (325,150) -- cycle ;
		\draw [color={rgb, 255:red, 208; green, 2; blue, 27 }  ,draw opacity=1 ]   (110,180) .. controls (204.57,126.42) and (261.07,133.92) .. (310,170) .. controls (358.93,206.08) and (447.57,126.42) .. (520,150) ;
		\draw  [color={rgb, 255:red, 208; green, 2; blue, 27 }  ,draw opacity=1 ] (200,165) .. controls (200,156.72) and (206.72,150) .. (215,150) .. controls (223.28,150) and (230,156.72) .. (230,165) .. controls (230,173.28) and (223.28,180) .. (215,180) .. controls (206.72,180) and (200,173.28) .. (200,165) -- cycle ;
		\draw (121,201.4) node [anchor=north west][inner sep=0.75pt]  [color={rgb, 255:red, 208; green, 2; blue, 27 }  ,opacity=1 ]  {$a_{*}$};
		\draw (295,95.4) node [anchor=north west][inner sep=0.75pt]    {$W_{o}$};
		\draw (85.5,239.4) node [anchor=north west][inner sep=0.75pt]    {$Y_{-}$};
		\draw (84,139.9) node [anchor=north west][inner sep=0.75pt]  [color={rgb, 255:red, 208; green, 2; blue, 27 }  ,opacity=1 ]  {$\mathbf{p}_{-}$};
		\draw (522,138.4) node [anchor=north west][inner sep=0.75pt]  [color={rgb, 255:red, 208; green, 2; blue, 27 }  ,opacity=1 ]  {$\mathbf{p}_{+}$};
		\draw (142.5,273.9) node [anchor=north west][inner sep=0.75pt]    {$Y_{-+}$};
		\draw (466.5,273.4) node [anchor=north west][inner sep=0.75pt]    {$Y_{+-}$};
		\draw (522,239.4) node [anchor=north west][inner sep=0.75pt]    {$Y_{+}$};
\end{tikzpicture}}

\newcommand{\figVVdag}{
	\tikzset{every picture/.style={line width=0.75pt}} 
	\begin{tikzpicture}[x=0.75pt,y=0.75pt,yscale=-1,xscale=1,scale=0.8]
		\begin{scope}
		\clip
		(10,50) -- (130,50)
		.. controls (130.96,80.44) and (179.46,80.44) .. (180,50)
		-- (300,50) -- (300,250) -- (180,250)
		.. controls (179.04,219.56) and (130.54,219.56) .. (130,250)
		-- (10,250) -- cycle;
		\fill[gray!18,even odd rule]
		(10,50) rectangle (300,250)
		(250,65) -- (60,65)
		.. controls (95.89,174.92) and (80.89,153.42) .. (60,235)
		-- (250,235)
		.. controls (205.89,153.42) and (238.39,102.42) .. (250,65)
		-- cycle;
		\end{scope}
		\fill[gray!18,even odd rule]
		(340,50) rectangle (630,250)
		(390,150)
		.. controls (390,105.82) and (432.73,70) .. (485.44,70)
		.. controls (538.16,70) and (580.89,105.82) .. (580.89,150)
		.. controls (580.89,194.18) and (538.16,230) .. (485.44,230)
		.. controls (432.73,230) and (390,194.18) .. (390,150)
		-- cycle;
		\node[anchor=north west,inner sep=0.75pt]
		at (256,102.4) {$V$};
		\node[anchor=north west,inner sep=0.75pt]
		at (368,98.4) {$V^{\dagger}$};
\node[anchor=north west,inner sep=0.75pt]
		at (125,134.4) {$\mathsf{T}_1 \sqcup \mathsf{T}_2$};
\node[anchor=north west,inner sep=0.75pt]
		at (450,90.4) {$\mathsf{S}_1 \sqcup \mathsf{S}_2$};
		\node[anchor=south] at (70,45)  {$Y$};
		\node[anchor=south] at (240,45) {$Z$};
		\node[anchor=north] at (70,255)  {$Y$};
		\node[anchor=north] at (240,255) {$Z$};
		\draw    (130,250) -- (10,250) ;
		\draw    (130,250) .. controls (130.54,219.56) and (179.04,219.56) .. (180,250) ;
		\draw    (300,250) -- (180,250) ;
		\draw    (10,150) -- (10,250) ;
		\draw [color={rgb, 255:red, 208; green, 2; blue, 27 }  ,draw opacity=1 ]   (30,150) -- (30,250) ;
		\draw [color={rgb, 255:red, 208; green, 2; blue, 27 }  ,draw opacity=1 ]   (100,250) .. controls (100.04,200.06) and (209.04,200.06) .. (210,250) ;
		\draw    (300,150) -- (300,250) ;
		\draw [color={rgb, 255:red, 208; green, 2; blue, 27 }  ,draw opacity=1 ]   (280,150) -- (280,250) ;
		\draw [color={rgb, 255:red, 208; green, 2; blue, 27 }  ,draw opacity=1 ]   (70,250) .. controls (71.04,160.06) and (240.54,160.56) .. (240,250) ;
		\draw    (180,50) -- (300,50) ;
		\draw    (180,50) .. controls (179.46,80.44) and (130.96,80.44) .. (130,50) ;
		\draw    (10,50) -- (130,50) ;
		\draw    (300,150) -- (300,50) ;
		\draw [color={rgb, 255:red, 208; green, 2; blue, 27 }  ,draw opacity=1 ]   (280,150) -- (280,50) ;
		\draw [color={rgb, 255:red, 208; green, 2; blue, 27 }  ,draw opacity=1 ]   (210,50) .. controls (209.96,99.94) and (100.96,99.94) .. (100,50) ;
		\draw    (10,150) -- (10,50) ;
		\draw [color={rgb, 255:red, 208; green, 2; blue, 27 }  ,draw opacity=1 ]   (30,150) -- (30,50) ;
		\draw [color={rgb, 255:red, 208; green, 2; blue, 27 }  ,draw opacity=1 ]   (240,50) .. controls (238.96,139.94) and (69.46,139.44) .. (70,50) ;
		\draw    (630,50) -- (340,50) ;
		\draw    (460,150) .. controls (460.54,119.56) and (509.04,119.56) .. (510,150) ;
		\draw    (340,50) -- (340,150) ;
		\draw [color={rgb, 255:red, 208; green, 2; blue, 27 }  ,draw opacity=1 ]   (360,50) -- (360,150) ;
		\draw [color={rgb, 255:red, 208; green, 2; blue, 27 }  ,draw opacity=1 ]   (430,150) .. controls (430.04,100.06) and (539.04,100.06) .. (540,150) ;
		\draw    (630,50) -- (630,150) ;
		\draw [color={rgb, 255:red, 208; green, 2; blue, 27 }  ,draw opacity=1 ]   (610,50) -- (610,150) ;
		\draw [color={rgb, 255:red, 208; green, 2; blue, 27 }  ,draw opacity=1 ]   (410,150) .. controls (411.04,60.06) and (560.54,60.56) .. (560,150) ;
		\draw    (340,250) -- (630,250) ;
		\draw    (510,150) .. controls (509.46,180.44) and (460.96,180.44) .. (460,150) ;
		\draw    (630,250) -- (630,150) ;
		\draw [color={rgb, 255:red, 208; green, 2; blue, 27 }  ,draw opacity=1 ]   (610,250) -- (610,150) ;
		\draw [color={rgb, 255:red, 208; green, 2; blue, 27 }  ,draw opacity=1 ]   (540,150) .. controls (539.96,199.94) and (430.96,199.94) .. (430,150) ;
		\draw    (340,250) -- (340,150) ;
		\draw [color={rgb, 255:red, 208; green, 2; blue, 27 }  ,draw opacity=1 ]   (360,250) -- (360,150) ;
		\draw [color={rgb, 255:red, 208; green, 2; blue, 27 }  ,draw opacity=1 ]   (560,150) .. controls (558.96,239.94) and (409.46,239.44) .. (410,150) ;
		\draw  [dash pattern={on 0.84pt off 2.51pt}] (390,150) .. controls (390,105.82) and (432.73,70) .. (485.44,70) .. controls (538.16,70) and (580.89,105.82) .. (580.89,150) .. controls (580.89,194.18) and (538.16,230) .. (485.44,230) .. controls (432.73,230) and (390,194.18) .. (390,150) -- cycle ;
		\draw  [dash pattern={on 0.84pt off 2.51pt}]  (250,65) -- (60,65) .. controls (95.89,174.92) and (80.89,153.42) .. (60,235) ;
		\draw  [dash pattern={on 0.84pt off 2.51pt}]  (60,235) -- (250,235) .. controls (205.89,153.42) and (238.39,102.42) .. (250,65) ;
		\node[anchor=south] at (485,45)  {$Y\#_2 Z$};
		\node[anchor=north] at (485,255) {$Y\#_2 Z$};
\end{tikzpicture}}

\newcommand{\figcomposites}{
	\tikzset{every picture/.style={line width=0.75pt}} 
	\begin{tikzpicture}[x=0.75pt,y=0.75pt,yscale=-1,xscale=1,scale=0.85]
		\draw    (130,250) -- (10,250) ;
		\draw    (130,250) .. controls (130.54,219.56) and (179.04,219.56) .. (180,250) ;
		\draw    (300,250) -- (180,250) ;
		\draw    (10,150) -- (10,250) ;
		\draw [color={rgb, 255:red, 208; green, 2; blue, 27 }  ,draw opacity=1 ]   (30,150) -- (30,250) ;
		\draw [color={rgb, 255:red, 208; green, 2; blue, 27 }  ,draw opacity=1 ]   (100,250) .. controls (100.04,200.06) and (209.04,200.06) .. (210,250) ;
		\draw    (300,150) -- (300,250) ;
		\draw [color={rgb, 255:red, 208; green, 2; blue, 27 }  ,draw opacity=1 ]   (280,150) -- (280,250) ;
		\draw [color={rgb, 255:red, 208; green, 2; blue, 27 }  ,draw opacity=1 ]   (70,250) .. controls (71.04,160.06) and (240.54,160.56) .. (240,250) ;
		\draw  [color={rgb, 255:red, 0; green, 0; blue, 0 }  ,draw opacity=1 ][fill={rgb, 255:red, 0; green, 0; blue, 0 }  ,fill opacity=1 ] (67,250) .. controls (67,248.34) and (68.34,247) .. (70,247) .. controls (71.66,247) and (73,248.34) .. (73,250) .. controls (73,251.66) and (71.66,253) .. (70,253) .. controls (68.34,253) and (67,251.66) .. (67,250) -- cycle ;
		\draw  [color={rgb, 255:red, 0; green, 0; blue, 0 }  ,draw opacity=1 ][fill={rgb, 255:red, 0; green, 0; blue, 0 }  ,fill opacity=1 ] (97,250) .. controls (97,248.34) and (98.34,247) .. (100,247) .. controls (101.66,247) and (103,248.34) .. (103,250) .. controls (103,251.66) and (101.66,253) .. (100,253) .. controls (98.34,253) and (97,251.66) .. (97,250) -- cycle ;
		\draw  [color={rgb, 255:red, 0; green, 0; blue, 0 }  ,draw opacity=1 ][fill={rgb, 255:red, 0; green, 0; blue, 0 }  ,fill opacity=1 ] (207,250) .. controls (207,248.34) and (208.34,247) .. (210,247) .. controls (211.66,247) and (213,248.34) .. (213,250) .. controls (213,251.66) and (211.66,253) .. (210,253) .. controls (208.34,253) and (207,251.66) .. (207,250) -- cycle ;
		\draw  [color={rgb, 255:red, 0; green, 0; blue, 0 }  ,draw opacity=1 ][fill={rgb, 255:red, 0; green, 0; blue, 0 }  ,fill opacity=1 ] (237,250) .. controls (237,248.34) and (238.34,247) .. (240,247) .. controls (241.66,247) and (243,248.34) .. (243,250) .. controls (243,251.66) and (241.66,253) .. (240,253) .. controls (238.34,253) and (237,251.66) .. (237,250) -- cycle ;
		\draw    (180,50) -- (300,50) ;
		\draw    (180,50) .. controls (179.46,80.44) and (130.96,80.44) .. (130,50) ;
		\draw    (10,50) -- (130,50) ;
		\draw    (300,150) -- (300,50) ;
		\draw [color={rgb, 255:red, 208; green, 2; blue, 27 }  ,draw opacity=1 ]   (280,150) -- (280,50) ;
		\draw [color={rgb, 255:red, 208; green, 2; blue, 27 }  ,draw opacity=1 ]   (210,50) .. controls (209.96,99.94) and (100.96,99.94) .. (100,50) ;
		\draw    (10,150) -- (10,50) ;
		\draw [color={rgb, 255:red, 208; green, 2; blue, 27 }  ,draw opacity=1 ]   (30,150) -- (30,50) ;
		\draw [color={rgb, 255:red, 208; green, 2; blue, 27 }  ,draw opacity=1 ]   (240,50) .. controls (238.96,139.94) and (69.46,139.44) .. (70,50) ;
		\draw  [color={rgb, 255:red, 0; green, 0; blue, 0 }  ,draw opacity=1 ][fill={rgb, 255:red, 0; green, 0; blue, 0 }  ,fill opacity=1 ] (243,50) .. controls (243,51.66) and (241.66,53) .. (240,53) .. controls (238.34,53) and (237,51.66) .. (237,50) .. controls (237,48.34) and (238.34,47) .. (240,47) .. controls (241.66,47) and (243,48.34) .. (243,50) -- cycle ;
		\draw  [color={rgb, 255:red, 0; green, 0; blue, 0 }  ,draw opacity=1 ][fill={rgb, 255:red, 0; green, 0; blue, 0 }  ,fill opacity=1 ] (213,50) .. controls (213,51.66) and (211.66,53) .. (210,53) .. controls (208.34,53) and (207,51.66) .. (207,50) .. controls (207,48.34) and (208.34,47) .. (210,47) .. controls (211.66,47) and (213,48.34) .. (213,50) -- cycle ;
		\draw  [color={rgb, 255:red, 0; green, 0; blue, 0 }  ,draw opacity=1 ][fill={rgb, 255:red, 0; green, 0; blue, 0 }  ,fill opacity=1 ] (103,50) .. controls (103,51.66) and (101.66,53) .. (100,53) .. controls (98.34,53) and (97,51.66) .. (97,50) .. controls (97,48.34) and (98.34,47) .. (100,47) .. controls (101.66,47) and (103,48.34) .. (103,50) -- cycle ;
		\draw  [color={rgb, 255:red, 0; green, 0; blue, 0 }  ,draw opacity=1 ][fill={rgb, 255:red, 0; green, 0; blue, 0 }  ,fill opacity=1 ] (73,50) .. controls (73,51.66) and (71.66,53) .. (70,53) .. controls (68.34,53) and (67,51.66) .. (67,50) .. controls (67,48.34) and (68.34,47) .. (70,47) .. controls (71.66,47) and (73,48.34) .. (73,50) -- cycle ;
		\draw    (630,50) -- (340,50) ;
		\draw    (460,150) -- (340,150) ;
		\draw    (460,150) .. controls (460.54,119.56) and (509.04,119.56) .. (510,150) ;
		\draw    (630,150) -- (510,150) ;
		\draw    (340,50) -- (340,150) ;
		\draw [color={rgb, 255:red, 208; green, 2; blue, 27 }  ,draw opacity=1 ]   (360,50) -- (360,150) ;
		\draw [color={rgb, 255:red, 208; green, 2; blue, 27 }  ,draw opacity=1 ]   (430,150) .. controls (430.04,100.06) and (539.04,100.06) .. (540,150) ;
		\draw    (630,50) -- (630,150) ;
		\draw [color={rgb, 255:red, 208; green, 2; blue, 27 }  ,draw opacity=1 ]   (610,50) -- (610,150) ;
		\draw [color={rgb, 255:red, 208; green, 2; blue, 27 }  ,draw opacity=1 ]   (400,150) .. controls (401.04,60.06) and (570.54,60.56) .. (570,150) ;
		\draw  [color={rgb, 255:red, 0; green, 0; blue, 0 }  ,draw opacity=1 ][fill={rgb, 255:red, 0; green, 0; blue, 0 }  ,fill opacity=1 ] (397,150) .. controls (397,148.34) and (398.34,147) .. (400,147) .. controls (401.66,147) and (403,148.34) .. (403,150) .. controls (403,151.66) and (401.66,153) .. (400,153) .. controls (398.34,153) and (397,151.66) .. (397,150) -- cycle ;
		\draw  [color={rgb, 255:red, 0; green, 0; blue, 0 }  ,draw opacity=1 ][fill={rgb, 255:red, 0; green, 0; blue, 0 }  ,fill opacity=1 ] (427,150) .. controls (427,148.34) and (428.34,147) .. (430,147) .. controls (431.66,147) and (433,148.34) .. (433,150) .. controls (433,151.66) and (431.66,153) .. (430,153) .. controls (428.34,153) and (427,151.66) .. (427,150) -- cycle ;
		\draw  [color={rgb, 255:red, 0; green, 0; blue, 0 }  ,draw opacity=1 ][fill={rgb, 255:red, 0; green, 0; blue, 0 }  ,fill opacity=1 ] (537,150) .. controls (537,148.34) and (538.34,147) .. (540,147) .. controls (541.66,147) and (543,148.34) .. (543,150) .. controls (543,151.66) and (541.66,153) .. (540,153) .. controls (538.34,153) and (537,151.66) .. (537,150) -- cycle ;
		\draw  [color={rgb, 255:red, 0; green, 0; blue, 0 }  ,draw opacity=1 ][fill={rgb, 255:red, 0; green, 0; blue, 0 }  ,fill opacity=1 ] (567,150) .. controls (567,148.34) and (568.34,147) .. (570,147) .. controls (571.66,147) and (573,148.34) .. (573,150) .. controls (573,151.66) and (571.66,153) .. (570,153) .. controls (568.34,153) and (567,151.66) .. (567,150) -- cycle ;
		\draw    (340,250) -- (630,250) ;
		\draw    (510,150) .. controls (509.46,180.44) and (460.96,180.44) .. (460,150) ;
		\draw    (630,250) -- (630,150) ;
		\draw [color={rgb, 255:red, 208; green, 2; blue, 27 }  ,draw opacity=1 ]   (610,250) -- (610,150) ;
		\draw [color={rgb, 255:red, 208; green, 2; blue, 27 }  ,draw opacity=1 ]   (540,150) .. controls (539.96,199.94) and (430.96,199.94) .. (430,150) ;
		\draw    (340,250) -- (340,150) ;
		\draw [color={rgb, 255:red, 208; green, 2; blue, 27 }  ,draw opacity=1 ]   (360,250) -- (360,150) ;
		\draw [color={rgb, 255:red, 208; green, 2; blue, 27 }  ,draw opacity=1 ]   (570,150) .. controls (568.96,239.94) and (399.46,239.44) .. (400,150) ;
		\draw (31,132.4) node [anchor=north west][inner sep=0.75pt]    {$\textcolor[rgb]{0.82,0.01,0.11}{I\times }\mathbf{\textcolor[rgb]{0.82,0.01,0.11}{p}}$};
		\draw (145,140.9) node [anchor=north west][inner sep=0.75pt]    {$X$};
		\draw (568.5,61.4) node [anchor=north west][inner sep=0.75pt]    {$X^{\dagger }$};
		\draw (35,253.4) node [anchor=north west][inner sep=0.75pt]    {$Y$};
		\draw (254,252.9) node [anchor=north west][inner sep=0.75pt]    {$Z$};
		\draw (467,255.9) node [anchor=north west][inner sep=0.75pt]    {$Y\#_{2} Z$};
		\draw (241,132.4) node [anchor=north west][inner sep=0.75pt]    {$\textcolor[rgb]{0.82,0.01,0.11}{I\times }\textcolor[rgb]{0.82,0.01,0.11}{\mathbf{q}}$};
		\draw (571,212.4) node [anchor=north west][inner sep=0.75pt]    {$\textcolor[rgb]{0.82,0.01,0.11}{I\times }\textcolor[rgb]{0.82,0.01,0.11}{\mathbf{q}}$};
		\draw (364,212.4) node [anchor=north west][inner sep=0.75pt]    {$\textcolor[rgb]{0.82,0.01,0.11}{I\times }\mathbf{\textcolor[rgb]{0.82,0.01,0.11}{p}}$};
		\draw (36,27.9) node [anchor=north west][inner sep=0.75pt]    {$Y$};
		\draw (253.5,27.9) node [anchor=north west][inner sep=0.75pt]    {$Z$};
		\draw (461.5,25.4) node [anchor=north west][inner sep=0.75pt]    {$Y\#_{2} Z$};
		\draw (181,167.9) node [anchor=north west][inner sep=0.75pt]  [color={rgb, 255:red, 208; green, 2; blue, 27 }  ,opacity=1 ]  {$b_{1}$};
		\draw (174,197.9) node [anchor=north west][inner sep=0.75pt]  [color={rgb, 255:red, 208; green, 2; blue, 27 }  ,opacity=1 ]  {$b_{2}$};
		\draw (152.5,86.9) node [anchor=north west][inner sep=0.75pt]  [color={rgb, 255:red, 208; green, 2; blue, 27 }  ,opacity=1 ]  {$\bar{b}_{2}$};
		\draw (172,116.4) node [anchor=north west][inner sep=0.75pt]  [color={rgb, 255:red, 208; green, 2; blue, 27 }  ,opacity=1 ]  {$\bar{b}_{1}$};
		\draw (417,78.4) node [anchor=north west][inner sep=0.75pt]  [color={rgb, 255:red, 208; green, 2; blue, 27 }  ,opacity=1 ]  {$b_{1}$};
		\draw (420.5,200.9) node [anchor=north west][inner sep=0.75pt]  [color={rgb, 255:red, 208; green, 2; blue, 27 }  ,opacity=1 ]  {$\bar{b}_{1}$};
		\draw (455,98.9) node [anchor=north west][inner sep=0.75pt]  [color={rgb, 255:red, 208; green, 2; blue, 27 }  ,opacity=1 ]  {$b_{2}$};
		\draw (455,183.9) node [anchor=north west][inner sep=0.75pt]  [color={rgb, 255:red, 208; green, 2; blue, 27 }  ,opacity=1 ]  {$\bar{b}_{2}$};
\end{tikzpicture}}

\newcommand{\figWWdag}{
\tikzset{every picture/.style={line width=0.75pt}} 
	\begin{tikzpicture}[x=0.75pt,y=0.75pt,yscale=-1,xscale=1,scale=0.8]
		\draw    (300,80) -- (10,80) ;
		\draw    (130,250) -- (10,250) ;
		\draw    (130,250) .. controls (129.98,159.24) and (178.98,158.74) .. (180,250) ;
		\draw    (300,250) -- (180,250) ;
		\draw    (10,80) -- (10,250) ;
		\draw [color={rgb, 255:red, 208; green, 2; blue, 27 }  ,draw opacity=1 ]   (30,80) -- (30,250) ;
		\draw [color={rgb, 255:red, 208; green, 2; blue, 27 }  ,draw opacity=1 ]   (100,250) .. controls (101.98,129.74) and (209.98,130.24) .. (210,250) ;
		\draw    (300,80) -- (300,250) ;
		\draw [color={rgb, 255:red, 208; green, 2; blue, 27 }  ,draw opacity=1 ]   (280,80) -- (280,250) ;
		\draw [color={rgb, 255:red, 208; green, 2; blue, 27 }  ,draw opacity=1 ]   (70,250) .. controls (70.48,80.74) and (239.98,80.74) .. (240,250) ;
		\draw    (360,250) -- (650,250) ;
		\draw    (530,80) .. controls (530.02,170.76) and (481.02,171.26) .. (480,80) ;
		\draw    (650,250) -- (650,80) ;
		\draw [color={rgb, 255:red, 208; green, 2; blue, 27 }  ,draw opacity=1 ]   (630,250) -- (630,80) ;
		\draw [color={rgb, 255:red, 208; green, 2; blue, 27 }  ,draw opacity=1 ]   (560,80) .. controls (558.02,200.26) and (450.02,199.76) .. (450,80) ;
		\draw    (360,250) -- (360,80) ;
		\draw [color={rgb, 255:red, 208; green, 2; blue, 27 }  ,draw opacity=1 ]   (380,250) -- (380,80) ;
		\draw [color={rgb, 255:red, 208; green, 2; blue, 27 }  ,draw opacity=1 ]   (590,80) .. controls (589.52,249.26) and (420.02,249.26) .. (420,80) ;
		\draw    (480,80) -- (360,80) ;
		\draw    (650,80) -- (530,80) ;
		\draw  [color={rgb, 255:red, 0; green, 0; blue, 0 }  ,draw opacity=1 ][fill={rgb, 255:red, 0; green, 0; blue, 0 }  ,fill opacity=1 ] (67,250) .. controls (67,248.34) and (68.34,247) .. (70,247) .. controls (71.66,247) and (73,248.34) .. (73,250) .. controls (73,251.66) and (71.66,253) .. (70,253) .. controls (68.34,253) and (67,251.66) .. (67,250) -- cycle ;
		\draw  [color={rgb, 255:red, 0; green, 0; blue, 0 }  ,draw opacity=1 ][fill={rgb, 255:red, 0; green, 0; blue, 0 }  ,fill opacity=1 ] (97,250) .. controls (97,248.34) and (98.34,247) .. (100,247) .. controls (101.66,247) and (103,248.34) .. (103,250) .. controls (103,251.66) and (101.66,253) .. (100,253) .. controls (98.34,253) and (97,251.66) .. (97,250) -- cycle ;
		\draw  [color={rgb, 255:red, 0; green, 0; blue, 0 }  ,draw opacity=1 ][fill={rgb, 255:red, 0; green, 0; blue, 0 }  ,fill opacity=1 ] (207,250) .. controls (207,248.34) and (208.34,247) .. (210,247) .. controls (211.66,247) and (213,248.34) .. (213,250) .. controls (213,251.66) and (211.66,253) .. (210,253) .. controls (208.34,253) and (207,251.66) .. (207,250) -- cycle ;
		\draw  [color={rgb, 255:red, 0; green, 0; blue, 0 }  ,draw opacity=1 ][fill={rgb, 255:red, 0; green, 0; blue, 0 }  ,fill opacity=1 ] (237,250) .. controls (237,248.34) and (238.34,247) .. (240,247) .. controls (241.66,247) and (243,248.34) .. (243,250) .. controls (243,251.66) and (241.66,253) .. (240,253) .. controls (238.34,253) and (237,251.66) .. (237,250) -- cycle ;
		\draw  [color={rgb, 255:red, 0; green, 0; blue, 0 }  ,draw opacity=1 ][fill={rgb, 255:red, 0; green, 0; blue, 0 }  ,fill opacity=1 ] (417,80) .. controls (417,78.34) and (418.34,77) .. (420,77) .. controls (421.66,77) and (423,78.34) .. (423,80) .. controls (423,81.66) and (421.66,83) .. (420,83) .. controls (418.34,83) and (417,81.66) .. (417,80) -- cycle ;
		\draw  [color={rgb, 255:red, 0; green, 0; blue, 0 }  ,draw opacity=1 ][fill={rgb, 255:red, 0; green, 0; blue, 0 }  ,fill opacity=1 ] (447,80) .. controls (447,78.34) and (448.34,77) .. (450,77) .. controls (451.66,77) and (453,78.34) .. (453,80) .. controls (453,81.66) and (451.66,83) .. (450,83) .. controls (448.34,83) and (447,81.66) .. (447,80) -- cycle ;
		\draw  [color={rgb, 255:red, 0; green, 0; blue, 0 }  ,draw opacity=1 ][fill={rgb, 255:red, 0; green, 0; blue, 0 }  ,fill opacity=1 ] (557,80) .. controls (557,78.34) and (558.34,77) .. (560,77) .. controls (561.66,77) and (563,78.34) .. (563,80) .. controls (563,81.66) and (561.66,83) .. (560,83) .. controls (558.34,83) and (557,81.66) .. (557,80) -- cycle ;
		\draw  [color={rgb, 255:red, 0; green, 0; blue, 0 }  ,draw opacity=1 ][fill={rgb, 255:red, 0; green, 0; blue, 0 }  ,fill opacity=1 ] (587,80) .. controls (587,78.34) and (588.34,77) .. (590,77) .. controls (591.66,77) and (593,78.34) .. (593,80) .. controls (593,81.66) and (591.66,83) .. (590,83) .. controls (588.34,83) and (587,81.66) .. (587,80) -- cycle ;
		\draw  [dash pattern={on 0.84pt off 2.51pt}]  (330,210) -- (330,122) ;
		\draw [shift={(330,120)}, rotate = 90] [color={rgb, 255:red, 0; green, 0; blue, 0 }  ][line width=0.75]    (10.93,-3.29) .. controls (6.95,-1.4) and (3.31,-0.3) .. (0,0) .. controls (3.31,0.3) and (6.95,1.4) .. (10.93,3.29)   ;
		\draw (63.5,250.4) node [anchor=north west][inner sep=0.75pt]    {$x_{1}$};
		\draw (31,102.4) node [anchor=north west][inner sep=0.75pt]    {$\textcolor[rgb]{0.82,0.01,0.11}{I\times \mathbf{p}}$};
		\draw (92.5,250.4) node [anchor=north west][inner sep=0.75pt]    {$x_{2}$};
		\draw (205.5,249.9) node [anchor=north west][inner sep=0.75pt]    {$z_{2}$};
		\draw (234.5,249.9) node [anchor=north west][inner sep=0.75pt]    {$z_{1}$};
		\draw (411,57.9) node [anchor=north west][inner sep=0.75pt]    {$x_{1}$};
		\draw (581,138.9) node [anchor=north west][inner sep=0.75pt]    {\textcolor[rgb]{0.82,0.01,0.11}{$b_{1}$}};
		\draw (545,138.9) node [anchor=north west][inner sep=0.75pt]    {\textcolor[rgb]{0.82,0.01,0.11}{$b_{2}$}};
		\draw (441,57.9) node [anchor=north west][inner sep=0.75pt]    {$x_{2}$};
		\draw (552,56.4) node [anchor=north west][inner sep=0.75pt]    {$z_{2}$};
		\draw (581,56.4) node [anchor=north west][inner sep=0.75pt]    {$z_{1}$};
	\draw (100,165.9) node [anchor=north west][inner sep=0.75pt]    {\textcolor[rgb]{0.82,0.01,0.11}{$b_{2}$}};
		\draw (65,165.9) node [anchor=north west][inner sep=0.75pt]    {\textcolor[rgb]{0.82,0.01,0.11}{$b_{1}$}};
		\draw (241,172.4) node [anchor=north west][inner sep=0.75pt]    {$W$};
		\draw (391,171.4) node [anchor=north west][inner sep=0.75pt]    {$W^{\dagger }$};
		\draw (27,262.4) node [anchor=north west][inner sep=0.75pt]    {$Y$};
		\draw (381,52.4) node [anchor=north west][inner sep=0.75pt]    {$Y$};
		\draw (267,262.4) node [anchor=north west][inner sep=0.75pt]    {$Z$};
		\draw (235,102.4) node [anchor=north west][inner sep=0.75pt]    {$\textcolor[rgb]{0.82,0.01,0.11}{I\times \mathbf{q}}$};
		\draw (381,212.4) node [anchor=north west][inner sep=0.75pt]    {$\textcolor[rgb]{0.82,0.01,0.11}{I\times }\mathbf{\textcolor[rgb]{0.82,0.01,0.11}{p}}$};
		\draw (591,212.4) node [anchor=north west][inner sep=0.75pt]    {$\textcolor[rgb]{0.82,0.01,0.11}{I\times }\textcolor[rgb]{0.82,0.01,0.11}{\mathbf{q}}$};
		\draw (617,52.4) node [anchor=north west][inner sep=0.75pt]    {$Z$};
		\draw (139,61.4) node [anchor=north west][inner sep=0.75pt]    {$Y\#_{2} Z$};
		\draw (491,252.4) node [anchor=north west][inner sep=0.75pt]    {$Y\#_{2} Z$};
\end{tikzpicture}}

\newcommand{\figfourballs}{
	\tikzset{every picture/.style={line width=0.75pt}} 
	\begin{tikzpicture}[x=0.75pt,y=0.75pt,yscale=-1,xscale=1,scale=0.75]
		\fill[black!12] (60,50) circle[radius=32];
		\fill[black!12] (60,147) circle[radius=32];
		\fill[black!12] (430,50) circle[radius=32];
		\fill[black!12] (430,147) circle[radius=32];
		\draw    (20,50) -- (220,50) ;
		\draw    (20,10) -- (20,50) ;
		\draw    (220,10) -- (220,50) ;
		\draw    (270,50) -- (470,50) ;
		\draw    (270,10) -- (270,50) ;
		\draw    (470,10) -- (470,50) ;
		\draw  [dash pattern={on 0.84pt off 2.51pt}]  (0,10) -- (490,10) ;
		\draw    (20,147) -- (220,147) ;
		\draw    (20,147) -- (20,187) ;
		\draw    (220,147) -- (220,187) ;
		\draw    (270,147) -- (470,147) ;
		\draw    (270,147) -- (270,187) ;
		\draw    (470,147) -- (470,187) ;
		\draw  [dash pattern={on 0.84pt off 2.51pt}]  (0,187) -- (490,187) ;
		\draw  [color=black,dashed  ,draw opacity=1 ] (398,147) .. controls (398,129.33) and (412.33,115) .. (430,115) .. controls (447.67,115) and (462,129.33) .. (462,147) .. controls (462,164.67) and (447.67,179) .. (430,179) .. controls (412.33,179) and (398,164.67) .. (398,147) -- cycle ;
		\draw  [fill={rgb, 255:red, 0; green, 0; blue, 0 }  ,fill opacity=1 ] (64,50) .. controls (64,47.792) and (62.208,46) .. (60,46) .. controls (57.792,46) and (56,47.792) .. (56,50) .. controls (56,52.208) and (57.792,54) .. (60,54) .. controls (62.208,54) and (64,52.208) .. (64,50) -- cycle ;
		\draw  [fill={rgb, 255:red, 0; green, 0; blue, 0 }  ,fill opacity=1 ] (64,147) .. controls (64,144.792) and (62.208,143) .. (60,143) .. controls (57.792,143) and (56,144.792) .. (56,147) .. controls (56,149.208) and (57.792,151) .. (60,151) .. controls (62.208,151) and (64,149.208) .. (64,147) -- cycle ;
		\draw  [fill={rgb, 255:red, 0; green, 0; blue, 0 }  ,fill opacity=1 ] (434,50) .. controls (434,47.792) and (432.208,46) .. (430,46) .. controls (427.792,46) and (426,47.792) .. (426,50) .. controls (426,52.208) and (427.792,54) .. (430,54) .. controls (432.208,54) and (434,52.208) .. (434,50) -- cycle ;
		\draw  [fill={rgb, 255:red, 0; green, 0; blue, 0 }  ,fill opacity=1 ] (434,147) .. controls (434,144.792) and (432.208,143) .. (430,143) .. controls (427.792,143) and (426,144.792) .. (426,147) .. controls (426,149.208) and (427.792,151) .. (430,151) .. controls (432.208,151) and (434,149.208) .. (434,147) -- cycle ;
		\draw  [fill={rgb, 255:red, 0; green, 0; blue, 0 }  ,fill opacity=1 ] (174,50) .. controls (174,47.792) and (172.208,46) .. (170,46) .. controls (167.792,46) and (166,47.792) .. (166,50) .. controls (166,52.208) and (167.792,54) .. (170,54) .. controls (172.208,54) and (174,52.208) .. (174,50) -- cycle ;
		\draw  [fill={rgb, 255:red, 0; green, 0; blue, 0 }  ,fill opacity=1 ] (324,50) .. controls (324,47.792) and (322.208,46) .. (320,46) .. controls (317.792,46) and (316,47.792) .. (316,50) .. controls (316,52.208) and (317.792,54) .. (320,54) .. controls (322.208,54) and (324,52.208) .. (324,50) -- cycle ;
		\draw  [color=black,dashed  ,draw opacity=1 ] (398,50) .. controls (398,32.33) and (412.33,18) .. (430,18) .. controls (447.67,18) and (462,32.33) .. (462,50) .. controls (462,67.67) and (447.67,82) .. (430,82) .. controls (412.33,82) and (398,67.67) .. (398,50) -- cycle ;
		\draw  [color=black,dashed  ,draw opacity=1 ] (28,50) .. controls (28,32.33) and (42.33,18) .. (60,18) .. controls (77.67,18) and (92,32.33) .. (92,50) .. controls (92,67.67) and (77.67,82) .. (60,82) .. controls (42.33,82) and (28,67.67) .. (28,50) -- cycle ;
		\draw  [color=black,dashed  ,draw opacity=1 ] (28,147) .. controls (28,129.33) and (42.33,115) .. (60,115) .. controls (77.67,115) and (92,129.33) .. (92,147) .. controls (92,164.67) and (77.67,179) .. (60,179) .. controls (42.33,179) and (28,164.67) .. (28,147) -- cycle ;
		\draw [color={rgb, 255:red, 0; green, 0; blue, 255}, opacity=1][line width=1.5]    (60,50) -- (170,50) ;
		\draw [color={rgb, 255:red, 0; green, 0; blue, 255}, opacity=1][line width=1.5]    (320,50) -- (430,50) ;
		\draw (172,53.4) node [anchor=north west][inner sep=0.75pt]    {$p_{1}$};
		\draw (301,52.4) node [anchor=north west][inner sep=0.75pt]    {$p_{2}$};
		\draw (42,51.4) node [anchor=north west][inner sep=0.75pt]    {$x_{1}$};
		\draw (432,51.4) node [anchor=north west][inner sep=0.75pt]    {$x_{2}$};
		\draw (63,148.4) node [anchor=north west][inner sep=0.75pt]    {$z_{1}$};
		\draw (413,148.4) node [anchor=north west][inner sep=0.75pt]    {$z_{2}$};
		\draw (91,62.4) node [anchor=north west][inner sep=0.75pt]  [color=black  ,opacity=1 ]  {$B_{x1}$};
		\draw (91,121.4) node [anchor=north west][inner sep=0.75pt]  [color=black  ,opacity=1 ]  {$B_{z1}$};
		\draw (374,62.4) node [anchor=north west][inner sep=0.75pt]  [color=black  ,opacity=1 ]  {$B_{x2}$};
		\draw (360,113.4) node [anchor=north west][inner sep=0.75pt]  [color=black  ,opacity=1 ]  {$B_{z2}$};
		\draw (131,32.4) node [anchor=north west][inner sep=0.75pt]  [color={rgb, 255:red, 0; green, 0; blue, 255}, opacity=1]  {$\Gamma _{1}$};
		\draw (341,31.4) node [anchor=north west][inner sep=0.75pt] [color={rgb, 255:red, 0; green, 0; blue, 255}, opacity=1]  {$\Gamma _{2}$};
\end{tikzpicture}}

\newcommand{\figfourhypersurfaces}{
\tikzset{every picture/.style={line width=0.75pt}} 
\begin{tikzpicture}[x=0.75pt,y=0.75pt,yscale=-1,xscale=1]
\draw  [color={rgb, 255:red, 128; green, 128; blue, 128 }  ,draw opacity=1 ] (70,50) -- (410,50) -- (410,190) -- (70,190) -- cycle ;
\draw    (190,50) -- (190,190) ;
\draw    (290,50) -- (290,190) ;
\draw [color={rgb, 255:red, 128; green, 128; blue, 128 },draw opacity=1] (240,50) -- (240,125) (240,135) -- (240,190) ;
\draw [color={rgb, 255:red, 208; green, 2; blue, 27 },draw opacity=1] (70,130) -- (235,130) (245,130) -- (410,130) ;
\draw    (100,50) -- (100,120) ;
\draw    (100,140) -- (100,190) ;
\draw    (100,120) -- (260,120) ;
\draw    (100,140) -- (260,140) ;
\draw    (260,120) -- (260,140) ;
\draw    (220,110) -- (220,150) ;
\draw    (220,110) -- (380,110) ;
\draw    (220,150) -- (380,150) ;
\draw    (380,50) -- (380,110) ;
\draw    (380,150) -- (380,190) ;
\draw   (235,130) .. controls (235,127.24) and (237.24,125) .. (240,125) .. controls (242.76,125) and (245,127.24) .. (245,130) .. controls (245,132.76) and (242.76,135) .. (240,135) .. controls (237.24,135) and (235,132.76) .. (235,130) -- cycle ;
\draw (59,192.4) node [anchor=north west][inner sep=0.75pt]  [color={rgb, 255:red, 128; green, 128; blue, 128 }  ,opacity=1 ]  {$Y_{0}$};
\draw (409,191.4) node [anchor=north west][inner sep=0.75pt]  [color={rgb, 255:red, 128; green, 128; blue, 128 }  ,opacity=1 ]  {$Y_{2}$};
\draw (91,191.4) node [anchor=north west][inner sep=0.75pt]    {$Y_{0+}$};
\draw (371,192.4) node [anchor=north west][inner sep=0.75pt]    {$Y_{2-}$};
\draw (231,191.4) node [anchor=north west][inner sep=0.75pt]  [color={rgb, 255:red, 128; green, 128; blue, 128 }  ,opacity=1 ]  {$Y_{1}$};
\draw (181,192.4) node [anchor=north west][inner sep=0.75pt]    {$Y_{1-}$};
\draw (281,192.4) node [anchor=north west][inner sep=0.75pt]    {$Y_{1+}$};
\end{tikzpicture}}

\newcommand{\figskeinF}{
\tikzset{every picture/.style={line width=0.75pt}} 
\begin{tikzpicture}[x=0.75pt,y=0.75pt,yscale=-1,xscale=1,scale=0.7]
	\definecolor{skeinSurfaceFront}{RGB}{225,225,225}
	\definecolor{skeinSurfaceBack}{RGB}{185,185,185}
	\begin{scope}
		\clip (115,125) circle[radius=72.5];
		\fill[skeinSurfaceFront]
		(114.227,124.806) -- (61.76,52.43) -- (0,30)
		-- (0,230) -- (40,210) -- cycle;
		\fill[skeinSurfaceBack]
		(115.080,123.827) -- (177.94,51.68) -- (205,30)
		-- (205,230) -- (180,210) -- cycle;
	\end{scope}
	\begin{scope}
		\clip (292.5,125) circle[radius=72.5];
		\fill[skeinSurfaceFront]
		(294.356,122.475) -- (241.8,52.4) -- (205,30)
		-- (205,230) -- (230,210) -- cycle;
		\fill[skeinSurfaceBack]
		(295,123.333) -- (348.2,52.4) -- (385,30)
		-- (385,230) -- (360,210) -- cycle;
	\end{scope}
	\begin{scope}
		\clip (477.5,122.5) circle[radius=72.5];
		\fill[skeinSurfaceFront]
		(420,210) .. controls (459.4,180.45) and (473.49,76.12)
		.. (424.77,51.09) -- (390,30) -- (390,230) -- cycle;
		\fill[skeinSurfaceBack]
		(530,210) .. controls (478.2,205) and (476.94,72)
		.. (527.65,50.89) -- (565,30) -- (565,230) -- cycle;
	\end{scope}
	\draw    (40,210) -- (177.94,51.68) ;
	\draw [shift={(179.91,49.42)}, rotate = 131.07] [fill={rgb, 255:red, 0; green, 0; blue, 0 }  ][line width=0.08]  [draw opacity=0] (8.93,-4.29) -- (0,0) -- (8.93,4.29) -- cycle    ;
	\draw    (111.41,120.92) -- (61.76,52.43) ;
	\draw [shift={(60,50)}, rotate = 54.06] [fill={rgb, 255:red, 0; green, 0; blue, 0 }  ][line width=0.08]  [draw opacity=0] (8.93,-4.29) -- (0,0) -- (8.93,4.29) -- cycle    ;
	\draw    (117.41,126.92) -- (180,210) ;
	\draw    (297.5,120) -- (348.2,52.4) ;
	\draw [shift={(350,50)}, rotate = 126.87] [fill={rgb, 255:red, 0; green, 0; blue, 0 }  ][line width=0.08]  [draw opacity=0] (8.93,-4.29) -- (0,0) -- (8.93,4.29) -- cycle    ;
	\draw    (360,210) -- (241.8,52.4) ;
	\draw [shift={(240,50)}, rotate = 53.13] [fill={rgb, 255:red, 0; green, 0; blue, 0 }  ][line width=0.08]  [draw opacity=0] (8.93,-4.29) -- (0,0) -- (8.93,4.29) -- cycle    ;
	\draw    (292.5,125) -- (230,210) ;
	\draw    (420,210) .. controls (459.4,180.45) and (473.49,76.12) .. (424.77,51.09) ;
	\draw [shift={(422.5,50)}, rotate = 24.03] [fill={rgb, 255:red, 0; green, 0; blue, 0 }  ][line width=0.08]  [draw opacity=0] (8.93,-4.29) -- (0,0) -- (8.93,4.29) -- cycle    ;
	\draw    (530,210) .. controls (478.2,205) and (476.94,72) .. (527.65,50.89) ;
	\draw [shift={(530,50)}, rotate = 161.35] [fill={rgb, 255:red, 0; green, 0; blue, 0 }  ][line width=0.08]  [draw opacity=0] (8.93,-4.29) -- (0,0) -- (8.93,4.29) -- cycle    ;
	\draw  [dash pattern={on 0.84pt off 2.51pt}] (42.5,125) .. controls (42.5,84.96) and (74.96,52.5) .. (115,52.5) .. controls (155.04,52.5) and (187.5,84.96) .. (187.5,125) .. controls (187.5,165.04) and (155.04,197.5) .. (115,197.5) .. controls (74.96,197.5) and (42.5,165.04) .. (42.5,125) -- cycle ;
	\draw  [dash pattern={on 0.84pt off 2.51pt}] (220,125) .. controls (220,84.96) and (252.46,52.5) .. (292.5,52.5) .. controls (332.54,52.5) and (365,84.96) .. (365,125) .. controls (365,165.04) and (332.54,197.5) .. (292.5,197.5) .. controls (252.46,197.5) and (220,165.04) .. (220,125) -- cycle ;
	\draw  [dash pattern={on 0.84pt off 2.51pt}] (405,122.5) .. controls (405,82.46) and (437.46,50) .. (477.5,50) .. controls (517.54,50) and (550,82.46) .. (550,122.5) .. controls (550,162.54) and (517.54,195) .. (477.5,195) .. controls (437.46,195) and (405,162.54) .. (405,122.5) -- cycle ;
	\draw [preaction={draw=white,line width=3pt}] (105.02,119.57) .. controls (105.02,97.55) and (110.09,79.71) .. (116.36,79.71) .. controls (122.62,79.71) and (127.69,97.55) .. (127.69,119.57) .. controls (127.69,141.58) and (122.62,159.42) .. (116.36,159.42) .. controls (110.09,159.42) and (105.02,141.58) .. (105.02,119.57) -- cycle ;
	\draw [preaction={draw=white,line width=3pt}] (281.16,125) .. controls (281.16,102.99) and (286.24,85.14) .. (292.5,85.14) .. controls (298.76,85.14) and (303.84,102.99) .. (303.84,125) .. controls (303.84,147.01) and (298.76,164.86) .. (292.5,164.86) .. controls (286.24,164.86) and (281.16,147.01) .. (281.16,125) -- cycle ;
	\draw [preaction={draw=white,line width=3pt}] (123.403,114.275) -- (131.285,105.228);
	\draw [preaction={draw=white,line width=3pt}] (122.689,133.927) -- (129.910,143.511);
	\draw [preaction={draw=white,line width=3pt}] (299.654,117.127) -- (306.854,107.527);
	\draw [preaction={draw=white,line width=3pt}] (307.104,139.472) -- (299.904,129.872);
	\node[font=\large] at (115,215) {$L_+$};
	\node[font=\large] at (292.5,215) {$L_-$};
	\node[font=\large] at (477.5,215) {$L_0$};
	\node at (70,125) {$F_+$};
	\node at (250,125) {$F_-$};
	\node at (427,125) {$F_0$};
	\node at (530,125) {$F_0$};
	\node at (116.36,174) {$\gamma$};
	\node at (292.5,179.5) {$\gamma$};
	\end{tikzpicture}
}

\newcommand{\figcapcommute}{
\tikzset{every picture/.style={line width=0.75pt}} 
\begin{tikzpicture}[x=0.75pt,y=0.75pt,yscale=-1,xscale=1,scale=1]
\draw    (60,90) -- (60,190) ;
\draw    (160,90) -- (160,190) ;
\draw    (160,90) -- (260,90) ;
\draw    (260,90) -- (260,118.97) ;
\draw    (260,158.97) -- (260,190) ;
\draw    (260,118.97) .. controls (188,118.97) and (191.5,159.47) .. (260,158.97) ;
\draw    (310,90.03) -- (410,90.03) ;
\draw    (310,190.03) -- (410,190.03) ;
\draw    (410,90.03) -- (410,119.03) ;
\draw    (410,160.03) -- (410,190.03) ;
\draw    (410,120.03) .. controls (338,120.03) and (341.5,160.53) .. (410,160.03) ;
\draw    (310,90.03) -- (310,190.03) ;
\draw    (410,120.03) -- (510,120.03) ;
\draw    (410,160.03) -- (510,160.03) ;
\draw    (510,90.03) -- (510,120.03) ;
\draw    (510,159.03) -- (510,191.03) ;
\draw [color={rgb, 255:red, 0; green, 0; blue, 255}, opacity=1]   (160,124.5) .. controls (224.91,125.92) and (218.41,154.42) .. (160,154.5) ;
\draw [color={rgb, 255:red, 0; green, 0; blue, 255}, opacity=1]  (60,124.5) -- (160,124.5) ;
\draw [color={rgb, 255:red, 0; green, 0; blue, 255}, opacity=1]  (60,154.5) -- (160,154.5) ;
\draw [color={rgb, 255:red, 0; green, 0; blue, 255}, opacity=1]  (310,129.03) .. controls (374.91,130.24) and (368.41,154.49) .. (310,154.56) ;
\draw    (410,190.03) .. controls (443.41,192.45) and (438.41,209.95) .. (460,211.03) .. controls (481.59,212.11) and (480.41,190.95) .. (510,191.03) ;
\draw    (160,190) -- (260,190) ;
\draw    (410,90.03) .. controls (442.41,89.45) and (443.52,69.09) .. (460,70.03) .. controls (476.48,70.97) and (476.41,89.45) .. (510,90.03) ;
\draw    (60,90) .. controls (92.41,89.42) and (93.52,69.06) .. (110,70) .. controls (126.48,70.94) and (126.41,89.42) .. (160,90) ;
\draw    (60,189) .. controls (93.41,191.42) and (88.41,208.92) .. (110,210) .. controls (131.59,211.08) and (130.41,189.92) .. (160,190) ;
\draw (51,191.37) node [anchor=north west][inner sep=0.75pt]    {$Y_{\scone}$};
\draw (151,192.4) node [anchor=north west][inner sep=0.75pt]    {$Y_{\sctwo}$};
\draw (341,168.43) node [anchor=north west][inner sep=0.75pt]    {$W_{F_{\scone}}$};
\draw (301,191.43) node [anchor=north west][inner sep=0.75pt]    {$Y_{\scone}$};
\draw (401,191.43) node [anchor=north west][inner sep=0.75pt]    {$Y_{\sctwo}$};
\draw (191,162.4) node [anchor=north west][inner sep=0.75pt]    {$W_{F_{\sctwo}}$};
\draw (101,172.4) node [anchor=north west][inner sep=0.75pt]    {$W$};
\draw (451,101.43) node [anchor=north west][inner sep=0.75pt]    {$Z$};
\draw (451,171.43) node [anchor=north west][inner sep=0.75pt]    {$Z$};
\draw (179,112.4) node [anchor=north west][inner sep=0.75pt]  [color={rgb, 255:red, 0; green, 0; blue, 255}, opacity=1]  {$F_{\sctwo}$};
\draw (329,111.43) node [anchor=north west][inner sep=0.75pt] [color={rgb, 255:red, 0; green, 0; blue, 255}, opacity=1]  {$F_{\scone}$};
\draw (71,107.9) node [anchor=north west][inner sep=0.75pt]  [color={rgb, 255:red, 0; green, 0; blue, 255}, opacity=1]  {$[0,1] \times \partial F$};
\end{tikzpicture}}

\newcommand{\figtwofivehypersurface}{
\tikzset{every picture/.style={line width=0.75pt}} 
\begin{tikzpicture}[x=0.75pt,y=0.75pt,yscale=-1,xscale=1,scale=0.8]
	\draw   (10,139) .. controls (10,111.39) and (18.95,89) .. (30,89) .. controls (41.05,89) and (50,111.39) .. (50,139) .. controls (50,166.61) and (41.05,189) .. (30,189) .. controls (18.95,189) and (10,166.61) .. (10,139) -- cycle ;
	\draw   (270,139) .. controls (270,111.39) and (278.95,89) .. (290,89) .. controls (301.05,89) and (310,111.39) .. (310,139) .. controls (310,166.61) and (301.05,189) .. (290,189) .. controls (278.95,189) and (270,166.61) .. (270,139) -- cycle ;
	\draw    (30,89) .. controls (33.75,87.42) and (40.98,90.94) .. (50,89) .. controls (59.02,87.06) and (95.75,61.5) .. (110,49) .. controls (124.25,36.5) and (137.75,19) .. (160,19) .. controls (182.25,19) and (198.25,36) .. (210,49) .. controls (221.75,62) and (261.75,87) .. (270,89) .. controls (278.25,91) and (282.75,88) .. (290,89) ;
	\draw    (30,189) .. controls (35.46,190.21) and (39.25,186.56) .. (50,189) .. controls (60.75,191.44) and (93.75,214) .. (110,229) .. controls (126.25,244) and (140.25,259) .. (160,259) .. controls (179.75,259) and (194.75,247.5) .. (210,229) .. controls (225.25,210.5) and (241.25,188.5) .. (270,189) .. controls (298.75,189.5) and (278.25,187.5) .. (290,189) ;
	\draw    (190.71,29.92) -- (189.71,248.92) ;
	\draw    (110,49) .. controls (119.75,69) and (139.77,79.13) .. (160,79) .. controls (180.23,78.87) and (196.25,72.5) .. (210,49) ;
	\draw    (110,229) .. controls (128.25,204.5) and (145.13,199.33) .. (160,199) .. controls (174.87,198.67) and (196.41,202.11) .. (210,229) ;
	\draw    (50,89) .. controls (56.86,97.86) and (58.33,108.08) .. (60,119) .. controls (61.67,129.92) and (185.21,111.42) .. (185.71,136.92) .. controls (186.21,162.42) and (61.18,153.2) .. (60,159) .. controls (58.82,164.8) and (53.53,183.65) .. (50,189) ;
	\draw    (270,89) .. controls (251.75,142.5) and (141.46,100.92) .. (136.71,138.92) .. controls (131.96,176.92) and (251.25,139) .. (270,189) ;
	\draw   (148.5,138.5) .. controls (148.5,132.98) and (152.98,128.5) .. (158.5,128.5) .. controls (164.02,128.5) and (168.5,132.98) .. (168.5,138.5) .. controls (168.5,144.02) and (164.02,148.5) .. (158.5,148.5) .. controls (152.98,148.5) and (148.5,144.02) .. (148.5,138.5) -- cycle ;
	\draw   (350,140) .. controls (350,112.39) and (358.95,90) .. (370,90) .. controls (381.05,90) and (390,112.39) .. (390,140) .. controls (390,167.61) and (381.05,190) .. (370,190) .. controls (358.95,190) and (350,167.61) .. (350,140) -- cycle ;
	\draw   (610,140) .. controls (610,112.39) and (618.95,90) .. (630,90) .. controls (641.05,90) and (650,112.39) .. (650,140) .. controls (650,167.61) and (641.05,190) .. (630,190) .. controls (618.95,190) and (610,167.61) .. (610,140) -- cycle ;
	\draw    (370,190) .. controls (433.5,189.5) and (424,240) .. (501,240.5) .. controls (578,241) and (558.5,190) .. (630,190) ;
	\draw    (370,90) .. controls (461.5,85.5) and (388,29.5) .. (499,29.5) .. controls (610,29.5) and (552.5,73.5) .. (630,90) ;
	\draw    (414.5,77.5) .. controls (434.5,84.5) and (514,60.5) .. (522.5,30) ;
	\draw    (461,32.5) -- (459.5,233.5) ;
	\draw    (549.5,33) -- (549.5,230) ;
	\draw    (482.5,30.5) .. controls (503.5,73.5) and (557.5,92.5) .. (583,62.5) ;
	\draw    (391.5,88) .. controls (463,122.5) and (503.5,133.5) .. (608,83) ;
	\draw    (416.5,203) .. controls (424.5,190.5) and (489.5,184.5) .. (525,237.5) ;
	\draw    (476,237.5) .. controls (484,225) and (529.5,188) .. (582,203.5) ;
	\draw    (400.5,195.5) .. controls (460.5,128.5) and (558.5,147) .. (602,193.5) ;
	\draw    (129.71,30.42) -- (130.21,247.42) ;
	\draw [color={rgb, 255:red, 208; green, 2; blue, 27 }  ,draw opacity=1 ]   (50,139) -- (148.5,138.5) ;
	\draw [color={rgb, 255:red, 208; green, 2; blue, 27 }  ,draw opacity=1 ]   (168.5,138.5) -- (270,139) ;
	
	\draw (40,192.4) node [anchor=north west][inner sep=0.75pt]    {$Y_{i,+}$};
	\draw (260,191.4) node [anchor=north west][inner sep=0.75pt]    {$Y_{i+2,-}$};
	\draw (110.5,252.9) node [anchor=north west][inner sep=0.75pt]    {$Y_{i+1,-}$};
	\draw (201,62.4) node [anchor=north west][inner sep=0.75pt]    {$S_{i}$};
	\draw (191,191.4) node [anchor=north west][inner sep=0.75pt]    {$\tau S_{i}$};
	\draw (217.5,139.4) node [anchor=north west][inner sep=0.75pt]  [color={rgb, 255:red, 208; green, 2; blue, 27 }  ,opacity=1 ]  {$I\times p_{*}$};
	\draw (361,193.4) node [anchor=north west][inner sep=0.75pt]    {$Y_{i}$};
	\draw (628,192.4) node [anchor=north west][inner sep=0.75pt]    {$Y_{i}$};
	\draw (439.5,238.4) node [anchor=north west][inner sep=0.75pt]    {$Y_{i+1}$};
	\draw (547.5,233.4) node [anchor=north west][inner sep=0.75pt]    {$Y_{i+2}$};
	\draw (428.5,195.4) node [anchor=north west][inner sep=0.75pt]    {$\tau S_{i}$};
	\draw (433,57.4) node [anchor=north west][inner sep=0.75pt]    {$S_{i}$};
	\draw (583,47.9) node [anchor=north west][inner sep=0.75pt]    {$S_{i+1}$};
	\draw (578,204.4) node [anchor=north west][inner sep=0.75pt]    {$\tau S_{i+1}$};
	\draw (497,98.4) node [anchor=north west][inner sep=0.75pt]    {$R$};
	\draw (490,153.9) node [anchor=north west][inner sep=0.75pt]    {$\tau R$};
	\draw (187.5,250.9) node [anchor=north west][inner sep=0.75pt]    {$Y_{i+1,+}$};
\end{tikzpicture}}

\newcommand{\figbands}{
\tikzset{every picture/.style={line width=0.75pt}} 
\begin{tikzpicture}[x=0.75pt,y=0.75pt,yscale=-1,xscale=1,scale=0.7]
\fill[black!15] (150,20) rectangle (260,60);
\fill[black!15] (150,200) rectangle (260,240);
\fill[black!15] (570,20) rectangle (640,60);
\fill[black!15] (430,130) rectangle (640,170);
\draw    (10,10) -- (10,250) ;
\draw    (7.5,10) -- (150,10) ;
\draw    (7.5,250) -- (150,250) ;
\draw    (150,60) -- (150,200) ;
\draw    (260,10) -- (340,10) ;
\draw    (340,10) -- (340,250) ;
\draw    (260,60) -- (260,200) ;
\draw    (150,10) -- (150,20) ;
\draw    (150,240) -- (150,250) ;
\draw    (260,250) -- (340,250) ;
\draw    (260,240) -- (260,250) ;
\draw    (260,10) -- (260,20) ;
\draw    (150,20) -- (260,20) ;
\draw    (150,60) -- (260,60) ;
\draw    (150,200) -- (260,200) ;
\draw    (150,240) -- (260,240) ;
\draw[fill=black]   (147.5,40) .. controls (147.5,38.62) and (148.62,37.5) .. (150,37.5) .. controls (151.38,37.5) and (152.5,38.62) .. (152.5,40) .. controls (152.5,41.38) and (151.38,42.5) .. (150,42.5) .. controls (148.62,42.5) and (147.5,41.38) .. (147.5,40) -- cycle ;
\draw[fill=black]   (147.5,220) .. controls (147.5,218.62) and (148.62,217.5) .. (150,217.5) .. controls (151.38,217.5) and (152.5,218.62) .. (152.5,220) .. controls (152.5,221.38) and (151.38,222.5) .. (150,222.5) .. controls (148.62,222.5) and (147.5,221.38) .. (147.5,220) -- cycle ;
\draw[fill=black]   (257.5,40) .. controls (257.5,38.62) and (258.62,37.5) .. (260,37.5) .. controls (261.38,37.5) and (262.5,38.62) .. (262.5,40) .. controls (262.5,41.38) and (261.38,42.5) .. (260,42.5) .. controls (258.62,42.5) and (257.5,41.38) .. (257.5,40) -- cycle ;
\draw[fill=black]   (257.5,220) .. controls (257.5,218.62) and (258.62,217.5) .. (260,217.5) .. controls (261.38,217.5) and (262.5,218.62) .. (262.5,220) .. controls (262.5,221.38) and (261.38,222.5) .. (260,222.5) .. controls (258.62,222.5) and (257.5,221.38) .. (257.5,220) -- cycle ;
\draw[fill=black]   (7.5,120) .. controls (7.5,118.62) and (8.62,117.5) .. (10,117.5) .. controls (11.38,117.5) and (12.5,118.62) .. (12.5,120) .. controls (12.5,121.38) and (11.38,122.5) .. (10,122.5) .. controls (8.62,122.5) and (7.5,121.38) .. (7.5,120) -- cycle ;
\draw[fill=black]   (257.5,120) .. controls (257.5,118.62) and (258.62,117.5) .. (260,117.5) .. controls (261.38,117.5) and (262.5,118.62) .. (262.5,120) .. controls (262.5,121.38) and (261.38,122.5) .. (260,122.5) .. controls (258.62,122.5) and (257.5,121.38) .. (257.5,120) -- cycle ;
\draw  [dash pattern={on 0.84pt off 2.51pt}]  (150,20) -- (150,60) ;
\draw  [dash pattern={on 0.84pt off 2.51pt}]  (260,20) -- (260,60) ;
\draw  [dash pattern={on 0.84pt off 2.51pt}]  (150,200) -- (150,240) ;
\draw  [dash pattern={on 0.84pt off 2.51pt}]  (260,200) -- (260,240) ;
\draw    (390,10) -- (390,250) ;
\draw    (390,10) -- (570,10) ;
\draw    (390,250) -- (570,250) ;
\draw    (430,80) -- (430,130) ;
\draw    (640,10) -- (720,10) ;
\draw    (720,10) -- (720,250) ;
\draw    (540,80) -- (540,89.64) ;
\draw    (540,89.64) -- (540,99.28) ;
\draw    (540,99.28) -- (540,108.93) ;
\draw    (540,108.93) -- (540,118.57) ;
\draw    (540,118.57) -- (540,128.21) ;
\draw    (540,170.42) -- (540,190) ;
\draw    (640,250) -- (720,250) ;
\draw    (430,80) -- (540,80) ;
\draw    (430,190) -- (540,190) ;
\draw[fill=black]   (387.5,120) .. controls (387.5,118.62) and (388.62,117.5) .. (390,117.5) .. controls (391.38,117.5) and (392.5,118.62) .. (392.5,120) .. controls (392.5,121.38) and (391.38,122.5) .. (390,122.5) .. controls (388.62,122.5) and (387.5,121.38) .. (387.5,120) -- cycle ;
\draw    (570,60) -- (570,130) ;
\draw    (640,60) -- (640,130) ;
\draw    (570,10) -- (570,20) ;
\draw    (640,10) -- (640,20) ;
\draw    (570,20) -- (640,20) ;
\draw[fill=black]   (567.5,40) .. controls (567.5,38.62) and (568.62,37.5) .. (570,37.5) .. controls (571.38,37.5) and (572.5,38.62) .. (572.5,40) .. controls (572.5,41.38) and (571.38,42.5) .. (570,42.5) .. controls (568.62,42.5) and (567.5,41.38) .. (567.5,40) -- cycle ;
\draw  [dash pattern={on 0.84pt off 2.51pt}]  (570,20) -- (570,60) ;
\draw[fill=black]   (636.5,40) .. controls (636.5,38.62) and (637.62,37.5) .. (639,37.5) .. controls (640.38,37.5) and (641.5,38.62) .. (641.5,40) .. controls (641.5,41.38) and (640.38,42.5) .. (639,42.5) .. controls (637.62,42.5) and (636.5,41.38) .. (636.5,40) -- cycle ;
\draw  [dash pattern={on 0.84pt off 2.51pt}]  (639,20) -- (639,60) ;
\draw    (570,60) -- (640,60) ;
\draw  [dashed]  (570,130) -- (570,170) ;
\draw    (570,170.42) -- (570,250) ;
\draw    (640,170.42) -- (640,250) ;
\draw[fill=black]   (637.5,150) .. controls (637.5,148.62) and (638.62,147.5) .. (640,147.5) .. controls (641.38,147.5) and (642.5,148.62) .. (642.5,150) .. controls (642.5,151.38) and (641.38,152.5) .. (640,152.5) .. controls (638.62,152.5) and (637.5,151.38) .. (637.5,150) -- cycle ;
\draw  [dash pattern={on 0.84pt off 2.51pt}]  (640,130) -- (640,170) ;
\draw[fill=black]   (427.5,150) .. controls (427.5,148.62) and (428.62,147.5) .. (430,147.5) .. controls (431.38,147.5) and (432.5,148.62) .. (432.5,150) .. controls (432.5,151.38) and (431.38,152.5) .. (430,152.5) .. controls (428.62,152.5) and (427.5,151.38) .. (427.5,150) -- cycle ;
\draw    (430,170.42) -- (430,190) ;
\draw  [dash pattern={on 0.84pt off 2.51pt}]  (430,130) -- (430,170) ;
\draw  [dashed]  (540,130.42) -- (540,170.42) ;
\draw    (430,130) -- (640,130) ;
\draw    (430,170) -- (640,170) ;
\draw[fill=black]   (427.5,90) .. controls (427.5,88.62) and (428.62,87.5) .. (430,87.5) .. controls (431.38,87.5) and (432.5,88.62) .. (432.5,90) .. controls (432.5,91.38) and (431.38,92.5) .. (430,92.5) .. controls (428.62,92.5) and (427.5,91.38) .. (427.5,90) -- cycle ;
\draw (142,40) node [anchor=east][inner sep=0.75pt]    {$r_{1}$};
\draw (142,220) node [anchor=east][inner sep=0.75pt]    {$r_{2}$};
\draw (268,40) node [anchor=west][inner sep=0.75pt]    {$z_{1}$};
\draw (268,220) node [anchor=west][inner sep=0.75pt]    {$z_{2}$};
\draw (18,120) node [anchor=west][inner sep=0.75pt]    {$p_{a}$};
\draw (268,120) node [anchor=west][inner sep=0.75pt]    {$p_{b}$};
\draw (61,210) node [anchor=north west][inner sep=0.75pt]    {$F_{0}$};
\draw (300,210) node [anchor=north west][inner sep=0.75pt]    {$D^{2}$};
\draw (398,120) node [anchor=west][inner sep=0.75pt]    {$p_{a}$};
\draw (441,210) node [anchor=north west][inner sep=0.75pt]    {$F_{0}$};
\draw (680,210) node [anchor=north west][inner sep=0.75pt]    {$D^{2}$};
\draw (562,40) node [anchor=east][inner sep=0.75pt]    {$r_{1}$};
\draw (647,40) node [anchor=west][inner sep=0.75pt]    {$z_{1}$};
\draw (648,150) node [anchor=west][inner sep=0.75pt]    {$z_{2}$};
\draw (422,150) node [anchor=east][inner sep=0.75pt]    {$r_{2}$};
\draw (438,90) node [anchor=west][inner sep=0.75pt]    {$p_{b}'$};
\end{tikzpicture}}

\newcommand{\figskeinintro}{
\tikzset{every picture/.style={line width=0.75pt}} 
\begin{tikzpicture}[x=0.75pt,y=0.75pt,yscale=-1,xscale=1,scale=0.5]
		\draw    (56.99,188.91) -- (177.93,51.67) ;
		\draw [shift={(179.91,49.42)}, rotate = 131.39] [fill={rgb, 255:red, 0; green, 0; blue, 0 }  ][line width=0.08]  [draw opacity=0] (8.93,-4.29) -- (0,0) -- (8.93,4.29) -- cycle    ;
		\draw    (111.41,120.92) -- (61.76,52.43) ;
		\draw [shift={(60,50)}, rotate = 54.06] [fill={rgb, 255:red, 0; green, 0; blue, 0 }  ][line width=0.08]  [draw opacity=0] (8.93,-4.29) -- (0,0) -- (8.93,4.29) -- cycle    ;
		\draw    (117.41,126.92) -- (168.49,191.41) ;
		\draw    (297.5,120) -- (348.2,52.4) ;
		\draw [shift={(350,50)}, rotate = 126.87] [fill={rgb, 255:red, 0; green, 0; blue, 0 }  ][line width=0.08]  [draw opacity=0] (8.93,-4.29) -- (0,0) -- (8.93,4.29) -- cycle    ;
		\draw    (346.49,192.91) -- (241.79,52.41) ;
		\draw [shift={(240,50)}, rotate = 53.31] [fill={rgb, 255:red, 0; green, 0; blue, 0 }  ][line width=0.08]  [draw opacity=0] (8.93,-4.29) -- (0,0) -- (8.93,4.29) -- cycle    ;
		\draw    (292.5,125) -- (246.49,188.91) ;
		\draw    (426.99,191.41) .. controls (466.39,161.86) and (473.7,75.56) .. (424.77,51.08) ;
		\draw [shift={(422.5,50)}, rotate = 24.03] [fill={rgb, 255:red, 0; green, 0; blue, 0 }  ][line width=0.08]  [draw opacity=0] (8.93,-4.29) -- (0,0) -- (8.93,4.29) -- cycle    ;
		\draw    (531.49,194.91) .. controls (479.7,189.91) and (476.01,72.04) .. (526.65,51.38) ;
		\draw [shift={(529,50.5)}, rotate = 161.35] [fill={rgb, 255:red, 0; green, 0; blue, 0 }  ][line width=0.08]  [draw opacity=0] (8.93,-4.29) -- (0,0) -- (8.93,4.29) -- cycle    ;
		\draw  [dash pattern={on 0.84pt off 2.51pt}] (42.5,125) .. controls (42.5,84.96) and (74.96,52.5) .. (115,52.5) .. controls (155.04,52.5) and (187.5,84.96) .. (187.5,125) .. controls (187.5,165.04) and (155.04,197.5) .. (115,197.5) .. controls (74.96,197.5) and (42.5,165.04) .. (42.5,125) -- cycle ;
		\draw  [dash pattern={on 0.84pt off 2.51pt}] (220,125) .. controls (220,84.96) and (252.46,52.5) .. (292.5,52.5) .. controls (332.54,52.5) and (365,84.96) .. (365,125) .. controls (365,165.04) and (332.54,197.5) .. (292.5,197.5) .. controls (252.46,197.5) and (220,165.04) .. (220,125) -- cycle ;
		\draw  [dash pattern={on 0.84pt off 2.51pt}] (405,122.5) .. controls (405,82.46) and (437.46,50) .. (477.5,50) .. controls (517.54,50) and (550,82.46) .. (550,122.5) .. controls (550,162.54) and (517.54,195) .. (477.5,195) .. controls (437.46,195) and (405,162.54) .. (405,122.5) -- cycle ;
		\draw (101.5,201.4) node [anchor=north west][inner sep=0.75pt]    {$L_{+}$};
		\draw (281.5,202.4) node [anchor=north west][inner sep=0.75pt]    {$L_{-}$};
		\draw (469,200.9) node [anchor=north west][inner sep=0.75pt]    {$L_{0}$};
\end{tikzpicture}}

\newcommand{\figselfmixed}{
\tikzset{every picture/.style={line width=0.75pt}} 
\begin{tikzpicture}[x=0.75pt,y=0.75pt,yscale=-1,xscale=1,scale=0.9]
	\draw  [dash pattern={on 0.84pt off 2.51pt}] (21.82,152.68) .. controls (21.82,136.11) and (35.25,122.68) .. (51.82,122.68) .. controls (68.38,122.68) and (81.82,136.11) .. (81.82,152.68) .. controls (81.82,169.25) and (68.38,182.68) .. (51.82,182.68) .. controls (35.25,182.68) and (21.82,169.25) .. (21.82,152.68) -- cycle ;
	\draw    (61.82,122.68) .. controls (124.32,12.18) and (113.32,296.18) .. (61.82,182.68) ;
	\draw    (41.82,122.68) .. controls (-10.68,31.68) and (10.82,277.18) .. (41.82,182.68) ;
	\draw  [dash pattern={on 0.84pt off 2.51pt}] (127.15,152.68) .. controls (127.15,136.11) and (140.58,122.68) .. (157.15,122.68) .. controls (173.72,122.68) and (187.15,136.11) .. (187.15,152.68) .. controls (187.15,169.25) and (173.72,182.68) .. (157.15,182.68) .. controls (140.58,182.68) and (127.15,169.25) .. (127.15,152.68) -- cycle ;
	\draw    (167.15,122.68) .. controls (229.65,12.18) and (218.65,296.18) .. (167.15,182.68) ;
	\draw    (147.15,122.68) .. controls (94.65,31.68) and (116.15,277.18) .. (147.15,182.68) ;
	\draw  [dash pattern={on 0.84pt off 2.51pt}] (231.82,152.68) .. controls (231.82,136.11) and (245.25,122.68) .. (261.82,122.68) .. controls (278.38,122.68) and (291.82,136.11) .. (291.82,152.68) .. controls (291.82,169.25) and (278.38,182.68) .. (261.82,182.68) .. controls (245.25,182.68) and (231.82,169.25) .. (231.82,152.68) -- cycle ;
	\draw    (271.82,122.68) .. controls (334.32,12.18) and (323.32,296.18) .. (271.82,182.68) ;
	\draw    (251.82,122.68) .. controls (199.32,31.68) and (220.82,277.18) .. (251.82,182.68) ;
	\draw    (392.08,163) .. controls (367.94,200.2) and (325.28,199.99) .. (345.11,142.27) .. controls (364.94,84.55) and (436.46,64.93) .. (421.47,133.87) ;
	\draw  [dash pattern={on 0.84pt off 2.51pt}] (406.33,129.4) .. controls (417.6,127.03) and (428.8,134.95) .. (431.35,147.09) .. controls (433.9,159.23) and (426.84,170.99) .. (415.57,173.36) .. controls (404.3,175.73) and (393.09,167.81) .. (390.54,155.67) .. controls (387.99,143.53) and (395.06,131.77) .. (406.33,129.4) -- cycle ;
	\draw    (394.26,139.59) .. controls (374.33,108.97) and (341.38,158.35) .. (367.95,183.37) ;
	\draw    (427.63,163.17) .. controls (425.62,190.75) and (383.9,242.36) .. (369.49,190.69) ;
	\draw    (492.08,163) .. controls (467.94,200.2) and (425.28,199.99) .. (445.11,142.27) .. controls (464.94,84.55) and (536.46,64.93) .. (521.47,133.87) ;
	\draw  [dash pattern={on 0.84pt off 2.51pt}] (506.33,129.4) .. controls (517.6,127.03) and (528.8,134.95) .. (531.35,147.09) .. controls (533.9,159.23) and (526.84,170.99) .. (515.57,173.36) .. controls (504.3,175.73) and (493.09,167.81) .. (490.54,155.67) .. controls (487.99,143.53) and (495.06,131.77) .. (506.33,129.4) -- cycle ;
	\draw    (494.26,139.59) .. controls (474.33,108.97) and (441.38,158.35) .. (467.95,183.37) ;
	\draw    (527.63,163.17) .. controls (525.62,190.75) and (483.9,242.36) .. (469.49,190.69) ;
	\draw    (592.08,163) .. controls (567.94,200.2) and (525.28,199.99) .. (545.11,142.27) .. controls (564.94,84.55) and (636.46,64.93) .. (621.47,133.87) ;
	\draw  [dash pattern={on 0.84pt off 2.51pt}] (606.33,129.4) .. controls (617.6,127.03) and (628.8,134.95) .. (631.35,147.09) .. controls (633.9,159.23) and (626.84,170.99) .. (615.57,173.36) .. controls (604.3,175.73) and (593.09,167.81) .. (590.54,155.67) .. controls (587.99,143.53) and (595.06,131.77) .. (606.33,129.4) -- cycle ;
	\draw    (594.26,139.59) .. controls (574.33,108.97) and (541.38,158.35) .. (567.95,183.37) ;
	\draw    (627.63,163.17) .. controls (625.62,190.75) and (583.9,242.36) .. (569.49,190.69) ;
	\draw    (60.57,127.76) -- (42.07,180.43) ;
	\draw [shift={(61.57,124.93)}, rotate = 109.36] [fill={rgb, 255:red, 0; green, 0; blue, 0 }  ][line width=0.08]  [draw opacity=0] (8.93,-4.29) -- (0,0) -- (8.93,4.29) -- cycle    ;
	\draw    (43.24,127.91) -- (49,150.5) ;
	\draw [shift={(42.5,125)}, rotate = 75.7] [fill={rgb, 255:red, 0; green, 0; blue, 0 }  ][line width=0.08]  [draw opacity=0] (8.93,-4.29) -- (0,0) -- (8.93,4.29) -- cycle    ;
	\draw    (53.23,157.38) -- (62,180.5) ;
	\draw    (419.56,137.05) -- (392.3,164.06) ;
	\draw [shift={(421.7,134.94)}, rotate = 135.26] [fill={rgb, 255:red, 0; green, 0; blue, 0 }  ][line width=0.08]  [draw opacity=0] (8.93,-4.29) -- (0,0) -- (8.93,4.29) -- cycle    ;
	\draw    (396.65,141.4) -- (404,147) ;
	\draw [shift={(394.26,139.59)}, rotate = 37.28] [fill={rgb, 255:red, 0; green, 0; blue, 0 }  ][line width=0.08]  [draw opacity=0] (8.93,-4.29) -- (0,0) -- (8.93,4.29) -- cycle    ;
	\draw    (410.95,151.38) -- (427.63,163.17) ;
	\draw    (520.79,136.25) -- (511.5,147.5) ;
	\draw [shift={(522.7,133.94)}, rotate = 129.54] [fill={rgb, 255:red, 0; green, 0; blue, 0 }  ][line width=0.08]  [draw opacity=0] (8.93,-4.29) -- (0,0) -- (8.93,4.29) -- cycle    ;
	\draw    (496.71,141.32) -- (527.63,163.17) ;
	\draw [shift={(494.26,139.59)}, rotate = 35.25] [fill={rgb, 255:red, 0; green, 0; blue, 0 }  ][line width=0.08]  [draw opacity=0] (8.93,-4.29) -- (0,0) -- (8.93,4.29) -- cycle    ;
	\draw    (507.5,152) -- (492.08,163) ;
	\draw    (592.08,163) .. controls (598.98,160.24) and (599.69,152.05) .. (595.47,142.19) ;
	\draw [shift={(594.26,139.59)}, rotate = 63.43] [fill={rgb, 255:red, 0; green, 0; blue, 0 }  ][line width=0.08]  [draw opacity=0] (8.93,-4.29) -- (0,0) -- (8.93,4.29) -- cycle    ;
	\draw    (627.63,163.17) .. controls (613.02,160.67) and (610.75,146.96) .. (619.49,136.09) ;
	\draw [shift={(621.47,133.87)}, rotate = 134.58] [fill={rgb, 255:red, 0; green, 0; blue, 0 }  ][line width=0.08]  [draw opacity=0] (8.93,-4.29) -- (0,0) -- (8.93,4.29) -- cycle    ;
	\draw    (166.73,125.65) -- (162,159) ;
	\draw [shift={(167.15,122.68)}, rotate = 98.07] [fill={rgb, 255:red, 0; green, 0; blue, 0 }  ][line width=0.08]  [draw opacity=0] (8.93,-4.29) -- (0,0) -- (8.93,4.29) -- cycle    ;
	\draw    (148.1,125.53) -- (167.15,182.68) ;
	\draw [shift={(147.15,122.68)}, rotate = 71.57] [fill={rgb, 255:red, 0; green, 0; blue, 0 }  ][line width=0.08]  [draw opacity=0] (8.93,-4.29) -- (0,0) -- (8.93,4.29) -- cycle    ;
	\draw    (157,164.5) -- (147.15,182.68) ;
	\draw    (251.82,182.68) .. controls (260.73,146.17) and (254.24,136.48) .. (252.23,125.57) ;
	\draw [shift={(251.82,122.68)}, rotate = 84.51] [fill={rgb, 255:red, 0; green, 0; blue, 0 }  ][line width=0.08]  [draw opacity=0] (8.93,-4.29) -- (0,0) -- (8.93,4.29) -- cycle    ;
	\draw    (271.82,182.68) .. controls (265.97,158.46) and (266.82,137.4) .. (270.79,125.41) ;
	\draw [shift={(271.82,122.68)}, rotate = 113.05] [fill={rgb, 255:red, 0; green, 0; blue, 0 }  ][line width=0.08]  [draw opacity=0] (8.93,-4.29) -- (0,0) -- (8.93,4.29) -- cycle    ;
	\draw  [fill={rgb, 255:red, 0; green, 0; blue, 0 }  ,fill opacity=1 ] (348,187.25) .. controls (348,186.28) and (348.78,185.5) .. (349.75,185.5) .. controls (350.72,185.5) and (351.5,186.28) .. (351.5,187.25) .. controls (351.5,188.22) and (350.72,189) .. (349.75,189) .. controls (348.78,189) and (348,188.22) .. (348,187.25) -- cycle ;
	\draw  [fill={rgb, 255:red, 0; green, 0; blue, 0 }  ,fill opacity=1 ] (388,212.75) .. controls (388,211.78) and (388.78,211) .. (389.75,211) .. controls (390.72,211) and (391.5,211.78) .. (391.5,212.75) .. controls (391.5,213.72) and (390.72,214.5) .. (389.75,214.5) .. controls (388.78,214.5) and (388,213.72) .. (388,212.75) -- cycle ;
	\draw  [fill={rgb, 255:red, 0; green, 0; blue, 0 }  ,fill opacity=1 ] (451,188.25) .. controls (451,187.28) and (451.78,186.5) .. (452.75,186.5) .. controls (453.72,186.5) and (454.5,187.28) .. (454.5,188.25) .. controls (454.5,189.22) and (453.72,190) .. (452.75,190) .. controls (451.78,190) and (451,189.22) .. (451,188.25) -- cycle ;
	\draw  [fill={rgb, 255:red, 0; green, 0; blue, 0 }  ,fill opacity=1 ] (488,212.75) .. controls (488,211.78) and (488.78,211) .. (489.75,211) .. controls (490.72,211) and (491.5,211.78) .. (491.5,212.75) .. controls (491.5,213.72) and (490.72,214.5) .. (489.75,214.5) .. controls (488.78,214.5) and (488,213.72) .. (488,212.75) -- cycle ;
	\draw  [fill={rgb, 255:red, 0; green, 0; blue, 0 }  ,fill opacity=1 ] (550,188.25) .. controls (550,187.28) and (550.78,186.5) .. (551.75,186.5) .. controls (552.72,186.5) and (553.5,187.28) .. (553.5,188.25) .. controls (553.5,189.22) and (552.72,190) .. (551.75,190) .. controls (550.78,190) and (550,189.22) .. (550,188.25) -- cycle ;
	\draw  [fill={rgb, 255:red, 0; green, 0; blue, 0 }  ,fill opacity=1 ] (22.5,203.75) .. controls (22.5,202.78) and (23.28,202) .. (24.25,202) .. controls (25.22,202) and (26,202.78) .. (26,203.75) .. controls (26,204.72) and (25.22,205.5) .. (24.25,205.5) .. controls (23.28,205.5) and (22.5,204.72) .. (22.5,203.75) -- cycle ;
	\draw  [fill={rgb, 255:red, 0; green, 0; blue, 0 }  ,fill opacity=1 ] (128,202.75) .. controls (128,201.78) and (128.78,201) .. (129.75,201) .. controls (130.72,201) and (131.5,201.78) .. (131.5,202.75) .. controls (131.5,203.72) and (130.72,204.5) .. (129.75,204.5) .. controls (128.78,204.5) and (128,203.72) .. (128,202.75) -- cycle ;
	\draw  [fill={rgb, 255:red, 0; green, 0; blue, 0 }  ,fill opacity=1 ] (232.5,204.25) .. controls (232.5,203.28) and (233.28,202.5) .. (234.25,202.5) .. controls (235.22,202.5) and (236,203.28) .. (236,204.25) .. controls (236,205.22) and (235.22,206) .. (234.25,206) .. controls (233.28,206) and (232.5,205.22) .. (232.5,204.25) -- cycle ;
	\draw  [fill={rgb, 255:red, 0; green, 0; blue, 0 }  ,fill opacity=1 ]
	 (294.5,209.75) .. controls (294.5,208.78) and (295.28,208) .. (296.25,208) .. controls (297.22,208) and (298,208.78) .. (298,209.75) .. controls (298,210.72) and (297.22,211.5) .. (296.25,211.5) .. controls (295.28,211.5) and (294.5,210.72) .. (294.5,209.75) -- cycle ;
	\draw  [fill={rgb, 255:red, 0; green, 0; blue, 0 }  ,fill opacity=1 ] (254,159.25) .. controls (254,158.28) and (254.78,157.5) .. (255.75,157.5) .. controls (256.72,157.5) and (257.5,158.28) .. (257.5,159.25) .. controls (257.5,160.22) and (256.72,161) .. (255.75,161) .. controls (254.78,161) and (254,160.22) .. (254,159.25) -- cycle ;
	\draw  [fill={rgb, 255:red, 0; green, 0; blue, 0 }  ,fill opacity=1 ] (266.5,165.75) .. controls (266.5,164.78) and (267.28,164) .. (268.25,164) .. controls (269.22,164) and (270,164.78) .. (270,165.75) .. controls (270,166.72) and (269.22,167.5) .. (268.25,167.5) .. controls (267.28,167.5) and (266.5,166.72) .. (266.5,165.75) -- cycle ;
	\draw  [fill={rgb, 255:red, 0; green, 0; blue, 0 }  ,fill opacity=1 ] (596,156.75) .. controls (596,155.78) and (596.78,155) .. (597.75,155) .. controls (598.72,155) and (599.5,155.78) .. (599.5,156.75) .. controls (599.5,157.72) and (598.72,158.5) .. (597.75,158.5) .. controls (596.78,158.5) and (596,157.72) .. (596,156.75) -- cycle ;
	\draw  [fill={rgb, 255:red, 0; green, 0; blue, 0 }  ,fill opacity=1 ] (614,155.75) .. controls (614,154.78) and (614.78,154) .. (615.75,154) .. controls (616.72,154) and (617.5,154.78) .. (617.5,155.75) .. controls (617.5,156.72) and (616.72,157.5) .. (615.75,157.5) .. controls (614.78,157.5) and (614,156.72) .. (614,155.75) -- cycle ;
	\draw (145,70.5) node [anchor=north west][inner sep=0.75pt]   [align=left] {self};
	\draw (473,62) node [anchor=north west][inner sep=0.75pt]   [align=left] {mixed};
	\draw (338,186.9) node [anchor=north west][inner sep=0.75pt]    {$p_{a}$};
	\draw (440.5,185.9) node [anchor=north west][inner sep=0.75pt]    {$p_{a}$};
	\draw (540.5,185.9) node [anchor=north west][inner sep=0.75pt]  [color={rgb, 255:red, 208; green, 2; blue, 27 }  ,opacity=1 ]  {$p_{a}$};
	\draw (385.5,214.4) node [anchor=north west][inner sep=0.75pt]    {$p_{b}$};
	\draw (478,212.9) node [anchor=north west][inner sep=0.75pt]    {$p_{b}$};
	\draw (14.5,204.9) node [anchor=north west][inner sep=0.75pt]    {$p_{a}$};
	\draw (116.5,204.9) node [anchor=north west][inner sep=0.75pt]    {$p_{a}$};
	\draw (218,208.9) node [anchor=north west][inner sep=0.75pt]  [color={rgb, 255:red, 208; green, 2; blue, 27 }  ,opacity=1 ]  {$p_{a}$};
	\draw (290,214.9) node [anchor=north west][inner sep=0.75pt]  [color={rgb, 255:red, 208; green, 2; blue, 27 }  ,opacity=1 ]  {$p'_{b}$};
	\draw (235.5,151.4) node [anchor=north west][inner sep=0.75pt]  [color={rgb, 255:red, 0; green, 0; blue, 255 }  ,opacity=1 ]  {$r_{1}$};
	\draw (271,148.4) node [anchor=north west][inner sep=0.75pt]  [color={rgb, 255:red, 0; green, 0; blue, 255 }  ,opacity=1 ]  {$r_{2}$};
	\draw (570,144.9) node [anchor=north west][inner sep=0.75pt]  [color={rgb, 255:red, 0; green, 0; blue, 255 }  ,opacity=1 ]  {$r_{1}$};
	\draw (635.5,142.4) node [anchor=north west][inner sep=0.75pt]  [color={rgb, 255:red, 0; green, 0; blue, 255 }  ,opacity=1 ]  {$r_{2}$};
\pgfresetboundingbox
\path[use as bounding box] (8,60) rectangle (634,236);
\end{tikzpicture}}

\title[Real monopole Floer homology and the Alexander polynomial]{Real monopole Floer homology\\and the Alexander polynomial}
\begin{document}
\begin{abstract}
	We prove that the Euler characteristic of the framed real monopole Floer homology of the double branched cover of a link is given by the Alexander polynomial evaluated at $\sqrt{-1}$.
	For knots $K$ with determinant one,
    this formula shows that Miyazawa's degree invariant of $K$ can be computed in terms of the Alexander polynomial of $K$.
	The proof involves a new connected-sum formula and oriented skein exact triangles for some new variants of real monopole Floer homology, called the diagonal and cone real monopole Floer homology, which admit mod-two gradings.
	This is the monopole counterpart of Srivastava's Euler characteristic formula in real Heegaard Floer homology.
\end{abstract}

\author[Jiakai Li]{Jiakai Li}
\address{Department of Mathematics, Stanford University\\
	450 Jane Stanford Way, Building 380, Stanford, CA 94305, USA}
\email{lijk@stanford.edu}
\maketitle

\section{Introduction}
In 2023, Miyazawa~\cite{Miyazawa} famously showed that there exist smoothly knotted $\mathbb{RP}^2$'s in the 4-sphere that are topologically standard.
This result is ground-breaking, as all previously known examples of exotic surfaces have much larger genera.
In distinguishing these knotted $\mathbb{RP}^2$'s from the standard one,
Miyazawa introduced an invariant of an embedded surface $\Sigma$ called degree $|\deg(\Sigma)|$.
This is the mapping degree of the real Bauer--Furuta invariant~\cite{KMT:Ktheory}, which counts framed real Seiberg--Witten solutions on branched double covers of $\Sigma$.
Via cut-and-paste techniques, the computation reduces to a 3-dimensional framed real Seiberg--Witten invariant of knots, called the degree $|\deg(K)|$.
Since then, it has been asked whether this $\mathbb{Z}/\pm$-valued knot invariant is classical.
That is, can $|\deg(K)|$ be derived from well-known invariants of knots such as the Alexander module?
Very few computations of $|\deg(K)|$ have been known, except for a few families of knots.
Therefore, a better understanding of $|\deg(K)|$ is fundamental to more computations of $|\deg(\Sigma)|$,
which can be used to distinguish exotic surfaces.

The degree of a determinant-one knot $K$ can be interpreted~\cite{Miyazawa, KMT:Ktheory,Xiao,Li26} as the Euler characteristic of the framed real monopole Floer homology (\hmr)~\cite{Li24,Li26} of the double branched cover along $K$.
An early attempt of the author to compute the Euler characteristic was to develop an unoriented skein exact triangle~\cite{Li24}.
While this approach shows interesting connections to Khovanov homology, the maps in the triangles are not mod-two homogeneous.
In fact, the vanilla versions of the real monopole Floer homology~\cite{Li22,Li24} may not even admit~\cite{Li26} relative mod-two gradings.
Furthermore, there is an Euler characteristic for each real spin\textsuperscript{c} structure on $\dbcv_2(K)$, each contributing a $\pm$ sign ambiguity.
Thus adding Euler characteristics for distinct real spin\textsuperscript{c} structures is not well-defined.
Hence the problem of interpreting the degree of a link is not well-posed for links in general, without an absolute mod-two grading on the real monopole Floer homology.

Recent developments in the Heegaard Floer homology have shed new light on this problem. Guth and Manolescu~\cite{GM} defined real Heegaard Floer homology (\hfr),
including a hat version, which admits relative mod-two gradings.
Their construction requires possibly multiple basepoints in order to make the relevant Lagrangians orientable.
At that time, no analogue existed on the monopole side.
More recently, Srivastava~\cite{Srivastava} fixed the sign of the total Euler characteristic using Seifert surface data and
proved a formula for the Euler characteristic of $\hfr$ of the branched double cover $\dbcv_2(L)$ of a link $L$, in terms of the multivariable Alexander polynomial evaluated at $i=\sqrt{-1}$.
The same formula was expected from \hmr, by the general principle~\cite{KutluhanLeeTaubes2010HFHM1,ColinGhigginiHonda2012HFECH1} that any flavour of monopole Floer homology~\cite{KM} should have a counterpart in Heegaard Floer theory~\cite{OzSz}.

Inspired by real Heegaard Floer homology, the author defined a framed version~\cite{Li26} of the real monopole Floer homology with multiple basepoints.
The construction involves suitable finite covers of configuration spaces.
In this setting, the author found a criterion for the existence of a relative mod-two grading.
In \cite{LM}, Manolescu and the author constructed absolute mod-two gradings $\widehat{\gr}$ on the framed real monopole Floer homology.
This fixes the sign of Miyazawa's degree invariant $\deg(K) \in \mathbb{Z}$,
and the input data of the absolute gradings parallel those of~\cite{Srivastava}.
More importantly, the cobordism maps in the framed real monopole homology have well-defined mod-two degrees.
Based on computational evidence, \cite{LM} conjectures that there should be a mod-two grading-preserving isomorphism between framed real monopole Floer homology and the hat-flavour of real Heegaard Floer homology.

As the final installment in the series \cite{Li26} and \cite{LM}, this paper establishes the real monopole Euler characteristic formula:
\begin{thm}\label{thm:alexander-intro}
	For every oriented link $L\subset S^3$,
	\begin{equation}
	\label{eq:alexander-result-intro}
		\widehat\chi\left(\tH(\dbcv_2(L),\mathbf p)\right)=(-i\sqrt2)^{|L|-1}\Delta_L(i).
	\end{equation}
	Here, $\mathbf{p}$ is a set of basepoints for which every component of $L$ contains one basepoint.
\end{thm}

The strategy of the proof involves a skein exact triangle, analogous to e.g.~\cite{KMalex,OzSz2008multi}.
Specifically,
the maps in the exact triangle are homogeneous with respect to the mod-two gradings.
Thus taking the Euler characteristics of the triangles gives rise to skein relations.
We emphasize that the right side of Equation~\eqref{eq:alexander-result-intro} is a number rather than a polynomial.
Indeed, real monopole Floer homology does not have a full Alexander grading/filtration like the knot Heegaard/monopole/instanton homologies.

This leaves the question of interpreting the formula: Recall that the Alexander polynomial of a knot can be defined algebraically using the infinite cyclic cover of the knot complement.
On the other hand, the real monopole Floer homology uses the order-two branched cover.
Does there exist a skein-free proof of Theorem~\ref{thm:alexander-intro}?
What happens if the order-two branched cover is replaced with an order-four cover?

The Euler characteristic of $\widetilde{\hmr}$ is related to the 4-dimensional framed Seiberg--Witten invariant, as observed in~\cite{Miyazawa}.
Broadly speaking, there are two layers to the problems of orientation for enumerative invariants: orientability and (canonical) choice of orientation.
The first allows a signed count to be defined over $\mathbb{Z}/\pm$, as in the case of the degree invariants.
The second upgrades the count from $\mathbb{Z}/\pm$ to $\mathbb{Z}$.
From this point of view, \cite{Li26} resolves the orientability problem and \cite{LM} resolves the orientation problem, for the framed Seiberg--Witten invariants on $S^1 \times Y$,
whenever the hypotheses are met.

Combining Theorem~\ref{thm:alexander-intro} with \cite[Proposition~4.30]{Miyazawa} and fixing their signs using $\widehat{\chi}$,
we have the following corollary.

\begin{cor}
	For every determinant-one knot $K$,
	\begin{equation*}
	\deg(K) = \deg \left(\BF_R(S^1 \times \dbcv_2(K))\right) = \widetilde{\swr}(S^1 \times \dbcv_2(K)) =  \widehat\chi\left(\tH(\dbcv_2(K),\mathbf p)\right) = \Delta_K(i).
	\end{equation*}
Here, $\BF_R$ is the real Bauer--Furuta invariant~\cite{KMT:homology} with sign fixed by $\widehat{\gr}$,
and $\widetilde{\swr}$ is the framed real Seiberg--Witten invariant~\cite{Miyazawa,Baraglia:surface,Li26} with orientation fixed by $\widehat{\gr}$.
On $S^1 \times \dbcv_2(K)$, we use the pullback real spin\textsuperscript{c} structure from $\dbcv_2(K)$.
\end{cor}

\begin{rem}\label{rem:polynomial-conventions}
Our convention for the Alexander polynomial of links is
$\Delta_L(t)=\nabla_L\!\left(t^{1/2}-t^{-1/2}\right)$,
where $\nabla_L$ is the Conway polynomial satisfying
$\nabla_U(z)=1$ and $\nabla_{L_+}(z)-\nabla_{L_-}(z)=z\nabla_{L_0}(z)$.
For $|L|>1$, with $\doteq$ denoting up to multiplication by a unit (\cite[\S2]{Cimasoni} and~\cite{EN85}),
the multivariable Alexander polynomial satisfies
\[\Delta_L(t)\doteq(t^{1/2}-t^{-1/2})\Delta_L^{\mathrm{multi}}(t,\ldots,t).\]
\end{rem}

Indeed, the skein relations of the Euler characteristic appearing in the proof of Theorem~\ref{thm:alexander-intro} coincide with those of~\cite[Proposition~1.6]{Srivastava}.
We conclude that the Euler characteristics of the real monopole and real Heegaard Floer theories coincide:
\begin{cor}
	For every link $L \subset S^3$,
	\begin{equation*}
		\widehat\chi\left(\tH(\dbcv_2(L),\mathbf p)\right) =
		\widehat\chi\left(\widehat{\hfr}(\dbcv_2(L),\mathbf p)\right).
	\end{equation*}
\end{cor}

This serves as further evidence for \cite[Conjecture~1.3]{LM} that $\tH$ and $\widehat{\hfr}$ are isomorphic in the mod-two graded sense.

In what follows, we summarize the three results used in the proof of Theorem~\ref{thm:alexander-intro}.

\subsubsection*{Diagonal and cone HMR}
By working with a finite-index subgroup of the real gauge group, we define the \emph{diagonal} real monopole Floer homologies.
This comes in three flavours, \emph{to}, \emph{from}, and \emph{bar}:
\[\widecheck{\HMR}_{\Delta,\bullet}(Y,\tau,\mathbf{p}),\quad \widehat{\HMR}_{\Delta,\bullet}(Y,\tau,\mathbf{p}), \quad \overline{\HMR}_{\Delta,\bullet}(Y,\tau,\mathbf{p}).\]
This notation parallels that of ordinary monopole Floer homology~\cite{KM}, real monopole Floer homology~\cite{Li22}, and Heegaard Floer homology~\cite{OzSz2004} ($\mathit{HF}^+,\mathit{HF}^-,\mathit{HF}^{\infty}$).

There are two notable new features.
First, recall that the vanilla real monopole Floer homologies are modules over $\F[\upsilon_1,\dots,\upsilon_n]/(\upsilon_i^2=\upsilon_j^2)$,
where each $\upsilon_i$ corresponds to
a point operation (V-map) $\upsilon_{p_i}$ for each $p_i\in\mathbf p$.\footnote{For convenience, assume every component contains one basepoint.}
Instead of multiple V-maps, the diagonal monopole Floer homology admits a \emph{single} V-map, denoted by $\upsilon$.

Second, under the same assumption for the absolute mod-two grading in the framed version,
\emph{all} flavours of the diagonal real monopole Floer homology admit an absolute grading $\widehat{\gr}$.
This grading $\widehat{\gr}$ is defined analogously its the framed counterpart~\cite{LM}.

By considering the mapping cone of $\upsilon$ on the \emph{to}-flavour of the diagonal \HMR, we define a fourth flavour, the \emph{cone} real monopole Floer homology of $(Y,\tau,\mathbf{p})$.
This is denoted analogously to the framed version:
\[\tH_\Delta(Y,\tau,\mathbf p).\]
This is a graded finite $\F$-vector space,
which inherits an absolute mod-two grading from the \emph{to}-$\HMR$.
In particular:
\begin{thm}\label{thm:cone-framed-euler}
The Euler characteristics of the framed and cone real monopole Floer homology coincide:
	\begin{equation}
	\label{eq:equal-chi-intro}
	\widehat{\chi}\left(\tH_\Delta(Y,\tau,\mathbf p)\right)
=\widehat{\chi}\left(\tH(Y,\tau,\mathbf p)\right).
\end{equation}
\end{thm}
The mapping cone construction (\cite{Bloom, Li24}) used in $\tH_{\Delta}$ is an algebraic realization of taking homology of a double cover.
For contrast, the cone version $\tH_{\Delta}$ is constructed using blown-up configuration spaces, whereas the framed version $\tH$ works directly with downstairs configuration spaces.

\subsubsection*{Exact triangles}
\begin{figure}[htbp]
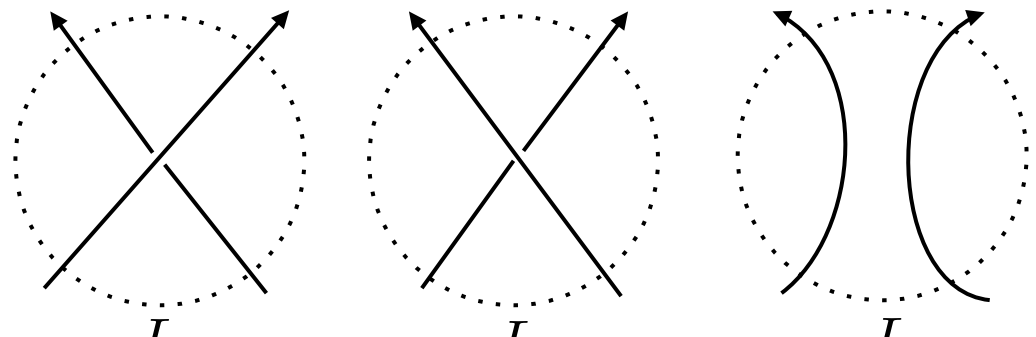

	\centering
	\figskeinintro
	\caption{An oriented skein triple of links.}
	\label{fig:skeinintro}
\end{figure}
An \emph{oriented skein triple} consists of three links $L_+,L_-,L_{0}$ that differ in a local 3-ball as in Figure~\ref{fig:skeinintro}.
\begin{thm}\label{thm:exact-triangle-intro}
	Given an oriented skein triple $(L_+,L_-,L_0)$,
	there exist cobordism maps $\widetilde{F}_+,\widetilde{F}_-,\widetilde{F}_{00}$ that satisfy an exact triangle
	\begin{equation}\label{eq:exact-triangle-intro}
		\begin{tikzcd}
			\tH_\Delta(\dbcv_2(L_+),\mathbf p)\arrow[rr,"\widetilde F_+"] &&
			\tH_\Delta(\dbcv_2(L_-),\mathbf p)\arrow[dl,"\widetilde F_{-}"]\\
			& \tH_\Delta(\dbcv_2(L_0) \# (S^1 \times S^2),\mathbf p)\arrow[ul,"\widetilde F_{00}"] &
		\end{tikzcd}.
	\end{equation}
Similarly, there is an exact triangle for the \emph{to}-version $\widecheck{\hmr}_{\Delta}$:
\begin{equation}
	\begin{tikzcd}
		\widecheck{\HMR}_\Delta(\dbcv_2(L_+),\mathbf p)\arrow[rr,"\widecheck F_+"] &&
		\widecheck{\HMR}_\Delta(\dbcv_2(L_-),\mathbf p)\arrow[dl,"\widecheck F_{-}"]\\
		& \widecheck{\HMR}_\Delta(\dbcv_2(L_0) \# (S^1 \times S^2),\mathbf p)\arrow[ul,"\widecheck F_{00}"] &
	\end{tikzcd}.
\end{equation}
Here, the basepoints lie away from the local balls and hence all vertices in the triangles inherit the same set of basepoints $\mathbf{p}$.
Finally, the maps $\widetilde{F}_+,\widetilde{F}_-,\widetilde{F}_{00},\widecheck{F}_+,\widecheck{F}_-,\widecheck{F}_{00}$ in the triangles above are all homogeneous with respect to the absolute mod-two grading.
\end{thm}
The advantage of working with a mapping cone is that the theory now allows mixed operations in the sense that elements of the ordinary monopole Floer homology can act on the real monopole Floer homology.
This allows us to apply the techniques in the ordinary monopole exact triangle to the real case.
The proof of Theorem~\ref{thm:exact-triangle-intro} also applies to the \emph{from} and \emph{bar} versions of the diagonal real monopole Floer homology.

\subsubsection*{Connected sum formula}
The third vertex $\dbcv_2(L_0) \# (S^1 \times S^2)$ in \eqref{eq:exact-triangle-intro} can be interpreted as the double connected sum $\dbcv_2(L_0) \#_2 S^3$ in the following sense.
Let $Y$ and $Z$ be two real 3-manifolds, with basepoint sets $\mathbf{p}$ and $\mathbf{q}$, respectively.
The set $\mathbf{q}$ is allowed to be empty.
Let $Y \#_2 Z$ be obtained by attaching two 1-handles from $x_i \in Y$ to $z_i \in Z$, for $i=1,2$.
Working directly with framed configurations avoids a full connected sum formula in the diagonal version.
We prove the following theorem in the framed real monopole homology.

\begin{thm}
	\label{thm:connected-sum-intro}
	Suppose $x_1,x_2$ both lie on components of the fixed set that contain basepoints in $\mathbf{p}$.
	Then there exists an isomorphism
	\begin{equation*}
		\widetilde{\HMR}(Y \#_2 Z, \mathbf{p} \cup \mathbf{q}) \cong
		\widetilde{\HMR}(Y, \mathbf{p}) \otimes \widetilde{\HMR}(Z, \mathbf{q} \cup \{z_1,z_2\}).
	\end{equation*}
\end{thm}

The mod-two degrees of the isomorphism can be computed on a case-by-case basis.
Applying Theorem~\ref{thm:connected-sum-intro} to $Y=\dbcv_2(L_0)$ and $Z = S^3$ with suitable choices of basepoints gives a relation between the Euler characteristics of the $\tH$ of $\dbcv_2(L_0)$ and $\dbcv_2(L_0) \#_2 S^3$.
Here the relation depends on whether the zeroth oriented resolution merges or splits the components of $L_+$.
Via the equality~\eqref{eq:equal-chi-intro} of the framed and cone Euler characteristics, Theorem~\ref{thm:exact-triangle-intro} implies that $\widehat{\chi}(\tH(\dbcv_2(L)))$ satisfies the same recursive formula as that of $(-i\sqrt2)^{|L|-1}\Delta_L(i)$.

\subsubsection*{Speculation about pointed Khovanov homology}
The author proves a spectral sequence~\cite{Li24} from the reduced Khovanov homology to the single-pointed real monopole Floer homology of a link.
The appearance of the cube of chain complexes in \S\ref{ssec:V-map}
\begin{equation*}
	\bigoplus_{\mathsf{v} \in \{0,1\}^{\mathbf{p}}} \widecheck{\cmr}_{\Delta,\bullet}(\dbcv_2(L), \mathbf{p}; \mathsf{v}\widetilde{})
\end{equation*}
suggests that the cone real monopole Floer homology should be related to the (perhaps reduced) pointed Khovanov homology~\cite{BLS}.
\begin{conjecture}
	\label{conj:ptKh}
	There is a spectral sequence from the pointed Khovanov homology $\mathit{Kh}(-L,\mathbf{p})$ of $-L$ to $\thmr_{\Delta}(\dbcv_2(L),\mathbf{p})$.
\end{conjecture}

The analytical foundation of the diagonal real monopole Floer homology is adapted from~\cite{Li22}, which builds upon \cite{KM}.
The foundation of the cone-version is adapted from its single basepoint case, treated in~\cite{Li24}, which in turn builds upon \cite{KM} and \cite{Bloom} in ordinary monopole theory.
In this article, we reduce the proof of the oriented exact triangles to the ordinary case \cite{KMOS} and the framed ordinary case \cite{Bloom}, using the \textsc{r}eal adaptation from~\cite{Li24}.

\subsection{Organization of the paper}
Section~\ref{sec:diagonal-hmr} defines the diagonal \HMR\ and its absolute mod-two gradings.
Section~\ref{sec:cone-hmr} defines the cone \hmr\
 and establishes the equality~\eqref{eq:equal-chi-intro} of Euler characteristics of the cone and framed real \hmr\ (Theorem~\ref{thm:cone-framed-euler}).
Section~\ref{sec:exact-triangle} proves oriented skein exact triangles (Theorem~\ref{thm:exact-triangle-intro}) in the diagonal and cone \hmr.
Section~\ref{sec:connected-sum} gives a connected sum formula (Theorem~\ref{thm:connected-sum-intro}) for the framed \hmr.
Finally, Section~\ref{sec:alex} deduces the Euler characteristic formula (Theorem~\ref{thm:alexander-intro}) for cone and framed \HMR, using the oriented skein exact triangle
in combination with Equation~\eqref{eq:equal-chi-intro} and the connected sum theorem.

\subsection{Statement of AI use}
AI helped with finding the suitable convention (Remark~\ref{rem:polynomial-conventions}) for the Alexander polynomial of links, and tested the normalization formula against several families of examples.
In Section~\ref{ssec:double-connected-sum}, AI suggested to book-keep the \RRSC\ structures with the maps $i_{\Gamma}$ and $p_{\Gamma}$, which simplifies the initial exposition.
During manuscript preparation, AI checked the arithmetics and the consistency of signs, left/right, ascending/descending indices, super/subscripts, etc.
The figures were first made by the author with a TikZ editor (AI-free); then Codex helped with shading and small adjustment.
AI was used for checking for mathematical errors and proof-reading.

\subsection{Acknowledgement}
The author is indebted to Deeparaj Bhat and Hongjian Yang for many discussions on a related project.
The author also thanks Ciprian Manolescu for comments on the draft.
Part of the research was conducted while the author was visiting the Okinawa Institute of Science and Technology (OIST) through the Thematic Program and Symposium on ``New Frontiers in Gauge Theory, Topology, and Physics'' (TP26GT).

\section{Diagonal real monopole Floer homology}\label{sec:diagonal-hmr}

The diagonal real monopole Floer homology uses the relative real structures of \cite[\S3]{Li26} as the multi-framed real monopole Floer homology.
However, we use larger gauge groups containing $-1$ and the blown-up construction of \cite{KM}.

\subsection{Relative real spin$^c$ structures and blown-up configuration spaces}
Let $Y$ be a closed connected oriented 3-manifold,
with an orientation-preserving involution $\tau$.
Assume $C=\fix(\tau)$ is non-empty.
Let $\mathbf{p} = \{p_1, \dots, p_n\} $ be a nonempty set of basepoints in $C$.
At each $p\in\mathbf p$, fix a standard Clifford module $(S_p^0,\rho_p^0)$ and a real structure $\sr_p^0$, using the local models of \cite[\S3.2]{Li26} in terms of Pauli matrices and complex conjugation.

\begin{defn}
A relative real spin$^c$ ($\RRSC$) structure is a triple $\widetilde{\sr}=(\mathfrak s,T,\sr)$,
consisting of a spin$^c$ structure $\mathfrak{s}$,
Clifford-module identifications $T$ at the points of $\mathbf{p}$,
and a compatible conjugate-linear lift $\sr$ with $\sr^2=1$.
Here, $T_p:S_p\to S_p^0$ is a unitary Clifford-module identification for which
$T_p\sr=\sr_p^0T_p$ for every $p\in\mathbf p$.
An \emph{isomorphism} of \rrsc\ structures $u \colon (S,T,\sr)\to(S',T',\sr')$ is a
spin\textsuperscript{c} isomorphism satisfying $u\sr=\sr'u$ and
$T'_pu_p=T_p$ over the basepoints $p \in \mathbf{p}$.
Denote by $\RRSC(Y,\tau,\mathbf{p})$ the set of isomorphism classes of $\RRSC$ structures.
\end{defn}

The \emph{(unframed) real gauge group} $\mathcal{G}$ consists of $u \colon Y \to S^1$ anti-invariant under $\tau$, i.e. $u\circ\tau=\bar u$.
Every element of $\mathcal{G}$ takes the values $\pm1$ on the fixed set and is constant on each fixed component.
The \emph{framed real gauge group} is $\mathcal{G}_{\mathbf{p}}=\{u\in\mathcal{G}:u|_{\mathbf{p}}=1\}$, consisting of real gauge transformations that evaluate to $+1$ on $\mathbf{p}$.
\begin{defn}
	The \emph{diagonal real gauge group} is the subgroup of $\mathcal{G}$ generated by $\mathcal{G}_{\mathbf{p}}$ and the constant $-1$.
	Equivalently,
	\[
	 \mathcal{G}_\Delta=
	 \left\{ u \in \mathcal{G} : u(p_1) = \cdots = u(p_n) \, \right\}.
 \]
\end{defn}

Let $\mathcal A^{\sr} = \mathcal A(Y,\widetilde{\sr})$ be the affine space of spin$^c$ connections compatible with the real structure $\sr$, modelled on $\tau$-anti-invariant $1$-forms.
The downstairs configuration space is $\mathcal{C}=\mathcal A^{\sr}\times\Gamma(S)^{\sr}$,
and the blow-up configuration space is
\[
	\mathcal{C}^\sigma:=\{(B,s,\phi): B\in\mathcal A^{\sr},\  s\geq0,\ \sr\phi=\phi, \ \|\phi\|_{L^2(Y)}=1\},
\]
and $\mathcal{B}_\Delta^\sigma:=\mathcal{C}^\sigma/\mathcal{G}_\Delta$.
The blow-down map $\pi \colon \mathcal{C}^{\sigma} \to \mathcal{C}$ sends $(B,s,\phi)$ to $(B,s\phi)$.
We work with $L^2_k$-completion of configuration spaces and $L^2_{k+1}$-gauge groups, for a sufficiently large $k$,
as in \cite[\S5--6]{Li22}. The Sobolev subscripts will be suppressed whenever no confusion arises.

Let $\ev_{\mathbf{p}}: \mathcal{G} \to \{\pm 1\}^{\mathbf{p}}$ be the homomorphism $g \mapsto (g(p_k))$ that evaluates every $g \in \mathcal{G}$ over $\mathbf{p}$,
and $\epsilon_{\mathbf{p}}$ be the homomorphism $\mathcal{G}_{\Delta} \to \{\pm 1\}$ that takes the common value of $g$ over $\mathbf{p}$.
In particular, $\mathcal{G}_\Delta$ is the inverse image of the \emph{diagonal} subgroup,
and $\mathcal{G}_{\mathbf{p}}$ is the kernel of the evaluation $\ev_{\mathbf{p}}$.
There are three versions of configuration spaces of gauge-equivalence classes,
arising from the nested gauge groups $\mathcal{G}_{\mathbf{p}} \le \mathcal{G}_{\Delta} \le \mathcal{G}$:
\[
\mathcal{B}_{\mathbf p}^\sigma=\mathcal{C}^\sigma/\mathcal{G}_{\mathbf p},\qquad
\mathcal{B}_\Delta^\sigma=\mathcal{C}^\sigma/\mathcal{G}_\Delta,\qquad
\mathcal{B}^\sigma=\mathcal{C}^\sigma/\mathcal{G}.
\]
Thus these configuration spaces fit in the following sequence of finite coverings:
\[
\mathcal{B}_{\mathbf p}^\sigma\longrightarrow\mathcal{B}_\Delta^\sigma
\longrightarrow\mathcal{B}^\sigma
\]
with deck transformations $\{\pm1\}$ and
$\im (\ev_{\mathbf{p}})/\langle(-1,\ldots,-1)\rangle$, respectively.
The double cover $\mathcal{B}_{\mathbf p}^\sigma \to \mathcal{B}^{\sigma}_{\Delta}$ defines a real line bundle and a cohomology class
\begin{equation*}
	\mathbb L_{\mathbf p}
	:=\mathcal{C}^\sigma\times_{\epsilon_{\mathbf p}}\mathbb{R}
	\longrightarrow\mathcal{B}_\Delta^\sigma,
	\qquad
	\upsilon_{\mathbf p}:=w_1(\mathbb L_{\mathbf p})
	\in H^1(\mathcal{B}_\Delta^\sigma;\F).
\end{equation*}
Finally, for the downstairs configurations, let $\mathcal{B}=\mathcal{C}/\mathcal{G}$, $\mathcal{B}_{\Delta} = \mathcal{C}/\mathcal{G}_{\Delta}$, and $\widetilde{\mathcal{B}}_{\mathbf{p}} = \mathcal{C}/\mathcal{G}_{\mathbf{p}}$.
These also fit into a sequence:
\begin{equation*}
	\widetilde{\mathcal{B}}_{\mathbf{p}} \to \mathcal{B}_{\Delta} \to \mathcal{B}.
\end{equation*}
The first space $\widetilde{\mathcal{B}}_{\mathbf{p}}$ is the framed real configuration space used in \cite{Li26} for the definition of the framed \HMR.
Furthermore, the map $\widetilde{\mathcal{B}}_{\mathbf{p}} \to \mathcal{B}_{\Delta}$ is 2-to-1 away from the reducible configurations, and 1-to-1 over the reducibles.

\begin{rem}
	The class $\upsilon_{\mathbf{p}}$ is the analogue of the $\upsilon_p$ in \cite{Li22}.
In the diagonal setting, the evaluation at every basepoint $p \in \mathbf{p}$ gives the same $\upsilon$-map.
Henceforth, if the basepoint set $\mathbf{p}$ is understood, we write $\upsilon$ in place of $\upsilon_{\mathbf{p}}$.
\end{rem}

\subsection{Perturbations and trajectory spaces}
Fix a \rrsc\ structure $(\mathfrak{s}, T, \sr)$.
To set up the moduli spaces of trajectories in $\mathcal{B}^{\sigma}_{\Delta}$,
we use the \emph{tame real perturbations} as in \cite[\S6]{Li22}.
It suffices to use tame perturbations that are invariant under the full \emph{unframed real gauge group $\mathcal{G}$},
and pull back the perturbations from $\mathcal{B}^{\sigma}$ to $\mathcal{B}_\Delta^\sigma$.
Let $\pertL \colon \mathcal{C}\to \mathbb{R}$ denote the perturbed Chern--Simons--Dirac (CSD) functional.
Thus the perturbed gradient vector field on the diagonal quotient space $\mathcal{B}_{\Delta}^{\sigma}$ is the lift of the usual blown-up gradient $\grad \pertL_{\mathfrak{q}}^{\sigma}$ on $\mathcal{B}^{\sigma}$.
In particular, all such perturbations are invariant under $\mathbb{Z}_2 \cong \{\pm1\}$.

Fix a tame perturbation $\mathfrak q$ whose critical points are nondegenerate.
The \emph{interior critical points} are the irreducible critical points.
The \emph{boundary-stable} and \emph{boundary-unstable} critical points are reducible critical points;
the former are of the form $(B,0,\psi)$ for which $\psi$ is a positive eigenvector of the perturbed Dirac operator $D_{B,\mathfrak{q}}$,
and the latter correspond to negative-eigenvalue eigenvector of $D_{B,\mathfrak{q}}$.
Let $\mathfrak C_\Delta^o$, $\mathfrak C_\Delta^s$ and $\mathfrak C_\Delta^u$ denote the interior, boundary-stable and boundary-unstable critical points in $\mathcal{B}_\Delta^\sigma$.

On the infinite cylinder $\mathbb{R} \times Y$, we use the Sobolev completions and asymptotic conditions as in~\cite[\S8]{Li22}.
For critical points $\mathfrak a,\mathfrak b$, let $\pi_1(\mathcal{B}_\Delta^\sigma;\mathfrak a,\mathfrak b)$ be the set of homotopy classes of paths from $\mathfrak a$ to $\mathfrak b$, with endpoints fixed.
Given such a class $z$, let
$M_{\Delta,z}(\mathfrak a,\mathfrak b)$ be the moduli space
of solutions to the perturbed Seiberg--Witten equations $\mathfrak{F}_{\mathfrak{q}}^{\sigma}$ on $\mathbb{R} \times Y$
with limit $\mathfrak a$ at $-\infty$ and limit $\mathfrak b$ at $+\infty$
in the class $z$.
Denote by $\breve M_{\Delta,z}(\mathfrak a,\mathfrak b)$ the space of nonconstant trajectories modulo translation,
and by $M^{\red}_{\Delta,z}$ the subspace of trajectories contained in the reducible locus $\partial\mathcal{B}_\Delta^\sigma$.
Let $\gr_z(\mathfrak{a}, \mathfrak{b})$ be the relative grading between $\mathfrak{a}$ and $\mathfrak{b}$.

Given a trajectory $\upgamma$,
let $Q_\upgamma$ be the linearization of the Seiberg--Witten operator with gauge-fixing on $\mathbb{R}\times Y$.
At a reducible trajectory, let $Q_\gamma^\partial$ denote the linearization of the (gauge-fixed) Seiberg--Witten operator restricted to the reducible boundary.
A trajectory is \emph{regular} if $Q_\gamma$ is surjective,
except when $\mathfrak a\in\mathfrak C_\Delta^s$ and $\mathfrak b\in\mathfrak C_\Delta^u$.
In this case, the trajectory $\upgamma$ is \emph{boundary-obstructed} and $\upgamma$ is \emph{regular} if and only if $Q_\gamma^\partial$ is surjective.
A trajectory space is \emph{regular} if each of its trajectories is regular.
See~\cite[Definition~8.20]{Li22}.

\begin{defn}
A perturbation $\mathfrak{q}$ is \emph{admissible} if all the critical points of $\grad \pertL_{\mathfrak{q}}$ are nondegenerate,
all corresponding trajectory spaces are regular,
and no reducible critical points exist when the underlying spin$^c$ structure $\mathfrak{s}$ is non-torsion.
\end{defn}

\begin{prop}\label{prop:diagonal-analysis}
Let $\mathfrak q$ be a tame perturbation for the non-diagonal configuration space $\mathcal{B}^{\sigma}$
for which all the corresponding critical points are nondegenerate and all trajectory spaces $M_z(\mathfrak{a},\mathfrak{b})$ are regular, as provided by~\cite[Theorem~8.24]{Li22}.
Then the following hold for the pullback of $\mathfrak{q}$ to $\mathcal{B}_\Delta^\sigma$, and the corresponding diagonal moduli spaces.

\begin{enumerate}[leftmargin=*, label=(\roman*)]
\item
The spaces $M_{\Delta,z}(\mathfrak a,\mathfrak b)$ are regular in the sense of \cite[Definition~8.20]{Li22}.
The smooth structures of $M_{\Delta,z}(\mathfrak a,\mathfrak b)$ admit the same description as in~\cite[Proposition~8.21]{Li22}.

\item
Let $\breve M_{\Delta,z}^{+}(\mathfrak a,\mathfrak b)$ the moduli space obtained by adding broken trajectories.
Then $\breve M_{\Delta,z}^{+}(\mathfrak a,\mathfrak b)$ is compact, as in \cite[Theorem~9.3]{Li22}.
Moreover, as in \cite[Propositions~9.4--9.5]{Li22},
only finitely homotopy classes $z$ admit nontrivial trajectories for a given energy upper bound.

\item
The compactified trajectory spaces and compactified reducible trajectory spaces have the stratifications described in \cite[Propositions~9.9--9.11]{Li22}.
In particular, the classification of codimension-one strata apply verbatim.

\item
The neighbourhoods of broken trajectories in $\breve M_{\Delta,z}^{+}(\mathfrak a,\mathfrak b)$ are described by the gluing constructions of \cite[Theorems~10.7--10.8]{Li22}, for both boundary-obstructed and unobstructed trajectories.
The corresponding codimension-one boundary $\delta$-structures coincide with the description of~\cite[Theorems~10.11--10.12]{Li22}.
\end{enumerate}
\end{prop}

\begin{proof}
Since $\pi_{\Delta} \colon \mathcal{B}_\Delta^\sigma\to\mathcal{B}^\sigma$ is a finite covering,
the pulled-back gradient flow equations in $\mathcal{B}_{\Delta}^{\sigma}$ have the same linearizations as the downstairs ones,
including their normal components at reducibles.
Thus the regularity and the dimension statements in (i) follow from their counterparts on $\mathcal{B}^\sigma$.
By the path-lifting property, every (broken) path in $\mathcal{B}^{\sigma}$ can be lifted to $\mathcal{B}_{\Delta}^{\sigma}$.
The projection $\pi_{\Delta}$ induces homeomorphisms
\[
\breve M_{\Delta,z}^{+}(\mathfrak a,\mathfrak b)
\;\cong\;
\breve M_{{\pi_{\Delta}}_*z}^{+}
\bigl({\pi_{\Delta}}(\mathfrak a),{\pi_{\Delta}}(\mathfrak b)\bigr).
\]
Therefore, the compactness and finiteness in (ii) hold on the cover $\mathcal{B}_{\Delta}^{\sigma}$.
The gluing constructions~\cite[Theorems~10.7--10.8]{Li22} of $\mathcal{B}_{\Delta}^{\sigma}$ and $\mathcal{B}^{\sigma}$ can be performed with the same obstruction maps.
This implies the stratifications in (iii) and boundary structures in (iv).
\end{proof}

\subsection{Floer complexes}
Let $C^o,C^s,C^u$ be the $\F$-vector spaces generated by the critical point sets
$\mathfrak C_\Delta^o,\mathfrak C_\Delta^s,\mathfrak C_\Delta^u$.
Let
\[\partial_o^o \colon {C}_{\Delta}^o
\to {C}_{\Delta}^o,\quad
\partial_s^o \colon {C}_{\Delta}^o
\to {C}_{\Delta}^s,\quad
\partial_o^u \colon {C}_{\Delta}^u
\to {C}_{\Delta}^o,\quad
\partial_s^u \colon {C}_{\Delta}^u
\to {C}_{\Delta}^s.\]
count zero-dimensional unparametrized trajectories.
Similarly, let
\[\bar\partial_u^u \colon {C}_{\Delta}^u
\to {C}_{\Delta}^u,\quad
\bar\partial_s^u \colon {C}_{\Delta}^u
\to {C}_{\Delta}^s,\quad
\bar\partial_u^s \colon {C}_{\Delta}^s
\to {C}_{\Delta}^u, \quad
\bar\partial_s^s \colon {C}_{\Delta}^s
\to {C}_{\Delta}^s.
\]
count the corresponding spaces of reducible trajectories between critical points.
The diagonal \emph{to}, \emph{from} and \emph{bar} complexes are
\begin{equation*}
\widecheck C_\Delta=C^o\oplus C^s,\quad\widehat C_\Delta=C^o\oplus C^u, \quad \overline C_\Delta =C^s\oplus C^u,
\end{equation*}
equipped with the differentials\footnote{The minus signs in the matrices are cosmetic, as we work over $\F$.
From now on, all minus signs in the matrix definitions of maps will be made plus.}:
\begin{equation*}
	\check\partial=
	\begin{pmatrix}
		\partial_o^o&-\partial_o^u\bar\partial_u^s\\
		\partial_s^o&\bar\partial_s^s-\partial_s^u\bar\partial_u^s
	\end{pmatrix},
\quad
	\hat\partial=
	\begin{pmatrix}
		\partial_o^o&\partial_o^u\\
		-\bar\partial_u^s\partial_s^o&
		-\bar\partial_u^u-\bar\partial_u^s\partial_s^u
	\end{pmatrix},
	\quad
	\bar\partial=
	\begin{pmatrix}
		\bar\partial_s^s&\bar\partial_s^u\\
		\bar\partial_u^s&\bar\partial_u^u
	\end{pmatrix}.
\end{equation*}
Proposition~\ref{prop:diagonal-analysis} says the boundary structures of 1-dimensional trajectories in the diagonal
configurations have the same boundary stratification.
Thus the boundary identities in \cite[\S22.1]{KM} show that $\check\partial^2 = 0$, $\hat\partial^2 = 0$,
and $\bar\partial^2 = 0$.

\begin{defn}
Suppose $\widetilde{\sr}=(\mathfrak{s},T,\sr)$ is an \rrsc\ structure.
The homologies of $(\check C_\Delta, \check\partial)$, $(\hat C_\Delta, \hat\partial)$ and
$(\bar C_\Delta, \bar\partial)$ are the \emph{to-, from-, and bar-flavours} of the \emph{diagonal real monopole Floer homology}, 
denoted as
\[\widecheck{\HMR}_\Delta(Y,\tau,\mathbf p;\widetilde{\sr}), \quad \widehat{\HMR}_\Delta(Y,\tau,\mathbf p;\widetilde{\sr}), \quad \overline{\HMR}_\Delta(Y,\tau,\mathbf p;\widetilde{\sr}).\]
When $\widetilde{\sr}$ is dropped from the notation,
we take the direct sum over all \rrsc\ structures on $(Y,\tau,\mathbf{p})$.
\end{defn}

Each flavour of diagonal Floer homology admits a relative grading by a set $\mathbb{J}(\widetilde{\sr})$ as in \cite[\S11.2]{Li22},
but using paths in the cover $\mathcal{B}_\Delta^\sigma$.
As usual, these three flavours fit into a long exact sequence:
\[\cdots \xrightarrow{i_*}\widecheck{\HMR}_{\Delta,*}(Y,\tau,\mathbf{p})
\xrightarrow{j_*}\widehat{\HMR}_{\Delta,*}(Y,\tau,\mathbf{p})
\xrightarrow{p_*}\overline{\HMR}_{\Delta,*}(Y,\tau,\mathbf{p})
\xrightarrow{i_*}\widecheck{\HMR}_{\Delta,*}(Y,\tau,\mathbf{p})
\xrightarrow{j_*} \cdots.
\]

To describe cobordism maps, we complete the gradings in the negative direction as in~\cite[\S3.1]{KM}.
Denote the resulting groups as
\[\widecheck{\HMR}_{\Delta,\bullet}, \quad \widehat{\HMR}_{\Delta,\bullet}, \quad \overline{\HMR}_{\Delta,\bullet}.\]
From now on, we focus on $\widecheck{\HMR}_{\Delta,\bullet}$,
as the constructions in the from- and bar-flavours are analogous.

\subsection{Cobordism maps}
\begin{figure}[htbp]
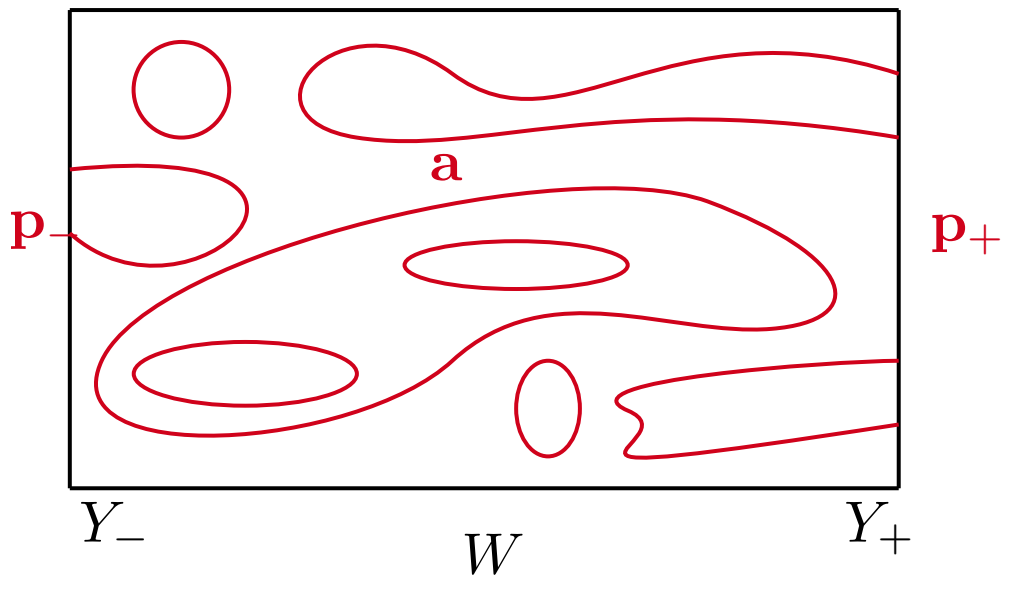

	\centering
	\figcobW
	\caption{A marked real cobordism $W$ from $(Y_-,\tau_-,\mathbf p_-)$ to $(Y_+,\tau_+,\mathbf p_+)$, with the marking set $\mathbf{a}$ coloured red.}
	\label{fig:cobW}
\end{figure}

We use the same cobordism data as in~\cite[\S3]{Li26}:
A \emph{marked real cobordism} is a triple
\[
(W,\tau_W,\mathbf a)\colon (Y_-,\tau_-,\mathbf p_-)
\longrightarrow(Y_+,\tau_+,\mathbf p_+)
\]
consisting of a cobordism $W\colon Y_- \to Y_+$, an involution $\tau_W$ on $W$, and a 1-dimensional subset $\mathbf{a} \subset \fix(\tau_W)$, which is a properly embedded 1-manifold
with boundary $\mathbf{p}_- \cup \mathbf{p}_+$ as in Figure~\ref{fig:cobW}.

Let the 4-dimensional framed and diagonal gauge groups be $\mathcal{G}_{\mathbf a}=\{u\in\mathcal{G}(W):u|_{\mathbf a}=1\}$
and $\mathcal{G}_{\Delta,\mathbf{a}}=\langle\mathcal{G}_{\mathbf a},-1\rangle$, respectively.
Let $\epsilon_{\mathbf{a}} \colon\mathcal{G}_{\Delta,\mathbf{a}} \to \{\pm 1\}$ be the homomorphism that evaluates a gauge transformation at the marking set $\mathbf{a}$.
Choose once and for all Clifford trivializations along $\mathbf a$,
and a relative real spin$^c$ structure $\widetilde{\sr}_W = (\mathfrak{t},T_{\mathbf{a}},\sr_W)$ on $W$,
that restricts to \rrsc\ structures $\widetilde{\sr}_{\pm}$ on $Y_{\pm}$.
After attaching cylindrical ends to $W$,
we consider perturbed Seiberg--Witten solutions on $W$ in the blown-up configuration modulo $\mathcal{G}_{\Delta}=\mathcal{G}_{\Delta,\mathbf{a}}$ (\cite[\S12]{Li22}).
There exists a residual set of perturbations for which the moduli spaces on $W$ are regular, by~\cite[Proposition~12.10]{Li22}.

Let $m^{*}_*$ and $\overline{m}^*_*$ denote the counts of zero-dimensional moduli spaces and reducible moduli spaces on $W$,
with the superscript and subscript denoting the input and output type.
Let $\partial_{\pm,*}^*$ denote the from-differential on the $Y_{\pm}$-end.
The cobordism map from $(Y_-; \widetilde{\sr}_-)$ to $(Y_+; \widetilde{\sr}_+)$ is given on the chain level by the matrix:
\begin{equation}\label{eq:diagonal-cobordism-map}
	\widecheck{m}_{\Delta}(W; \widetilde{\sr}_W) = \widecheck{\cmr}_{\Delta}(W; \widetilde{\sr}_W)=
	\begin{pmatrix}
		m_o^o&m_o^u\bar\partial_{-,u}^s+
		\partial_{+,o}^u\bar m_u^s\\
		m_s^o&\bar m_s^s+m_s^u\bar\partial_{-,u}^s+
		\partial_{+,s}^u\bar m_u^s
	\end{pmatrix}.
\end{equation}

The blownup moduli spaces on marked manifolds with boundary have the same $\delta$-structures and codimension-one boundary stratification as in \cite{Li22}.
By Proposition~\ref{prop:diagonal-analysis} and the same argument as in \cite[Proposition~22.1.4]{KM},
we have
	\begin{equation*}
	\widecheck{\partial}_+ \circ \widecheck{m}_{\Delta}(W; \widetilde{\sr}_W) +
	\widecheck{m}_{\Delta}(W; \widetilde{\sr}_W) \circ \widecheck\partial_- = 0.
	\end{equation*}
Hence we denote by
\[\widecheck{\hmr}_{\Delta}(W; \widetilde{\sr}_W) \colon\widecheck{\HMR}_{\Delta,\bullet}(Y_-; \widetilde{\sr}_-) \to \widecheck{\HMR}_{\Delta,\bullet}(Y_+; \widetilde{\sr}_+)\]
the induced map on homology,
and $\widecheck{\hmr}_{\Delta}(W)$ be the sum over all \rrsc\ structures.
	Furthermore, by Proposition~\ref{prop:diagonal-analysis} and the neck-stretching analysis \cite[Propositions~13.9 and
	13.11--13.13]{Li22}, we have:

\begin{prop}\label{prop:diagonal-functoriality}
	The chain-homotopy classes of the cobordism maps are independent of metrics and regular perturbations.
	Furthermore, the cobordism maps induced by the cylinders $I \times (Y,\tau, \mathbf{p})$ induce canonical isomorphisms on homology.
	Finally, given two marked real cobordisms $W_0 \colon Y_0 \to Y_1$ and $W_1 \colon Y_1 \to Y_2$, the composite cobordism map satisfies
	\begin{equation*}
		\widecheck{\hmr}_{\Delta}(W_1) \circ
		\widecheck{\hmr}_{\Delta}(W_0)
		=
		\widecheck{\hmr}_{\Delta}(W_1 \cup_{Y_1} W_0).
	\end{equation*}
\end{prop}

Fix $\F$-cohomology classes $u_0 \in H^*(\boldsymbol{\mathcal{B}}_{\Delta}^{\sigma}(W_0))$ and $u_1 \in H^*(\boldsymbol{\mathcal{B}}_{\Delta}^{\sigma}(W_1))$,
where $\boldsymbol{\mathcal{B}}_{\Delta}^{\sigma}(W_i)$ denotes diagonal blown-up configuration space on $W_i$ over all \RRSC\ structures and homotopy classes of $W$-paths.
Evaluation of $u_i$ can be coupled with cobordism maps as in \cite[\S25.3]{KM} to
obtain a chain map $\widecheck{m}_{\Delta}(W_i|u_i)$.
Denote by $u_1u_0$ the resulting class on $W_1 \cup W_0$.
Then one can prove the composition law (\cite[Proposition~13.13]{Li22}):
\begin{equation*}
		\widecheck{\hmr}_{\Delta}(u_1|W_1) \circ
		\widecheck{\hmr}_{\Delta}(u_0|W_0)
		=
		\widecheck{\hmr}_{\Delta}(u_1 u_0|W_1 \cup_{Y_1} W_0).
\end{equation*}

\subsection{The V-map}\label{ssec:V-map}
Given a marked real cobordism $(W,\tau_W,\mathbf{a})$,
evaluation of the gauge group $\mathcal{G}_{\Delta}$ along $\mathbf{a}$ restricts to the same evaluation on $\mathbf{p}_-$ and $\mathbf{p}_+$.
Let $\upsilon_{\mathbf{a}}$ be the corresponding class in $\boldsymbol{\mathcal{B}}_{\Delta}^{\sigma}(W)$.
This defines
\begin{equation*}
	\widecheck m(\upsilon_{\mathbf{a}}) =\widecheck{m}_{\Delta}(W, \mathbf{a}_W\, | \, \upsilon_{\mathbf{a}}) \colon \widecheck{\HMR}_{\Delta,\bullet}(Y_-, \mathbf{p}_-; \widetilde{\sr}_-) \to \widecheck{\HMR}_{\Delta,\bullet}(Y_+, \mathbf{p}_+; \widetilde{\sr}_+).
\end{equation*}
Fix $(Y,\tau,\mathbf{p})$,
and let $\upsilon = \upsilon_{\mathbf{p}}$ be the corresponding classes in $\boldsymbol{\mathcal{B}}_{\Delta}^{\sigma}(Y)$.
The product cobordism induces a chain endomorphism:
\begin{equation*}
	\widecheck m(\upsilon) = \widecheck{m}_{\Delta}(I \times Y, I \times \mathbf{p}\, | \, \upsilon_{I \times \mathbf{p}}) \colon \widecheck{\HMR}_{\Delta,\bullet}(Y; \widetilde{\sr}) \to \widecheck{\HMR}_{\Delta,\bullet}(Y; \widetilde{\sr}),
\end{equation*}
of degree $-1$.
We also denote this chain endomorphism and its induced action on homology by $\upsilon$.
In particular, $\upsilon \widecheck{\partial} + \widecheck{\partial} \upsilon = 0$.
See \cite[\S13.4]{Li22}.
Unlike the original $\hmr$, the operator $\upsilon$ is independent of the component on which the basepoint $p$ lies.

\begin{figure}[htbp]
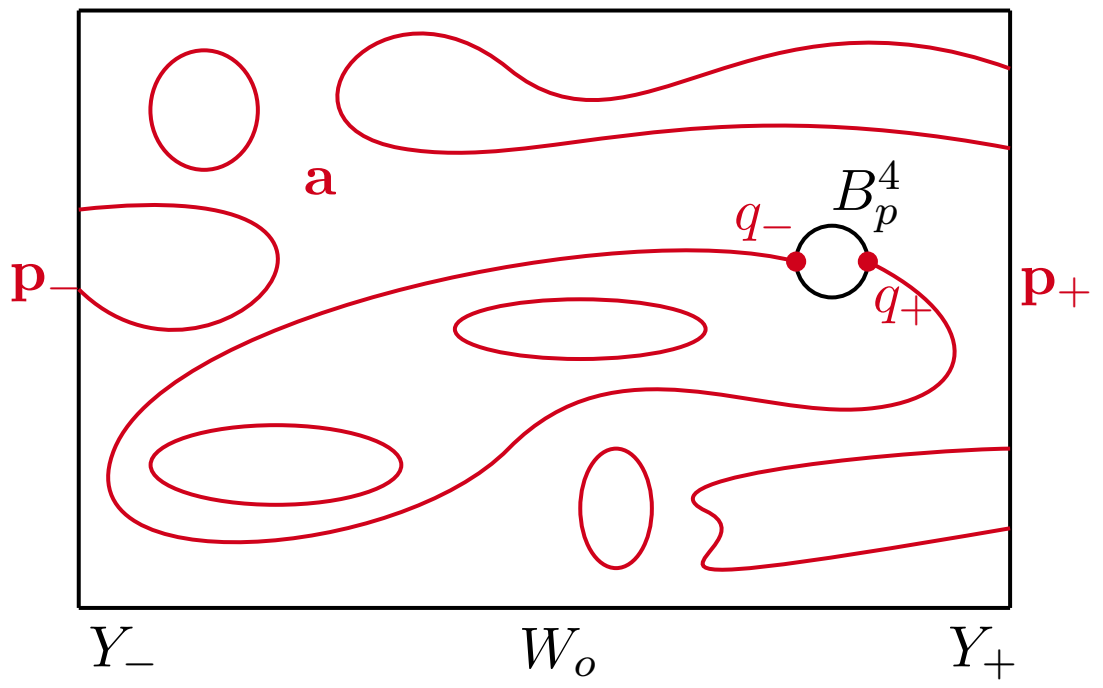

	\centering
	\figcobWo
	\caption{A real marked cobordism that is punctured at a point along the marking set.}
	\label{fig:cobWo}
\end{figure}

An alternative definition of the point class,
which goes back to Bloom~\cite{Bloom} in the ordinary monopole theory,
can be adapted to the real case~\cite{Li24} as follows.
One can turn the point evaluation map $\upsilon_p$ into a cobordism map with additional $(S^3, \tau_U)$ ends,
where $\tau_U$ is the deck transformation for the double branched cover of the unknot $U_1$.
In the current setup, given a marked cobordism $W$,
let $W_{o} = W \setminus B^4_p$,
where $B^4_{p}$ is a small $\tau$-invariant $4$-ball around $p \in \mathbf{a}$.
Furthermore, $B^4_p$ is viewed as a marked manifold $(B^4,\tau_B,I^1)$, where $\fix(\tau_B) = D^2$ containing the marking set $I^1 \subset D^2$.
The boundary of $B^4_p$ is a marked real 3-manifold $(S^3, \tau_U, \{q_-,q_+\})$.
There is a unique real \text{spin\textsuperscript{c}} structure on $S^3$,
but two \RRSC\ structures $\widetilde{\sr}_0,\widetilde{\sr}_1$.
See Figure~\ref{fig:cobWo},
and also \cite[Example~3.22--26]{Li26}.

Let $\widetilde{\sr}_0$ be the \RRSC\ structure that extends to $(B^4,\tau_B,I^1)$.
Now, view $W_o$ as a cobordism with three ends $-Y_-, -S^3, Y_+$.
Over $(S^3; \widetilde{\sr}_0)$, choose perturbations as in \cite{Li24} so that the set of critical points is a tower $\{\mathfrak{c}_i : i \in \mathbb{Z}\}$.
Define
	\[
	\widecheck m(\upsilon) =
\begin{bmatrix}
	m^{uo}_o(\mathfrak{c}_{-2} \otimes \ \cdot ) &
	m^{uu}_o(\mathfrak{c}_{-2} \otimes \overline{\partial}^s_u (\cdot) ) + \partial^u_o\bar m^{us}_u(\mathfrak{c}_{-2} \otimes \ \cdot )\\
	m^{uo}_s(\mathfrak{c}_{-2} \otimes \ \cdot ) &
	    \bar m_s^{us}(\mathfrak c_{-2}\otimes -)
	+ m_s^{uu}(\mathfrak c_{-2}\otimes\bar\partial_u^s(-))
	+ \partial_s^u\bar m_u^{us}(\mathfrak c_{-2}\otimes -).
\end{bmatrix}.
\]
Here, an entry of the form $m^{**}_*(\mathfrak{c}_{-2} \otimes \cdot)$ denotes the count of zero-dimensional moduli spaces on the triple-ended cobordism, with types specified by the superscripts and subscript,
and with the input $\mathfrak{c}_{-2}$ inserted at the $S^3$-end.
The neck-stretching argument in \cite{Li24} shows that this defines a chain map and is chain-homotopic to $\widecheck{m}(\upsilon | W)$.

Specialize to $Z = I \times (Y,\tau,\mathbf{p})$,
as in \cite[Example~3.24]{Li26}.
Let $\widetilde{\sr}_1$ be the \RRSC\ structure on $S^3$ that does not bound an \RRSC\ structure on $B^4$.
Inserting $\mathfrak{c}_{-2}$ at the $S^3$-end on the punctured cobordism $Z_o$,
but for $\widetilde{\sr}_1$ on $S^3$,
we obtain:
\begin{equation*}
	\upsilon_p^{\mathrm{ext}} \colon \widecheck{\HMR}_{\Delta,\bullet}(Y; \widetilde{\sr}) \to \widecheck{\HMR}_{\Delta,\bullet}(Y; \e_{p} \widetilde{\sr}),
\end{equation*}
Here, $\e_p \in \{0,1\}^{\mathbf{p}}$ denotes the vector with entry $1$ at $p$,
and $0$ everywhere else.
The \rrsc\ structure $\e_p\widetilde{\sr}$ is the result of a \emph{marking change}~\cite[Definition~3.2]{LM}.
Suppose $\widetilde{\sr}=(\mathfrak s, (T_p)_{p \in \mathbf{p}},\sr)$ is any \rrsc\ structure.
\begin{defn}
\label{defn:marking-change}
The \emph{marking change at $p$} acts on $\widetilde{\sr}$ by negating the Clifford trivialization $T_p$ at $p$. That is,
\begin{equation}
\label{eq:kronecker}
\e_p\widetilde{\sr}
:=\left(\mathfrak s,(\ldots,-T_p,\ldots),\sr\right).
\end{equation}
When the basepoints are ordered as $(p_1,\dots,p_n)$, we write $\e_j$ for the marking change at the $j$th point, $\e_{p_j}$.
\end{defn}

The maps $\upsilon_p^{\mathrm{ext}}$ have identical analytical definition for every $p \in \mathbf{p}$, as the internal V-map $\upsilon=\upsilon_p$.
Algebraically, one may view this set of external V-maps as the edge maps on the cube of chain complexes:
\begin{equation*}
	\bigoplus_{\mathsf{v} \in \{0,1\}^{\mathbf{p}}} \widecheck{\cmr}_{\Delta,\bullet}(Y; \mathsf{v}\widetilde{\sr}).
\end{equation*}

\begin{rem}
	These external V-maps seem to be related to the Kozul resolution in pointed Khovanov homology~\cite{BLS}.
	See Conjecture~\ref{conj:ptKh}.
\end{rem}

\subsection{Absolute gradings}
Unlike its ordinary counterpart, the real monopole Floer homology in~\cite{Li22} needs not have mod-two gradings.
This is because the Dirac operators in the ordinary case are complex,
and the Dirac spectral flows are always even.
The key idea in \cite{Li26} to guarantee relative mod-two gradings is to impose enough basepoints and take finite covers of the real configuration spaces so that no loop gives rise to odd Dirac spectral flows.
A sufficient condition (\cite[Proposition~9.3]{Li26}) is the following:

\begin{hypotheses}
	\label{hyp:abs-grad}
	Every component of the fixed point set $C$ contains precisely one basepoint in $\mathbf{p}$.
	Furthermore, the branched double cover $Y \to Y' = Y/\tau$ is induced by cutting along an oriented connected Seifert surface $F \subset Y'$ of $C$.
\end{hypotheses}

The existence of mod-two gradings~\cite{Li26} applies to the diagonal setting.
Indeed, while the discussion of \cite{Li26} is limited to \emph{framed} configuration spaces,
where one does not blow up along the reducible configurations,
the framed real Picard torus $H^1(Y;\mathbb{R})^{-\tau^*}/\Gamma(Y,\tau,\mathbf{p})$ in \cite{Li26} is identical to its diagonal counterpart.
Indeed, $\mathcal{G}_{\mathbf{p}}$ and $\mathcal{G}_{\Delta}$ differ only by the $\pm 1$ constant gauge transformations,
Moreover,
the order-two element of $\pi_1(\partial \mathcal{B}^{\sigma}_{\Delta})$ always has zero Dirac spectral flow.

\begin{figure}[htbp]
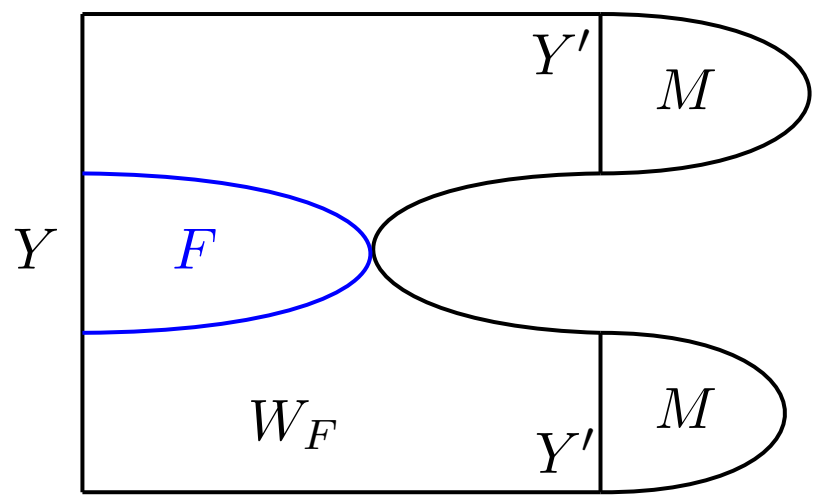

	\centering
	\figseifertcap
	\caption{The Seifert cobordism $W_F$ and a Seifert cap $W_{F,M}$. Here, the fixed point set $F$ (coloured blue) is viewed as the push-in copy of $F$ from the boundary $Y$.}
	\label{fig:seifertcap}
\end{figure}

We adapt \cite[Definition~4.1]{LM} from the framed setting as follows.
Fix a connected Seifert surface $F \subset Y'$ in Hypotheses~\ref{hyp:abs-grad}.
Let $W_F$ be the double cover of $[0,1]\times Y'$ branched
along a pushed-in copy of $F$ from $\{0\} \times Y'$ into the cylinder.
In particular, $W_F \colon Y \to Y' \sqcup Y'$ is a cobordism with involution $\tau$ that exchanges the two outgoing ends $Y'$.
Choose a 4-manifold $M$ with $\partial M = -Y'$ satisfying $H_1(M) = 0$.
Cap $W_F$ off at $Y' \sqcup Y'$ using two copies of $M$,
and let $W_{F,M} = W_F \cup_{Y' \sqcup Y'}(M \sqcup M)$ be the resulting 4-manifold that bounds $Y$.
Then $W_F$ is the \emph{Seifert cobordism} and $W_{F,M}$ is a \emph{Seifert cap}.
See Figure~\ref{fig:seifertcap}.

The 4-manifold $M$ is an \emph{outgoing cap}.
Over the two copies of $M$'s that are glued to $W_F$, take conjugate \text{spin\textsuperscript{c}} data so that
$\tau$ pulls back every \text{spin\textsuperscript{c}} structure to its conjugate upon identifying the two $M$'s.
Fix an \RRSC\ structure $\widetilde{\sr} = (\mathfrak{s},T,\sr)$ on $(Y,\tau,\mathbf{p})$,
and choose a real spin\textsuperscript{c} structure $(\widehat{\mathfrak{t}}, \sr_W)$ on $W_{F,M}$ rel $Y$ extending $(\mathfrak{s},\sr)$.

Given a configuration $\mathfrak{a}$ and perturbation $\mathfrak{q}$,
let $\gr^{\sigma}_{SW}(\mathfrak{a};W_{F,M},\widehat{\mathfrak{t}})$ be the formal dimension of the blown-up Seiberg--Witten moduli space on $W_{F,M}$ with APS boundary conditions specified by $(\mathfrak{a}, \mathfrak{q})$,
using the convention of the spectral projection in \cite[\S20.3]{KM}.
Let $-\iota_R(W_{F,M})$ be the APS index of $(d^* \oplus d^+)_R$ with the same boundary condition:
\[
(d^*\oplus d^+)_R:
\Omega^1(W_{F,M};\mathbb{R})^{-\tau^*_W}\longrightarrow
\Omega^0(W_{F,M};\mathbb{R})^{-\tau^*_W}\oplus\Omega^{+}(W_{F,M};\mathbb{R})^{-\tau_{W}^*}.
\]
If $Y$ is connected, then
\[\iota_R(W_{F,M}) = b^+_R(W_{F,M}) - b^1_R(W_{F,M}),\]
where the subscript $R$ denotes taking $\tau^*_W$-anti-invariant part of the Betti numbers.

By \cite[\S3]{LM}, an extending real \text{spin\textsuperscript{c}} structure on $W_{F,M}$ rel $Y$ and the Seifert surface $F$ determines a spin structure $\mathfrak{s}_{\mathfrak{t}}$ on $F$.
Fix an ordering of the components of $C$.
Then the data of $(\widehat{\mathfrak{t}}$, $F$, $\widetilde{\sr})$ together determine
\cite[\S3]{LM} a spin structure $\widehat{\mathfrak{s}}_{\mathfrak{t}}$ on $\hat{F} = F \cup G$, where $G$ is a planar surface of $n = |C|$ boundary components.
Following \cite[Definition~4.1]{LM},
the \emph{absolute mod-two grading} of the upstairs configuration $\mathfrak{a}$ is
\begin{equation}\label{eq:absolute-grading}
		\gr^{(2)}_{\Delta}(\mathfrak{a})
		= \gr^{\sigma}_{SW}(\mathfrak{a};W_{F,M},\widehat{\mathfrak{t}})
		+\iota_R(W_{F,M})
		+\Arf(\widehat{\mathfrak s}_{\mathfrak t})
		\pmod2.
\end{equation}
The original \cite[Definition~4.1]{LM}, using the downstairs version $\gr_{SW}$ of the Seiberg--Witten operator, applies to $\alpha \in \mathcal{B}_{\Delta}$:
\begin{equation}\label{eq:absolute-grading-down}
		\gr^{(2)}(\alpha)
		= \gr_{SW}(\alpha;W_{F,M},\widehat{\mathfrak{t}})
		+\iota_R(W_{F,M})
		+\Arf(\widehat{\mathfrak s}_{\mathfrak t})
		\pmod2.
\end{equation}
Let $\pi \colon \mathcal{B}^{\sigma}_{\Delta} \to \mathcal{B}_{\Delta}$ be the blow-down map.
The two mod-two gradings can be related, in the case of critical points as follows.
\begin{lemma}
	\label{lem:crit-same-grW}
If $\mathfrak{a}$ is an irreducible \emph{critical point}, then
\begin{equation}\label{eq:crit-same-grW-irr}
	\gr^{\sigma}_{\mathrm{SW}}(\mathfrak{a}; W, \widehat{\mathfrak{t}})= \gr_{\mathrm{SW}}(\pi(\mathfrak{a}); W, \widehat{\mathfrak{t}}).\end{equation}
Moreover, for a reducible \emph{critical point} downstairs, let $\vartheta_0$ be the blownup critical point corresponding to the smallest positive eigenvalue $\lambda_0$ of the perturbed Dirac operator $D_{\mathfrak{q}, \vartheta}$.
Then
\begin{equation}\label{eq:crit-same-grW-red}
	\gr^\sigma_{\mathrm{SW}}(\vartheta_0;W,\widehat{\mathfrak t})
	=
	\gr_{\mathrm{SW}}(\vartheta;W,\widehat{\mathfrak t}).
\end{equation}
\end{lemma}

\begin{proof}
The blow-down map $\pi$ is a local diffeomorphism at every irreducible configuration,
which on the 3-dimensional configuration spaces of $Y$, identifies $\grad \pertL$ with $\grad \pertL^{\sigma}$.
At an irreducible critical point $\mathfrak{a}$, the differential $d\pi$ intertwines the Hessians $\Hess_{\pi(\mathfrak{a})}$ and  $\Hess_{\mathfrak{a}}^{\sigma}$ as in~\cite[Lemma~12.4.3]{KM}.
This implies the APS boundary conditions for the upstairs and downstairs problems are conjugate.
Since $d\pi$ also identifies Seiberg--Witten operator up to compact terms over the interior, Equation~\eqref{eq:crit-same-grW-irr} follows.

The following argument in the reducible case is due to Lidman--Manolescu~\cite[Proof of Proposition~9.2.1]{LidmanManolescu} (see also \cite[Proposition~4.22]{Xiao}).
Suppose $\vartheta = (B,0)$ is a downstairs reducible critical point.
Let $\vartheta_0 = (B, 0, \psi_0)$ be the corresponding upstairs critical point,
where $\psi_0$ is a unit $\lambda_0$-eigenvector.
At a reducible configuration, the upstairs and downstairs gauge-fixing terms that both act on $T\mathcal{A}_k^{\sr}$ coincide:
$\mathbf d^{\sigma,\dagger}_{\vartheta_0}=\mathbf d^{*}_{\vartheta}=-d^*$.
Indeed, the averaging term $\mu_Y$ in the formula $\mathbf d^{\sigma,\dagger}$ vanishes for real spinors.
By the real counterparts of \cite[Proposition~12.2.5 and 12.3.1]{KM}, $\mathcal{D}_{\vartheta}(\grad \pertL)$ and $\mathcal{D}_{\vartheta_0}(\grad \pertL)^{\sigma}$ have the following matrix expressions:
\[\begin{bmatrix}
	\mathcal{D}_B(\grad \pertL)^{\mathrm{red}} & 0\\
	0 & D_{\mathfrak{q},B}
\end{bmatrix}
,\qquad
\begin{bmatrix}
	\mathcal{D}_B(\grad \pertL)^{\mathrm{red}} & 0 & 0\\
	0 & \lambda_0 & 0\\
	x & 0 & D_{\mathfrak{q},B} - \lambda_0
\end{bmatrix}.\]
The first operator acts on $T\mathcal{A}_k^{\sr} \oplus L^2_k(Y;S)^{\sr}$ and the second operator acts on $T\mathcal{A}_k^{\sr} \oplus \mathbb{R} \oplus \psi_0^{\perp}$, where $\psi_0^{\perp}$ denotes the $L^2$-orthogonal complement of $\psi_0$ in the space $L^2_k(Y;S)^{\sr}$ of real spinors.
Moreover, the following operator is invertible:
\[
L=\begin{bmatrix}
0 & -d^*\\
-d & \mathcal{D}_B(\grad \pertL)^{\mathrm{red}}
\end{bmatrix}.
\]
This operator acts on $L^2_k(Y,\mathbf{p};i\mathbb{R})^{-\tau^*} \oplus T\mathcal{A}_k^{\sr}$.
The invertibility of $L$ follows from non-degeneracy and the invertibility of the Laplacian on $\tau$-anti-invariant functions.

Consider the following family of operators for $0 \le t \le 1$ acting on $L^2_k(Y,\mathbf{p};i\mathbb{R})^{-\tau^*} \oplus T\mathcal{C}^{\sigma}_k$,
for which $K_0$ is the extended Hessian:
\[
K_t=\begin{bmatrix}
	L & 0 & 0\\
	0 & \lambda_0 & 0\\
	(1-t)(x+\psi_0 \cdot \,) & 0 & D_{\mathfrak{q},B} - (1-t)\lambda_0
\end{bmatrix}.\]
Since $\lambda_0$ is simple and $D_{\mathfrak{q},B}$ has no eigenvalues in $[0,\lambda_0)$,
every $K_t$ is hyperbolic.
Moreover, we have $T_Y K_1T_Y^{-1} = L \oplus D_{\mathfrak{q},B}$,
where
\[T_Y(c,b,r,\psi) := (c,b,r\psi_0+\psi).\]
This shows that the upstairs extended Hessian can be identified, via a homotopy through hyperbolic operators, with the downstairs Hessian.

Now, extend arbitrarily $\psi_0$ from $Y$ to a smooth spinor $\Phi$ on $W$,
viewed as a class of $\Gamma(S^+)^{\sr}/\mathbb{R}_{>0}$.
Denote, by $Q^{\sigma}_W$ and $Q_{W}$, the linearized Seiberg--Witten operators with gauge fixing in the definition of $\gr^{\sigma}_W$ and $\gr_{W}$, respectively.
The map $T_W(r,\phi)=r\Phi+\phi$, for $\phi\in\Phi^\perp$,
identifies the domains of $Q^{\sigma}_W$ and $Q_{W}$.
Furthermore, $T_W$ commutes with the upstairs/downstairs boundary restriction maps up to a rank-one term due to the projection onto $\psi_0^\perp$.
After identifying the domains and trivializing the codomains using $\Phi$, the operators $Q^{\sigma}_W$ and $Q_{W}$ differ by compact
lower-order and finite-rank terms,
as in \cite[proof of Proposition~24.3.1]{KM}.
Let $\Pi_t$ be the spectral projections of $T_Y K_tT_Y^{-1}$.
Then by \cite[Lemma~17.2.2, Propositions~17.2.4--5]{KM}, the $\Pi_t$'s are all commensurate.
After identifying the ranges of $\Pi_t$, the family $Q^{\sigma}_W\oplus\Pi_t$ is continuous and hence has constant index.
At $t=0$ it computes $\gr^\sigma_{\mathrm{SW}}
 (\vartheta_0;W,\widehat{\mathfrak t})$,
 and at $t=1$ it computes the downstairs $\gr_{\mathrm{SW}} (\vartheta;W,\widehat{\mathfrak t})$.
 Thus \eqref{eq:crit-same-grW-red} follows.
\end{proof}

\begin{prop}\label{prop:absolute-grading-property}
The grading \eqref{eq:absolute-grading} is independent of the outgoing cap $M$ and \text{spin\textsuperscript{c}} structure $\widehat{\mathfrak{t}}$.
Furthermore, for any \emph{configurations} $\mathfrak a,\mathfrak b$ and any homotopy class of path $z$ from $\mathfrak a$ to $\mathfrak b$, the absolute grading satisfies
\begin{equation}
	\label{eq:abs-to-rel-gr}
 \gr_\Delta^{(2)}(\mathfrak a)
 -\gr_\Delta^{(2)}(\mathfrak b)
 \equiv \gr_z(\mathfrak a,\mathfrak b)\bmod2.
\end{equation}
In addition, for any irreducible critical point $\mathfrak{a}$,
\[
 \gr^{(2)}_\Delta(\mathfrak a)
 =\gr^{(2)}(\pi\mathfrak a).
\]
Finally, if $\vartheta_0$ is a reducible critical point corresponding to the first positive eigenvalue of a downstairs reducible, then
\[
\gr^{(2)}_\Delta(\vartheta_0)
=\gr^{(2)}(\vartheta).
\]
\end{prop}

\begin{proof}
The claim about independence of $(M,\widehat{\mathfrak{t}})$ follows from the proof of \cite[Theorem~4.3(i)]{LM} verbatim.
Similarly, the proof of \cite[Proposition~4.2(i)]{LM} implies \eqref{eq:abs-to-rel-gr}.
Suppose now that $\mathfrak a$ is an irreducible critical
point.
The $\iota_R$ and Arf terms in \eqref{eq:absolute-grading} and \eqref{eq:absolute-grading-down} are identical.
By Lemma~\ref{lem:crit-same-grW}, the terms $\gr_{SW}^{\sigma}$ and $\gr_{SW}$ also agree.
\end{proof}

\begin{rem}\label{rem:inadmissible}
	In the definition of $\grtwo$, the outgoing cap is assumed to be admissible, i.e. $H_1(M;\mathbb{Z})$ vanishes.
	This guarantees every spin\textsuperscript{c} structure can be extended from $Y'$ to $M$, without changing the manifold $M$.
	The proof of independence of $(M,\widehat{\mathfrak{t}})$ does not require $M$ to be admissible.
	As far as grading computation is concerned, one may work with an inadmissible outgoing cap, such as in Lemma~\ref{lem:splitting-degree}.
\end{rem}
\begin{rem}
Using the mod-two grading above, the external V-maps $\upsilon^{\mathrm{ex}}_p$ may not have absolute degree $-1$.
Indeed, the external V-maps can shift the \rrsc\ structures,
which can result in additional grading changes, according to Corollary~\ref{cor:marking-change-euler}.
\end{rem}

\section{Cone real monopole Floer homology}
\label{sec:cone-hmr}
The construction of the cone real monopole Floer homology follows that of the vanilla single-pointed tilde-version~\cite{Li24}, except with diagonal gauge groups and V-maps.

\subsection{The cone complex}
Fix an admissible perturbation $\mathfrak{q}$ throughout this section.
Let $\widetilde{\sr}$ be an \RRSC\ structure on $Y$, and let $\check{C}$ denote the diagonal to-chain complex of $(Y, \tau, \mathbf{p}; \widetilde{\sr})$.
We conflate the notation $\upsilon$ on homology and $\widecheck m(\upsilon)$ of the chain representative of $\upsilon$.
\begin{defn}
	The \emph{cone real monopole Floer} chain complex and homology are, respectively,
	the mapping cone chain complex of $\widecheck m(\upsilon)$ and its homology:
	\begin{equation*}
		\widetilde C_{\Delta} = \widecheck{C}\oplus \widecheck{C},\quad
		\widetilde{\partial} =\begin{bmatrix}\widecheck{\partial}&0\\
		\widecheck m(\upsilon) & \widecheck{\partial} \end{bmatrix};\qquad
		\tH_\Delta(Y,\tau,\mathbf p;\widetilde{\sr})
		=H(\widetilde C_{\Delta},\widetilde{\partial}).
	\end{equation*}
Here, the two copies of $\widecheck{C}$ are graded the same way, without shift.
\end{defn}
Since $\upsilon\check{\partial}+\check{\partial}\upsilon = 0$, we have $\tilde{\partial}^2 = 0$ and hence $\tilde{\partial}$ is a well-defined differential.
Moreover, the absolute grading $\grtwo{}_{\Delta}$ induces an absolute mod-two grading on $\tilde{C}_{\Delta}$ by grading each copy of $\check{C}$ the same way.
When $\mathbf p=\{p\}$ consists of a single point, $\mathcal{G}_{\Delta}=\mathcal{G}$, so this construction recovers the single-pointed version \cite[Definition~2.1]{Li24}.

\begin{lem}
There is a natural short exact sequence
\begin{equation}
\label{eq:v-SES}
0\longrightarrow\coker\bigl(\upsilon_*\colon H_{j+1}(\widecheck C)\to H_j(\widecheck C)\bigr)
\longrightarrow\tH_{\Delta,j}
\longrightarrow\ker\bigl(\upsilon_*\colon H_j(\widecheck C)\to H_{j-1}(\widecheck C)\bigr)
\longrightarrow0.
\end{equation}
\end{lem}

\begin{proof}
Inclusion of the second copy of $\widecheck C$ into $\widetilde C_\Delta$ and projection onto the first define a short exact sequence:
$0\to\widecheck C\to\widetilde C_\Delta\to\widecheck C\to 0$.
The connecting homomorphism is $\upsilon$.
Indeed, the differential of the lift $(x,0)$ of a cycle $x$ is $(0,\upsilon x)$.
The above short exact sequence follows from the long exact sequence as in \cite[\S2.2]{Li24}.
\end{proof}

\subsection{Finiteness and Euler characteristic}
Denote by $\mathfrak{O}_\Delta$ and $\mathfrak{R}_\Delta$ the sets of downstairs irreducible and the reducible critical points, respectively, inside $\mathcal{B}_{\Delta}$.
By compactness downstairs, both $\mathfrak{O}_\Delta$ and $\mathfrak{R}_\Delta$ are finite sets.
Under the projection $\widetilde{B}_{\mathbf{p}} \to \mathcal{B}_{\Delta}$,
an irreducible configuration has two lifts of the same mod-two grading;
a reducible has one lift, because the constant gauge transformation $-1$ acts trivially on the reducibles.
The following generalizes the argument in \cite[\S2.3]{Li24} to the case when $b^1_R(Y) > 0$.
\begin{prop}\label{prop:finite-sector}
	The cone real monopole Floer homology $\tH_\Delta(Y,\tau,\mathbf{p}; \widetilde{\sr})$ is a finite-dimensional $\mathbb{F}_2$-vector space for each  \RRSC\ structure $\widetilde{\sr}$.
	Moreover, if the underlying \text{spin\textsuperscript{c}} structure of $\widetilde{\sr}$ is torsion, then
	\[
	\dim_{\F}\tH_\Delta(Y,\tau,\mathbf p;\widetilde{\sr})
	\leq 2|\mathfrak O_\Delta|+|\mathfrak R_\Delta|.
	\]
	Finally,
	\begin{equation}\label{eq:euler-count}
		\chi_{\Delta}\left(\tH_\Delta(Y,\tau,\mathbf p;\widetilde{\sr})\right)=
		2\sum_{\alpha\in\mathfrak O_\Delta}(-1)^{\gr^{(2)}(\alpha)}
		+\sum_{\vartheta \in\mathfrak R_\Delta}(-1)^{\gr^{(2)}(\vartheta)},
	\end{equation}
	whenever $\grtwo$ is well-defined.
\end{prop}

\begin{proof}
	If $c_1(\mathfrak s)$ is non-torsion, admissibility excludes reducible critical points, and hence there are only finitely many upstairs irreducible generators.
	Then~\eqref{eq:euler-count} follows from Proposition~\ref{prop:absolute-grading-property} and the fact that two $\check{C}$'s have the same grading.
	Suppose $\mathfrak{s}$ is torsion.
	Then $\pertL$ is single-valued, and there is a relative $\mathbb{Z}$-grading.
	We have
	\[
	\widecheck C=\F\mathfrak O_\Delta\oplus
	\bigoplus_{\vartheta\in\mathfrak R_\Delta,\ j\geq0}\F \vartheta_{j}.
	\]
	Here, $\vartheta_j = (\vartheta, \Psi_{\vartheta,j})$  for which $\Psi_{\vartheta,j}$ is the eigenvector with the $(j+1)$-th smallest positive eigenvalue of the perturbed Dirac operator $D_{\vartheta,\mathfrak{q}}$.
	We filter the \emph{to} chain complex by the values of the Chern--Simons--Dirac functional as follows.
	Let $x_1 < \dots < x_m$ be the values of the CSD functional, and let $\mathcal{F}_k \check C$ be the subgroup generated by critical points of value at most $x_k$.
	Let $\Gr_{\mathcal{F}} = \bigoplus_{k} \mathcal{F}_k \check C/\mathcal{F}_{k-1} \check C$.
Choose the representative $\widecheck m(\upsilon)$ by evaluating $\upsilon$ on moduli spaces of trajectories of $\grad\pertL$, so that $\widecheck m(\upsilon)$ preserves $\mathcal{F}_k$.
	Furthermore, $\check{\partial}$ also preserves this filtration.
	Indeed, the differential $\check{\partial}$ strictly decreases $\pertL$ along trajectories with nontrivial downstairs image.
	For reducible trajectories lying above the same $\vartheta$ downstairs, there are precisely two index-1 trajectories from $\vartheta_j$ to $\vartheta_{j-1}$, contributing zero to the differential $\check{\partial}$.
	Thus $\check{\partial} \mathcal{F}_k \subset \mathcal{F}_{k-1}$.
	Moreover,
	$\Gr_{\mathcal{F}}(\upsilon)(\vartheta_j)=\vartheta_{j-1}$ for $j \ge 1$, and
    $\Gr_{\mathcal{F}}(\upsilon)(\vartheta_0)=0$;
    for any irreducible $\mathfrak{a}$, we have $\Gr_{\mathcal{F}}(\upsilon)(\mathfrak{a})=0$ .

	Now, filter the cone complex $\tilde{C}$ via $\mathcal{F}_k\check{C} \oplus \mathcal{F}_k\check{C}$.
	The first differential in the filtration spectral sequence is given by
	\[
	\Gr_{\mathcal{F}}(\widetilde{\partial}) = \begin{bmatrix}
		0 & 0\\
		\Gr_{\mathcal{F}}(\upsilon) & 0
	\end{bmatrix}.\]
Every $\alpha \in \mathfrak{O}_{\Delta}$ contributes to the associated homology two generators, both of grading $\grtwo(\alpha)$.
For a reducible $\vartheta$, distinguish the two copies of $\check{C}$ by writing $\vartheta_j, \vartheta_j'$.
The graded differential sends $\vartheta_j$ to $\vartheta_{j-1}'$.
The two reducible towers $(\vartheta_j)_{j\ge 0}$ and $\{\vartheta_j'\}_{\ge 0}$ contribute precisely one element of homology, represented by $\vartheta_0$, with grading $\grtwo(\vartheta)$.
Thus the $E_1$ page of the spectral sequence is finite,
and the rank of the $E_{\infty}$ page is bounded by $2|\mathfrak{O}_{\Delta}| + |\mathfrak{R}_{\Delta}|$.
The Euler characteristic equality~\eqref{eq:euler-count} then follows from Lemma~\ref{lem:crit-same-grW}.
\end{proof}

The proof of \cite[Lemma~6.5]{Li26}, combined with Proposition~\ref{prop:finite-sector}, implies the following:

\begin{cor}
Suppose the tame perturbations over all \RRSC\ structures satisfy a uniform
bound as in \cite[Definition~4.11(iv)]{Li26}.
Then only finitely many
\RRSC\ structures admit nontrivial $\tH_\Delta(Y,\tau,\mathbf p)$.
In particular, the total rank is finite.\qed
\end{cor}

From now on, the \emph{total Euler characteristic} refers to the sum of Euler characteristics over all \RRSC\ structures:
\[\chi\left(\tH_\Delta(Y,\tau,\mathbf p)\right) :=
\sum_{\widetilde{\sr} \in \RRSC(Y,\tau,\mathbf{p})}
\chi\left(\tH(Y,\tau,\mathbf p;\widetilde{\sr})\right).\]

\subsection{Relation with the framed Euler characteristic}
The following proposition equates the Euler characteristics of the cone and framed real monopole Floer homologies.
Suppose Hypotheses~\ref{hyp:abs-grad} hold.
Fix a connected Seifert surface $F$ and an ordering of $\pi_0(C)$ throughout this subsection.
The following is the restatement of Theorem~\ref{thm:cone-framed-euler}.
\begin{prop}
For every \RRSC\ structure $\widetilde{\sr}$, the Euler characteristic of the cone real monopole Floer homology coincides with that of the framed real monopole Floer homology:
\[\chi\left(\tH_\Delta(Y,\tau,\mathbf p;\widetilde{\sr})\right)
=\chi\left(\tH(Y,\tau,\mathbf p;\widetilde{\sr})\right).\]
\end{prop}
\begin{proof}
	Choose a $\mathcal{G}$-invariant admissible perturbation
	$\mathfrak q$ as above.
	By \cite[Theorem~7.2]{Li22}, all the critical points in $\mathcal{C}/\mathcal{G}_{\Delta}$ and hence in $\mathcal{C}/\mathcal{G}_{\mathbf p}$ are nondegenerate.
	By the proof of \cite[Theorem~5.11]{Li26},
	choose a small framed perturbation $\tilde{\mathfrak q}$ equal to
	$\mathfrak q$ on neighborhoods of the critical points in $\tilde{\mathcal{B}}_{\mathbf{p}}$,
	for which all moduli spaces of trajectories are regular.
Moreover, for this $\tilde{\mathfrak{q}}$, no new critical points appear.
Hence the perturbations ${\mathfrak{q}}$ and $\tilde{\mathfrak{q}}$ have the same critical points and Hessians on the framed blown-down configuration space.
	Each diagonal irreducible critical point downstairs has two lifts, of the same mod-two grading $\grtwo$.
	Each diagonal reducible critical point in $\mathcal{B}_{\Delta}$ has one lift.
	By Proposition~\ref{prop:absolute-grading-property} and Proposition~\ref{prop:finite-sector},
	\[
	\chi\left(\tH(Y,\tau,\mathbf p;\widetilde{\sr})\right)
	=2\sum_{\alpha\in\mathfrak O_\Delta}(-1)^{\gr^{(2)}(\alpha)}
	+\sum_{\vartheta \in\mathfrak R_\Delta}(-1)^{\gr^{(2)}(\vartheta)}= \chi \left(\tH_{\Delta}(Y,\tau,\mathbf p;\widetilde{\sr})\right).\qedhere
	\]
\end{proof}

\begin{cor}
	\label{cor:marking-change-euler}
	Given a \RRSC\ structure $\widetilde{\sr}$, let $w_{1,j}=w_{1,j}(\widetilde{\sr})$ denote the Stiefel--Whitney number of the eigenline bundle of $S^{\sr}$ along $C_j$, determined by the Seifert framing as in \cite[\S3.3]{Li26}.
 	Every marking change $\e_j$ as in~\eqref{eq:kronecker} induces a sign change of $\chi$:
	\begin{equation}
	\label{eq:ej-on-rrsc-sign}
		\chi \left(\tH_\Delta(Y,\tau,\mathbf p;\e_j\widetilde{\sr})\right)
	=(-1)^{w_{1,j}+1}
	\chi\left(\tH_\Delta(Y,\tau,\mathbf p;\widetilde{\sr})\right).
	\end{equation}
	Hence the total Euler characteristic has only contributions from \RRSC\ structures satisfying $w_{1,j}=1$, \emph{for every} $j$.
	That is,
	\begin{equation}\label{eq:lambda-contribution-total}
		\chi\left(\tH_\Delta(Y,\tau,\mathbf p)\right)
		=\sum_{\substack{[\widetilde{\sr}]\in\RRSC(Y,\tau,\mathbf p)\\
				w_{1,j}(\widetilde{\sr})=1, \, \forall i}}
		\chi\left(\tH_\Delta(Y,\tau,\mathbf p;\widetilde{\sr})\right).
	\end{equation}
\end{cor}
\begin{proof}
The marking change $\e_j$ from $T_{p_j}$ to $-T_{p_j}$ preserves the diagonal gauge group, the configuration space, and the perturbed Seiberg--Witten equations.
The identity map on the configuration space hence induces an ungraded equivalence between the chain complexes associated to $\widetilde{\sr}$ and $\e_j\widetilde{\sr}$.
Thus the $\gr(\mathfrak{a},W_{F,M};\hat{\mathfrak{t}})$-term and $\iota_R$-term in \eqref{eq:absolute-grading} remain unchanged.
By \cite[Lemma~3.12]{LM}, on the planar closure of $F$, 
the spin structure changes by the Poincar\'{e} dual of $[C_j]$ under $\e_j$.
Then \cite[Lemma~3.14]{LM} says
\[\Arf(\widehat{\mathfrak{s}}_{\mathfrak{t}}, \varepsilon(\e_j\widetilde{\sr})) - \Arf(\widehat{\mathfrak{s}}_{\mathfrak{t}}, \varepsilon(\widetilde{\sr}))=\lambda([C_j]) =1+ w_{1,j}.\]
Here, we have suppressed the notation $F$ and $\mathsf{m}$ from \cite[Lemma~3.14]{LM} as they remain the same under marking changes.
Since only the Arf term contributes to change of the mod-two grading, this proves the first equation.

For \eqref{eq:lambda-contribution-total}, suppose $\vec{w} \in \{0,1\}^ n$ is a vector whose $j$-th component is not equal to $1$.
The marking change $\widetilde{\sr} \mapsto \e_j\widetilde{\sr}$ acts on the set of \RRSC\ structures with $(w_{1,j})_{1 \le j \le n}=\vec{w}$.
Equation~\eqref{eq:ej-on-rrsc-sign} implies that contribution of this set of \RRSC\ structures is zero to the total Euler characteristic.
Thus the only contributions arise from those with $w_{1,j}=1$ for all $1 \le j \le n$.
\end{proof}

\begin{rem}
\label{rem:w1-lambda}
The quantity $1+w_{1,j}$ is equal to $\lambda_j$ as in \cite[\S3.5]{LM},
referring to the boundary spin structure type on the circle $C_j$ in $\partial F$.
Specifically, $\lambda_j = 0$ indicates the bounding spin structure and $\lambda_j = 1$ indicates the Lie group spin structure.
For $\lambda=0$, the Arf invariant on the planar closure of $F$ is equal to the Arf invariant of the quadratic refinement induced on the quotient of $H_1(F;\F)$ by the radical of its intersection form.
\end{rem}

\subsection{Cobordism maps}
\label{ssec:cob-maps}
In this subsection, we define a restricted class of cobordism maps using punctured marked cobordisms.
Fix a connected marked real cobordism between two connected marked real 3-manifolds:
\[(W,\tau_W,\mathbf a)\colon(Y_-,\tau_-,\mathbf p_-)\longrightarrow(Y_+,\tau_+,\mathbf p_+).\]
The homomorphism $\epsilon_{\mathbf a} \colon \mathcal{G}_{\Delta,\mathbf{a}} \to \{\pm1\}$ restricts to
$\epsilon_{\mathbf p_-}$ and $\epsilon_{\mathbf p_+}$ on the $Y_{\pm}$ ends.
This defines a real line bundle $\mathbb L_{\mathbf a}=\mathcal{C}^\sigma(W)\times_{\epsilon_{\mathbf a}}\mathbb{R}/\mathcal{G}_{\Delta, \mathbf{a}}$.
Let $\upsilon_{\pm} = \upsilon_{\mathbf{p}_{\pm}}$,
and $\upsilon_{\mathbf a}=w_1(\mathbb L_{\mathbf a})$, which restricts to $\upsilon_{\pm}$.
\emph{Assume} there exists at least one arc $a_*$ with two endpoints on $Y_-$ and $Y_+$, respectively.
We define cobordism map
\[\thmr_{\Delta, a_*}(W; \widetilde{\sr}_W) \colon \tHMR_{\Delta,\bullet}(Y_-; \widetilde{\sr}_-) \to \tHMR_{\Delta,\bullet}(Y_+; \widetilde{\sr}_+),\]
as follows.

\begin{rem}
Such an $a_*$ may not exist for general marked real cobordisms.
However, this class of maps suffices for the proofs of the exact triangles, in which all basepoints stay constant during the skein/surgery cobordisms.
(For example, the cobordism in Figure~\ref{fig:cobW}.)
Moreover, it allows us to appeal to the identities of the boundary stratifications from \cite{Bloom} for parametrized moduli spaces.
This simplifies the proofs in \S\ref{sec:exact-triangle}.
In principle, one may define the cobordism map in the cone real monopole Floer homology using the classes $\upsilon_{\mathbf{a}}, \upsilon_{\mathbf{p}_{\pm}}$ via the \v{C}ech approach~\cite{KM}.
\end{rem}
\begin{figure}[htbp]
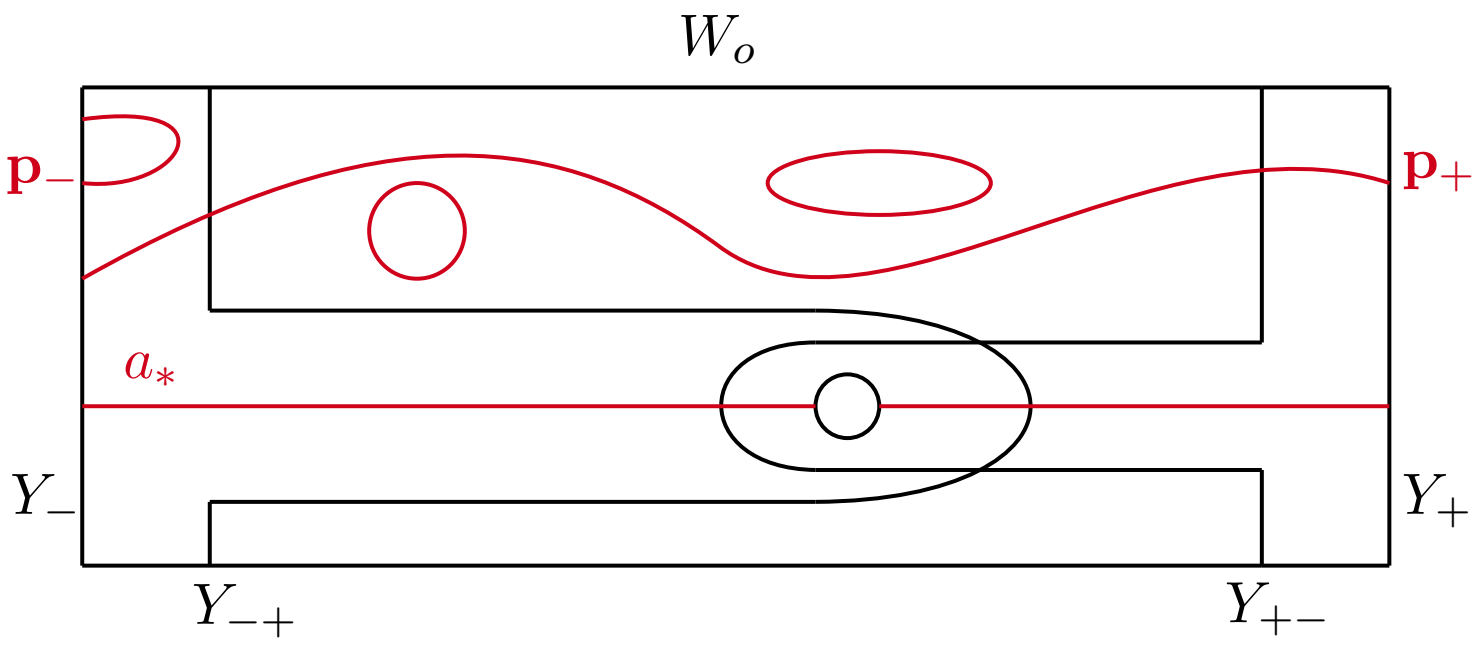

	\centering
	\figWplusminus
	\caption{Two hypersurfaces $\{Y_{-+},Y_{+-}\}$ in $W_o$, with a puncture at $p \in a_{*}$. The union of all red arcs is the marking set $\mathbf{a} \setminus \mathrm{nbhd}(p)$.}
	\label{fig:Wplusminus}
\end{figure}

Fix an arc $a_*$.
For each interior point $p \in a_*$, let $W_o = W \setminus B^4_p$ as in \S\ref{ssec:V-map}.
Following \cite[\S2.2 and Figure~2]{Li24}, introduce two hypersurfaces $Y_-+$ and $Y_{+-}$ diffeomorphic to $Y_{\pm}$
as in Figure~\ref{fig:Wplusminus}.
Consider a $[-\infty, +\infty]$-family $g_T$ of $\tau_W$-invariant Riemannian metrics on $W_o$ that neck-stretch along $Y_{-+}$ for $T \in [-\infty,0]$ and then along $Y_{+-}$ for $T \in [0,+\infty]$.
We refer to this type of metrics as \emph{puncture-moving}.

Fix a \RRSC\ structure $\widetilde\sr_W=(\mathfrak t,T_{\mathbf a},\sr_W)$ on $W$ restricting to $\widetilde\sr_\pm$ on $Y_\pm$.
Let $(\check C_\pm,\check\partial_\pm)$ be the \HMR-\emph{to} complexes of $(Y_\pm,\tau_\pm,\mathbf p_\pm;\widetilde\sr_\pm)$,
respectively.
Counting zero-dimensional family moduli spaces over $g_T$ defines an operator
\begin{equation}
\label{eq:checkKdefn}
	\widecheck K = \widecheck{K}_W \colon \widecheck{C}_{-} \to \widecheck{C}_{+}.
\end{equation}
Here, we are using the punctured cylinder model for $\widecheck m(\upsilon_{\pm})$.
By the same boundary stratification analysis
as in \cite[Remark~2.3 \& Lemma~4.8]{Li24},
\begin{equation*}
	\widecheck\partial_+ \widecheck K+ \widecheck K\widecheck\partial_-=\upsilon_+\widecheck m(W)+\widecheck m(W)\,\upsilon_-.
\end{equation*}
Here, we use
\cite[Proposition~25.3.8, Lemma~26.2.2]{KM}, and \cite[Propositions~12.10, 13.11--13.13]{Li22}.
These results apply to the diagonal setting, as $\mathcal{B}^\sigma_{\Delta}(W)$ is a finite cover
over $\mathcal{B}^\sigma(W)$.
Since $\check\partial_+\widecheck m(W)+\widecheck m(W)\check\partial_-=0$, the following is a chain map (\cite[Definition~2.2]{Li24}) that defines the cobordism map on $\thmr_{\Delta}$:
\begin{equation}\label{eq:cone-map-summary}
	\widetilde m_{\Delta,a_*}(W;\widetilde\sr_W)=\begin{bmatrix}\widecheck m(W; \widetilde{\sr}_W)&0\\ \widecheck K&\widecheck m(W; \widetilde{\sr}_W)\end{bmatrix}
	\colon\widetilde C_-\longrightarrow\widetilde C_+.
\end{equation}
Denote the resulting cobordism map as
\begin{equation*}
	\thmr_{\Delta,a_*}(W) \colon \thmr_{\Delta,\bullet}(Y_-) \to \thmr_{\Delta,\bullet}(Y_+).
\end{equation*}
Here, we sum up all \RRSC\ structures on $W$ and $Y_{\pm}$, whenever the composition makes sense.
In future notation, we suppress the marking data if no confusion arises.

\begin{prop}
Fix $a_*$.
The chain map \eqref{eq:cone-map-summary} is independent, up to canonical isomorphisms,
of the metrics, the perturbations, and the chain-level representatives of $\upsilon$.
In particular, the cone Floer homology is independent of metrics and perturbations.
\end{prop}

\begin{proof}
First, we fix the metrics, perturbations, and chain representatives of $\upsilon_\pm$ on the 3-manifolds $Y_{\pm}$.
Given two choices of metrics and perturbations of $\upsilon$ on $W$, let
$\widetilde m^{(0)},\widetilde m^{(1)},\widecheck m^{(0)},\widecheck m^{(1)}$ be the corresponding maps as in~\eqref{eq:cone-map-summary}.
An interpolation between the metrics and perturbations on the unpunctured cobordism defines a homotopy $H$, for which
\[\check{\partial}_+ H + H \check{\partial}_- = \widecheck m^{(1)} + \widecheck m^{(0)}.\]
Combining this family with the puncture-moving family of metrics, as in the definition of $\check K$, gives rise to a two-dimensional family of moduli spaces~\cite[\S2.2]{Li24}.
Counting zero-dimensional moduli spaces in this family defines an operator $\check L$.
The boundary of 1-dimensional family moduli spaces gives the identity:
\[\check{\partial}_+ \widecheck L + \widecheck L\check{\partial}_- = \check K^{(1)} -\check K^{(0)} + \upsilon_+ H + H\upsilon_-.\]
Therefore
\[
\widetilde m^{(1)}+\widetilde m^{(0)}
=\widetilde\partial_+\begin{bmatrix}H&0\\\widecheck L&H\end{bmatrix}
+\begin{bmatrix}H&0\\\widecheck  L&H\end{bmatrix}\widetilde\partial_-.
\]

Next, we vary the metrics and perturbations on $Y_\pm$.
Puncture the cylinder $I \times Y_{\pm}$ used in the corresponding continuation map.
The diagonal entries of the resulting cone cobordism maps induce isomorphisms the \HMR-\emph{to} of $Y_{\pm}$ by \cite[Corollary~13.5]{Li22}.
Thus the induced maps on the homology of the mapping cones are also isomorphisms.
The neck-stretching argument in the proof of Proposition~\ref{prop:composition-law} below shows that these isomorphisms compose compatibly with the cobordism maps.

Finally, we vary the chain representatives of $\upsilon_\pm$.
Interpolation between the representatives gives rise to chain homotopies between $\upsilon_{\pm}$, and hence quasi-isomorphisms between their mapping cones.
A similar argument as above shows that the cobordism maps agree up to chain homotopy under these isomorphisms.
\end{proof}

\begin{figure}[htbp]
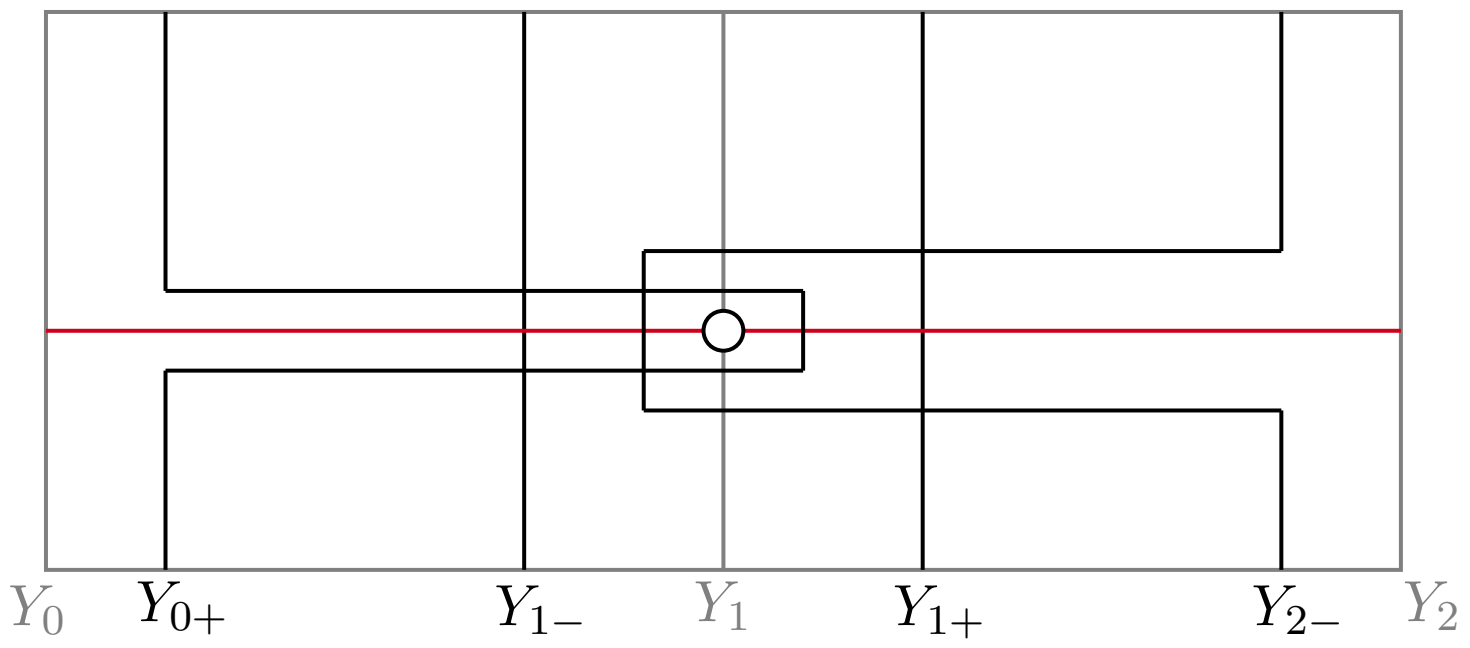

	\centering
	\figfourhypersurfaces
	\caption{Four hypersurfaces $\{Y_{0+},Y_{1-},Y_{1+},Y_{2-}\}$ in $W_1 \circ W_0$, with a puncture along $a_{*1} \cup a_{*0}$.}
	\label{fig:fourhypersurfaces}
\end{figure}

\begin{prop}\label{prop:composition-law}
	Suppose $W_0 \colon Y_0 \to Y_1$ and $W_1 \colon Y_1 \to Y_2$ are real marked cobordisms,
	together with distinguished arcs $a_{*0} \subset \mathbf{a}_0$ and $a_{*1} \subset \mathbf{a}_{1}$.
	Suppose also the outgoing endpoint of $a_{*0}$ coincides with the incoming endpoint of $a_{*1}$.
	Then
	\begin{equation*}
		\thmr_{\Delta,a_{*1} \cup a_{*0}}(W_1 \circ W_0) = \thmr_{\Delta,a_{*1}}(W_1) \circ \thmr_{\Delta,a_{*0}}(W_0).
\end{equation*}
\end{prop}

\begin{proof}
We give a proof sketch, as~\cite[\S8]{Bloom} explained a more general neck-stretching argument needed in this proposition.
Let $\widecheck m_j,\check K_j$ be the operators defining the cone cobordism maps of $W_j$,
and $\widecheck m,\check K$ be those of $W_1 \circ W_0$ with arc $a_{*0}\cup a_{*1}$.
Let $\check\partial_j$ and $\upsilon_j$ be the corresponding differential and V-map for $Y_j$.
The ordinary version of the composition law, using a 1-parameter family of moduli space, defines a homotopy $H$:
\[\widecheck\partial_2 H+H\widecheck\partial_0=\widecheck m+\widecheck m_1\widecheck m_0.\]

We define a 2-dimensional family of metrics that combines stretching along $Y_1$ and moving the puncture,
by stretching pairs of disjoint hypersurfaces in Figure~\ref{fig:fourhypersurfaces}.
The boundary of the resulting parameter space consists of five edges: the puncture-moving family on $W_1 \circ W_0$, the two families where the punctures are stretched off the $Y_0,Y_2$ ends, and the two families broken along $Y_1$, with puncture inside $W_0$ or $W_1$.
Counting zero-dimensional families defines an operator $\check L$, satisfying the following identity:
\[\widecheck\partial_2 \widecheck L+\widecheck  L\widecheck\partial_0
      =\widecheck K+\widecheck K_1\widecheck m_0+\widecheck m_1\widecheck K_0+\upsilon_2H+H\upsilon_0.\]
Define $\tilde \partial_j$ as the cone differential on $Y_j$.
Then the above identities imply
\[\begin{bmatrix}
	\widecheck m & 0\\\widecheck K & \widecheck m
\end{bmatrix}
+
\begin{bmatrix}
	\widecheck m_1 & 0\\\widecheck K_1 & \widecheck m_1
\end{bmatrix}
\begin{bmatrix}
	\widecheck m_0 & 0\\\widecheck K_0 & \widecheck m_0
\end{bmatrix}=
\widetilde{\partial}_2
\begin{bmatrix}
	H & 0\\\widecheck L & H
\end{bmatrix}
+
\begin{bmatrix}
	H & 0\\ \widecheck L & H
\end{bmatrix}\widetilde{\partial}_0.\]
The cone-version of composition law follows.
\end{proof}

\section{Exact triangles}
\label{sec:exact-triangle}
In this section,
we prove oriented exact triangles.
The branched covers of the skein cobordisms of $L_+$ are surgeries on $\dbcv_2(S^3,L_+)$,
which are equivariant in the sense that the surgery curves come in $\tau$-exchanged pairs.
This is in contrast to the unoriented skein triangle in \cite{Li23},
where surgeries are performed along $\tau$-invariant circles.
We adapt the ordinary exact triangle~\cite{KMOS} to prove the \emph{to}-version exact triangle in diagonal \hmr.
Combining this triangle with the arguments in \cite{Li24} (based on~\cite{Bloom}), we establish an exact triangle in the cone version.
\begin{figure}[htbp]
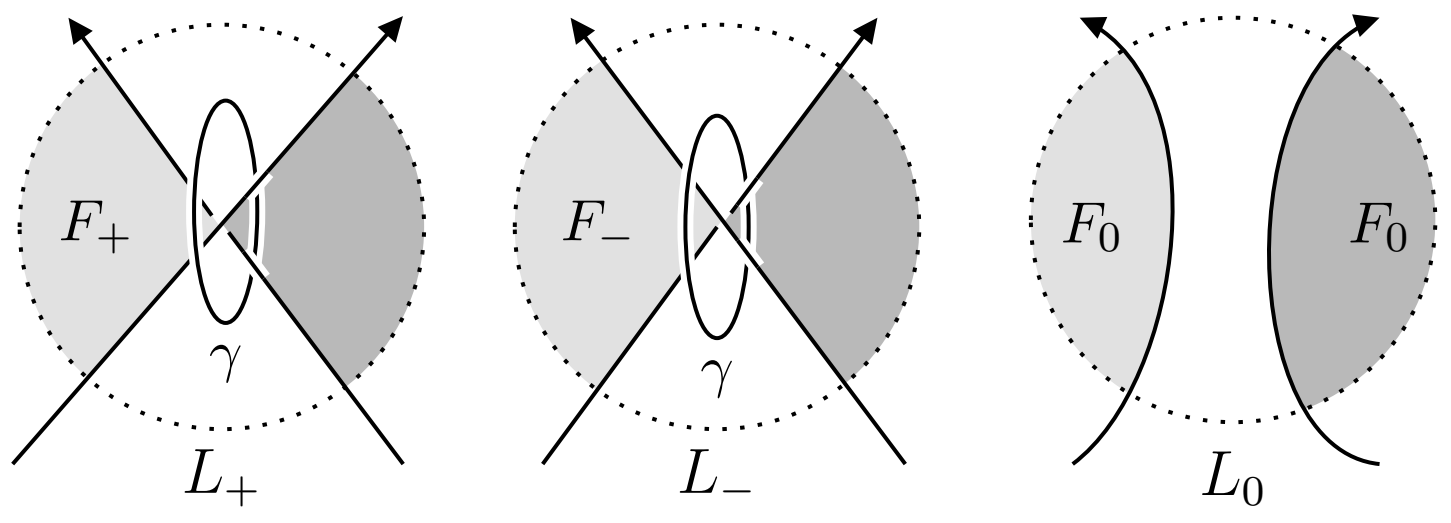

	\centering
	\figskeinF
	\caption{An oriented skein triple, with shaded Seifert surfaces.}
	\label{fig:skeinF}
\end{figure}
\subsection{Topology of the oriented skein triple}
\label{ssec:top-or-triple}
Let $(L_+,L_-,L_0)$ be an oriented skein triple,
and assume $(F_+, F_-, F_0)$ are Seifert surfaces.
In particular, $F_-$ is obtained from a full-twist of the band $F_+ \cap B^3$,
where $B^3$ is the local 3-ball around the crossing.
And $F_0$ is obtained from deleting the band, as in Figure~\ref{fig:skeinF}.
The crossing is \emph{mixed} when the two arcs of $L_+ \cap B^3$ lie on different components of $L_\pm$, and it is a \emph{self} crossing otherwise.
In other words,
$|L_0|=|L_\pm|-1$ in the mixed case, and $|L_0|=|L_\pm|+1$ in the self case.
See Figure~\ref{fig:selfmixed}.
Denote the double branched covers of the oriented skein triple as
\begin{equation*}
	Y_+=\dbcv_2(S^3,L_+),\quad
	Y_-=\dbcv_2(S^3,L_-),\quad
	Y_0 = \dbcv_2(S^3, L_0).
\end{equation*}
Let $\tau_+,\tau_-,\tau_0$ be the corresponding deck transformations.

\begin{figure}[htbp]
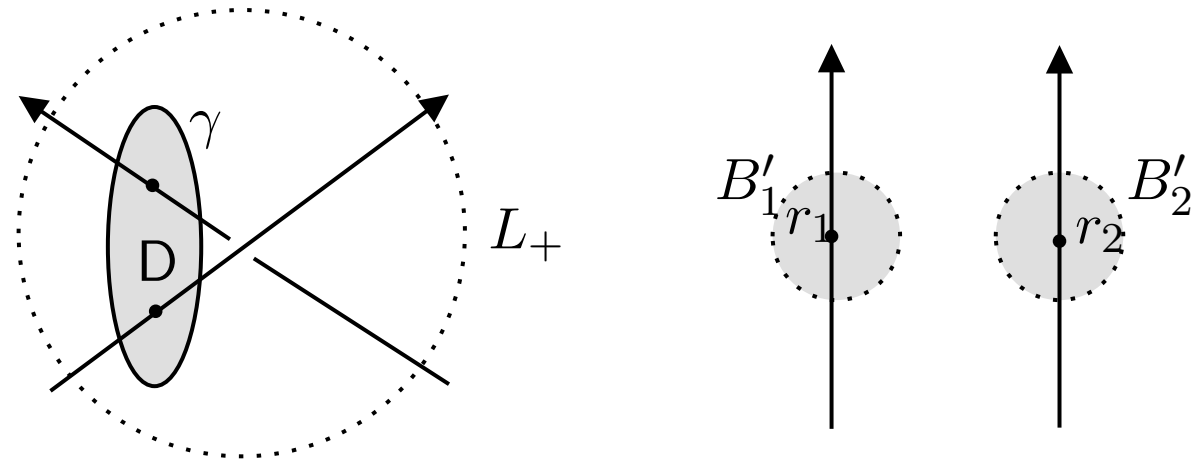

	\centering
	\figdiscandball
	\caption{The disc $\mathsf{D}$ bounding $\gamma$ (left) and the two 3-balls on $L_0$ (right), both shaded.}
	\label{fig:discandball}
\end{figure}

Let $\gamma$ be a circle around the crossing inside $B^3$,
which bounds a 2-disc $\mathsf{D}$ that intersects $L_+$ twice with opposite signs, as in Figure~\ref{fig:discandball}.
Assume that the $-1$ surgery along $\gamma$ turns $L_+$ into $L_-$.
Then the three Dehn fillings of $S^3 \setminus \nu(\gamma)$ with slopes $\infty,-1,0$ are $S^3,S^3$ and $S^3_0(\gamma)$, respectively.
The preimage of $\gamma$ in $Y_+$ is a pair of circles swapped by $\tau_+$.
Therefore, the double cover of these three surgeries are $Y_+$, $Y_-$, and $Y_{00}$,
where
\begin{equation*}
	Y_{00}=\dbcv_2(S^3_0(\gamma),L_+;[F_{0}] \bmod 2).
\end{equation*}
The above notation indicates taking the double branched cover determined by the mod-two class of the Seifert surface $[F_0]$.
In other words, $Y_+,Y_-,Y_{00}$ are related by $\tau$-equivariant (double) surgeries.

The manifold $Y_{00}$ can be obtained from $Y_{0}$ as follows.
Let $B_1', B_2'$ be two 3-balls in $S^3$ around two points $r_1, r_2$ on the two arcs of $L_0$, as in Figure~\ref{fig:discandball}.
Let $B_i$ be the preimage of $B_i'$ under $Y_0 \to S^3$.
The involution $\tau_0$ restricts to the 3-ball $B_i$ as a $180^{\circ}$-rotation around an axis, so $\tau|_{\partial B_i}$ fixes two points.
Denote such a rotation on $B^3$ or $S^2$ by $\rot$.
Define an involution on ${S^2\times[1,2]}$ by $\rot\times\id$.
Then
\begin{equation}\label{eq:tube}
	Y_{00}=\bigl(Y_0\setminus\operatorname{int}(B_1 \sqcup B_2)\bigr)
	\cup_{
		\partial B_i = S^2 \times \{i\}}(S^2\times[1,2]),
\end{equation}
In other words, $Y_{00}$ is a boundary component of the cobordism obtained by attaching a 4-dimensional 1-handle to $I\times Y_0$.
Conversely, cap off the disc $\mathsf{D}$ with the meridional disc of $S^1 \times D^2$ to obtain a 2-sphere $\Sigma' \subset S^3_0(\gamma)$.
Then the preimage $\Sigma$ of $\Sigma'$ in $Y_{00}$ is a 2-sphere with involution $\rot$.
Cutting $Y_{00}$ along $\Sigma$ and filling the two resulting boundaries by $(B^3, \rot)$ recovers $Y_0$.

Let $\mathbf{p}$ be a set of basepoints away from the local 3-ball, with one point on each component of $L_+$.
In particular, $Y_+,Y_-,Y_{00}$ all satisfy Hypotheses~\ref{hyp:abs-grad}.
Furthermore, they are related by traces of the double surgeries:
\[
W_+\colon Y_+\longrightarrow Y_-,\qquad
W_- \colon Y_-\longrightarrow Y_{00},\qquad
W_{00} \colon Y_{00}\longrightarrow Y_+,
\]
respectively.
In other words, each 3-manifold above is obtained by attaching $\tau$-exchanged pairs of $2$-handles to $I \times Y_{\bullet}$, away from the fixed sets $I \times \fix(\tau_{\bullet})$.
Here, $\bullet \in \{+,-,00\}$.
Since $\mathbf{p}$ is assumed to be away from $\gamma$, the marking set in each $W_{\bullet}$ is simply $\mathbf{a} = I \times \mathbf{p}$.
Hence the involutions on $I \times Y_{\bullet}$ extend to $W_{\bullet}$,
making $W_{\bullet}$ a real marked cobordism.

In what follows, let
\[Y_1 = Y_+, \quad Y_2 = Y_-,\quad Y_3 = Y_{00},\quad W_1 = W_+,\quad W_2 = W_-,\quad W_3 = W_{00}\]
extend these labels three-periodically.
Let $\tau_i$ denote the corresponding involutions, and $\widecheck F_i,\widetilde F_i$ denote the corresponding cobordism maps.
The surgeries of $\gamma \sqcup \tau \gamma$ take place in the trivially covered region.
Thus we have the following doubled version of the description in \cite[\S5.2]{KMOS} about the topology of composite cobordisms:

(i) The double composition
\[X_i = W_{i+1}\circ W_i\colon Y_i\to Y_{i+2}\]
contains a pair of $2$-spheres $E_i\sqcup \tau E_i$ of self-intersection $-1$,
exchanged by the involution.
Each $E_i$ is disjoint from the marking set,
and $\tau$-equivariantly:
\[W_{i+1}\circ W_i\;\cong\;
B_i\,\#\,\overline{\mathbb{CP}}^{\,2}\,\#\,\overline{\mathbb{CP}}^{\,2}.\]
Here, $\tau$ swaps the $\overline{\mathbb{CP}}^{\,2}$ summands,
and $B_i \colon Y_i \to Y_{i+2}$ is the surgery trace for the surgery on the same attaching circle as $W_i$, but with framing increased by one.
Let $Z_i$ be a neighbourhood of $E_i$, diffeomorphic to $\overline{\mathbb{CP}}^{\,2}\setminus B^4$.
Moreover,
let $B_i^\circ=B_i\setminus(B^4\sqcup\tau B^4)$,
where $B^4 \sqcup \tau B^4$ are the connected sum region above.
Then
\[X_i = B_i^\circ\cup_{S^3\sqcup\tau S^3}(Z_i\sqcup\tau Z_i).\]

(ii) The triple composition
\[V_i = W_{i+2}\circ W_{i+1}\circ W_i\colon Y_i\to Y_i\]
contains two pairs of 2-spheres $E_i\sqcup\tau E_i$ and
$E_{i+1}\sqcup\tau E_{i+1}$ for which $E_i\cdot E_{i+1}=1$.
These four spheres are pairwise disjoint otherwise.
A regular neighbourhood $N_i$ of
$E_i\cup E_{i+1}$ is diffeomorphic to
$\overline{\mathbb{CP}}^{\,2}\setminus\nu(S^1)$,
for which $\partial N_i= R \cong S^1\times S^2$.
Under this identification, $[\{\mathrm{pt}\}\times S^2]=E_i+E_{i+1}$.
The complement
$V_i \setminus(N_i\sqcup\tau N_i)$ is $\tau$-equivariantly diffeomorphic to $I\times Y_i$ with the tubular neighbourhoods of a pair of circles removed.
Thus
\[(V_i \setminus (N_i \cup \tau N_i))\cup_{R\sqcup\tau R}
\bigl((S^1\times B^3)\sqcup\tau(S^1\times B^3)\bigr)\;\cong\;I\times Y_i .\]

Similarly, we can describe relative real spin\textsuperscript{c} structures on the composites:
\begin{lem}\label{lem:spinc-composites}
We have the following bijections:
\begin{enumerate}[label=(\roman*)]
	\item \[\RRSC(W_{i+1}\circ W_i)\cong\RRSC(B_i^\circ)\times\mathbb{Z},\]
sending $\mathfrak{t}$ to $(\mathfrak{t}|_{B_i^{\circ}},k)$,
for which $\langle c_1,[E_i]\rangle=2k+1$;
	\item furthermore,
\[\left\{\mathfrak{t} \in \mathrm{Spin}^c(N_i): \mathfrak{t}|_{\partial N_i}\, \text{is torsion} \right\} \cong \mathbb{Z},\]
sending $\mathfrak{t}$ to $k$,
for which $\langle c_1,[E_i]\rangle=2k+1=-\langle c_1,[E_{i+1}]\rangle$.
\end{enumerate}
In both cases, conjugation of \text{spin\textsuperscript{c}} structures gives $k\mapsto-1-k$.
\end{lem}

\begin{proof}
On a pair $M \sqcup M$ of identical manifolds for which the involution $\tau$ acts by swapping the two connected components, a real spin\textsuperscript{c} structure on $M \sqcup M$ is equivalent to an ordinary spin\textsuperscript{c} structure $\mathfrak{t}$.
The real structure provides a conjugate-linear identification of the spin\textsuperscript{c} structures on the two copies of $M$.
Thus $\tau^*\mathfrak{t} \cong \bar{\mathfrak{t}}.$
This applies to $B_i^{\circ} \sqcup \tau B_i^{\circ}$ and $N_i \sqcup \tau N_i$, so (i) and (ii) follow from the ordinary counterpart as in~\cite[\S5.1-5.2]{KMOS}.
\end{proof}

\subsection{Mixed operations}
A \emph{mixed operation} refers to an operation on real monopole Floer homology that involves the ordinary Seiberg--Witten configurations.
The local computations in the proof of exact triangle all involve $\tau$-interchanged pairs of manifolds.
By the following lemma, the real monopole Floer theory on these pairs reduces to the ordinary theory of \cite{KM,KMOS}.
\begin{lem}\label{lem:exchanged-pairs}
Let $M$ be an oriented $3$- or $4$-manifold that is either compact or with cylindrical ends.
Suppose $\tau$ is the involution acting on $M \sqcup M$ by exchanging the two copies.
Choose the real structure $\sr$ induced by the identification of the spinor bundle on $M$ and $\tau M$, so that $\tau^*\mathfrak{s}|M= \overline{\mathfrak{s}}|_{\tau M}$.
Then:
\begin{enumerate}[leftmargin=*,label=(\roman*)]
\item Restriction to $M$ defines a bijection between real \textup{spin}$^c$ structures on $M \sqcup M$ and ordinary \textup{spin}$^c$ structures on $M$.
\item Furthermore, the restriction identifies the real configuration space, the real gauge group, and the perturbed Seiberg--Witten operators on $M \sqcup M$, with the ordinary configuration space, gauge group, and perturbed Seiberg--Witten operators on $M$.
This identification intertwines with the blow-up construction,
up to the rescaling $\|\Phi\|_{L^2(M\sqcup M)}=\sqrt2\,\|\Phi|_M\|_{L^2(M)}$.
Thus the restriction also identifies all the moduli spaces of (perturbed) Seiberg--Witten solutions.

\item Finally, when $M$ is a 4-manifold, the restriction identifies the real monopole chain complexes on $\partial M\sqcup \partial M$ with the corresponding ordinary complexes of $\partial M$.
Under these identifications,
the spectral projections and gluing maps along
$\partial M \sqcup \partial M$ coincide with the ordinary counterparts on $\partial M$.
\end{enumerate}
\end{lem}

\begin{proof}
Statement (i) has been proved in the proof of Lemma~\ref{lem:spinc-composites}.
Statement (ii) is clear by the definitions of real configuration spaces and gauge groups.
	For (iii), the gluing theorems in
	\cite[\S10, \S12--13]{Li22} are formulated for arbitrary
	$\tau$-invariant 3-manifolds.
	Hence they apply to the boundaries 
	$\partial M\sqcup\tau\partial M$,
	with gluing data identified with the ordinary data on
	$\partial M$ thanks to (ii).
\end{proof}

Let $(W,\tau_W,\mathbf{a}) \colon (Y_-,\tau_-, \mathbf{p}_-) \to (Y_+,\tau_+,\mathbf{p}_+)$ be a real marked cobordism, with two additional incoming ends $R \sqcup R$ exchanged by $\tau_W$.
For a real \text{spin\textsuperscript{c}} structure on $R \sqcup R$, identified with a \text{spin\textsuperscript{c}} structure on $R$.
Write $\widehat{C}_{\bullet}(R)=\widehat{C}_{\HM,\bullet}(R)$.
Then counting zero-dimensional moduli spaces on $W$ defines the following chain map
\begin{equation*}
	\mu_W \colon \widehat{C}_{\bullet}(R) \otimes \widecheck{C}_{\Delta}(Y_-,\tau_-,\mathbf{p}_-) \to \widecheck{C}_{\Delta}(Y_+,\tau_+,\mathbf{p}_+).
\end{equation*}
The combinatorics of the boundary degeneration~\cite{Bloom13} (see also \cite[\S4.3]{Lin18}),
together with Lemma~\ref{lem:exchanged-pairs} shows that this is a well-defined map on homology.
By Lemma~\ref{lem:exchanged-pairs}(iii), $\mu_W$ commutes with point operation, i.e.,
$\mu_W(z\otimes\upsilon x)$ and $\upsilon\mu_W(z\otimes x)$ are chain
homotopic, by moving the puncture through the connected cobordism $W$.

We extend the gluing setup in \cite{Li22} to the case where the involution can swap pairs of 3-manifolds as follows.
Let $(W,\tau_W)$ be a real cobordism whose boundary components are either $\tau_W$-invariant or pairs of two disconnected components that are swapped by $\tau_W$.
Fix admissible perturbations on some collar neighourhoods of the boundaries of $W$.
\begin{lem}
The boundary stratification \cite[Theorem~12.16]{Li22} holds for the moduli spaces on $W$.
Furthermore, the real analogue of \cite[Proposition~26.1.6]{KM} holds for moduli spaces over neck-stretching families of metrics.
\end{lem}

\begin{proof}
Let $R \sqcup R$ be a boundary component interchanged by $\tau_W$.
By Lemma~\ref{lem:exchanged-pairs}, the configuration spaces, gauge groups, and Seiberg--Witten equations on $R \sqcup R$ can be identified with their restriction on a single $R$.
The moduli space of $W$ is regular, by the transversality argument in \cite[Proposition 24.4.7]{KM}.
The family version follows from the argument in \cite[Proposition 12.10]{Li22}.
The gluing construction follows from \cite[\S24.7 and \S26.1]{KM} and Lemma~\ref{lem:exchanged-pairs},  using one copy of the interchanged pair.
This includes the fibre-product description and constructions of $\delta$ and $\mathbf{S}$ as in~\cite[\S24.7]{KM}.
\end{proof}

\subsection{Statement of the theorems}

For $i\in \{1,2,3\}$, denote the total diagonal \emph{to}-complex over all \rrsc\ structures as
\[\widecheck C_i=\widecheck{\cmr}_\Delta(Y_i,\tau_i,\mathbf p)=
\bigoplus_{\widetilde{\sr} \in \RRSC(Y_i,\tau_i,\mathbf{p}_i)}\widecheck{\cmr}_\Delta(Y_i,\tau_i,\mathbf p;\widetilde{\sr}),\]
with differential $\check\partial_i$.
Represent the internal $V$-map $\upsilon_i$ on $Y_i$ using the punctured cylinder in~\S\ref{ssec:cob-maps}.
Denote the cone complex as
\[\widetilde C_i=\widecheck C_i\oplus \widecheck C_i,\quad
\widetilde{\partial}_i=\begin{pmatrix}\widecheck\partial_i&0\\ \upsilon_i&\widecheck\partial_i\end{pmatrix}.\]

Let $\check F_i=\widecheck m_\Delta(W_i)$ denote the to-cobordism map, summed over all \RRSC\
structures on $W_i$.
Moreover, let $\check K_i$ be the point-moving operator as in \S\ref{ssec:cob-maps}, so that the cone cobordism map of $(W_i, I \times p_*)$ is given by the matrix
\[\widetilde F_i=\begin{bmatrix}\check F_i&0\\\check K_i&\check F_i\end{bmatrix}.\]
Here, $p_*$ is any point in $\mathbf{p}$, away from the local 3-ball in the skein relation, which we fix once and for all.
This is needed for the cone-cobordism maps, but does not affect the statement of the cone version of the exact triangle.
Henceforth we will suppress $p_*$ in notation.
Moreover, in the next two theorems, we use the marking set $I \times \mathbf{p}$ throughout.

\begin{thm}\label{thm:exact-triangle}
	The cobordism maps $\tilde{F}_1,\tilde{F}_2,\tilde{F}_3$ induced by $W_1,W_2,W_3$, respectively, define an exact triangle
	\begin{equation}\label{eq:exact-triangle}
		\begin{tikzcd}
			\tH_\Delta(Y_1,\tau_1,\mathbf p)\arrow[rr,"\widetilde F_1"] &&
			\tH_\Delta(Y_2,\tau_2,\mathbf p)\arrow[dl,"\widetilde F_2"]\\
			& \tH_\Delta(Y_3,\tau_3,\mathbf p)\arrow[ul,"\widetilde F_3"] &
		\end{tikzcd}
	\end{equation}
\end{thm}
\begin{thm}\label{thm:exact-triangle-to}
The cobordism maps $\check{F}_1,\check{F}_2,\check{F}_3$ induced by $W_1,W_2,W_3$, respectively, define an exact triangle
\begin{equation*}
	\begin{tikzcd}
		\toH_\Delta(Y_1,\tau_1,\mathbf p)\arrow[rr,"\widecheck F_1"] &&
		\toH_\Delta(Y_2,\tau_2,\mathbf p)\arrow[dl,"\widecheck F_2"]\\
		& \toH_\Delta(Y_3,\tau_3,\mathbf p)\arrow[ul,"\widecheck F_3"] &
	\end{tikzcd}
\end{equation*}
\end{thm}

Theorem~\ref{thm:exact-triangle-to} will be proved in the course of proving Theorem~\ref{thm:exact-triangle}.
We use the standard triangle detection lemma (see e.g.~{\cite[Lemma~5.1]{KMOS}}).

\begin{lem}
	Let $(C_i,\partial_i)_{i\in\mathbb{Z}/3}$ be chain complexes over $\F$, let
	$f_i\colon C_i\to C_{i+1}$ be chain maps, and let
	$H_i\colon C_i\to C_{i+2}$ be a chain homotopy of the following significance:
	\begin{enumerate}[leftmargin=*,label=(\alph*)]
		\item $\partial_{i+2}H_i+H_i\partial_i=f_{i+1}f_i$, and
		\item the maps $\psi_i:=f_{i+2}H_i+H_{i+1}f_i\colon C_i\to C_i$ induce
		isomorphisms on homology.
	\end{enumerate}
	Then the following sequence	is exact:
\[\cdots\longrightarrow H(C_i, \partial_i)\xrightarrow{(f_i)_*}H(C_{i+1},\partial_{i+1})
\xrightarrow{(f_{i+1})_*}H(C_{i+2},\partial_{i+2})\longrightarrow\cdots.\]
\end{lem}

To prove Theorem~\ref{thm:exact-triangle}, we plug in $C_i=\tilde C_i$, $\partial_i = \tilde \partial_i$, and $f_i=\tilde F_i$.
Similarly, for Theorem~\ref{thm:exact-triangle-to}, we plug in $C_i=\check C_i$, $\partial_i = \check \partial_i$, and $f_i=\check F_i$.

We construct the homotopies $\widetilde H_i$ and verify (a) in \S\ref{ssec:null-homotopies}.
Condition (b) is proved in \S\ref{ssec:detector}.
The homotopies $\check H_i$ and Conditions (a), (b) will be verified along the way.
The structure of the proof and notation follow \cite[\S5]{KMOS}.

\subsection{The null-homotopies}
\label{ssec:null-homotopies}
\begin{figure}[htbp]
	\centering
	\figtwofivehypersurface
	\caption{Five-tuples of hypersurfaces used in Lemmas~\ref{lem:null-homotopy-cone} and~\ref{lem:detector}. Compare with~\cite[Figure~2]{Li24} and~\cite[Figure~1-2]{KMOS}}
	\label{fig:twofivehyp}
\end{figure}
The double composite $X_i = W_{i+1}\circ W_i$ contains hypersurfaces $S_i \sqcup \tau S_i$ and $Y_{i+1}$, where $S_i = \partial Z_i$.
See Figure~\ref{fig:twofivehyp}.
Let $Q = [-\infty,+\infty]$ parametrize a one-parameter family of metrics on $X_i$ as in \cite[\S5.1]{KMOS},
where the neck along $Y_{i+1}$ is stretched for $T \in [-\infty,0]$,
and the neck of $S_i \sqcup \tau S_i$ is stretched for $T \in [0,+\infty]$.
Choose
a regular family of real perturbations over $Q$ by~ \cite[Proposition~24.4.7]{KM} and \cite[Proposition~12.10]{Li22}.
Indeed, the perturbation data on $Z_i \cup \tau Z_i$ reduce to the ordinary perturbations by Lemma~\ref{lem:exchanged-pairs}.

Let $H^a_b$ and $\bar H^a_b$, $a,b\in\{o,s,u\}$, be the counts of the zero-dimensional family moduli spaces and reducible moduli spaces over $Q$, respectively,
summed over all \RRSC\ structures on $X_i$.
Denote by $m^a_{b}(X_i),\bar m^a_{b}(X_i)$ the entries of $\widecheck m(W_i)$.
Define the following operator using the identical formula as in \cite[\S5.1]{KMOS}:
\begin{equation}
\label{eq:check-H-matrix}
	\widecheck H_i=\begin{bmatrix}
H_o^o & H_o^u\bar\partial_u^s+\partial_o^u\bar H_u^s
       +m_o^u(W_{i+1})\bar m_u^s(W_i)\\
H_s^o & \bar H_s^s+H_s^u\bar\partial_u^s+\partial_s^u\bar H_u^s
       +m_s^u(W_{i+1})\bar m_u^s(W_i)
	\end{bmatrix}
	\colon \widecheck C_i\longrightarrow \widecheck C_{i+2}.
\end{equation}

\begin{lem}\label{lem:h-identity}
	$\widecheck{\partial} \widecheck H_i+ \widecheck H_i \widecheck{\partial}_i=\widecheck F_{i+1} \widecheck F_i$.
\end{lem}

\begin{proof}
The structure of proof is identical to that of \cite[Proposition~5.2]{KMOS}, thanks to Lemma~\ref{lem:exchanged-pairs}.
In particular, \cite[Proposition~5.2]{KMOS} describes the boundary of the compactified one-dimensional parametrized moduli spaces over $Q$,
	which by Lemma~\ref{lem:exchanged-pairs}(ii) implies the real version.
	At $T=-\infty$, the degeneration of the solutions on $X_i$ along
	$Y_{i+1}$ can be analyzed using the \cite[Theorems~10.7--10.8 \&
	\S12]{Li22}.
	At the hypersurface $S^3\sqcup\tau S^3$, the gluing theorem reduces to the ordinary version~\cite[\S19]{KM} on $S^3$ by Lemma~\ref{lem:exchanged-pairs}(iii).	
	Finally,  by Lemma~\ref{lem:exchanged-pairs}(ii) and Lemma~\ref{lem:spinc-composites}(i),
	the contributions from $Z_i \sqcup\tau Z_i$ in the $k$-th spin\textsuperscript{c} structure add up to an element in $\widehat C_{\bullet}(S^3)$.
	By \cite[Lemma~5.3]{KMOS} the terms from $\mathfrak s_k$ and
	$\mathfrak s_{-1-k}$ coincide,
	on which the conjugation symmetry $k\mapsto-1-k$ acts freely.
	Therefore the contribution on the $T= \infty$ end, summed over all \RRSC\ structures on $X_i$, vanishes mod two.
\end{proof}

Following \cite{Bloom} and \cite[\S4]{Li24}, we upgrade $\widecheck H_i$ to the cone version.
This is achieved by a family of invariant metrics of one dimension higher,
obtained by dragging the puncture along the arc $a_*$ that we have fixed.
By definition, the arcs $a_* = I \times p_*$ concatenate along the cobordisms.
Let $p \in a_*$, and puncture $X_i$ equivariantly at $p$.
Over every metric in $Q$, move the puncture along $a_*$ by introducing two new hypersurfaces $Y_{i,+}, Y_{i+2,-}$ as in~\S\ref{ssec:cob-maps}.
See Figure~\ref{fig:twofivehyp} for the five hypersurfaces $\{Y_{i,+}, S_i \sqcup \tau S_i, Y_{i+2,+}, Y_{i+1,-}, Y_{i+1,+}\}$ used in the neck stretching.
This results in a pentagon-family $\tilde{Q}$ of metrics with five boundary edges, corresponding to degeneration of the metrics, that are in cyclic order:
\begin{enumerate}[leftmargin=*,label=(\arabic*)]
	\item the puncture is stretched off to the incoming end $I \times Y_i$;
	\item the metric is broken along $S_i \sqcup \tau S_i$.
		\item the puncture is stretched off to the outgoing end $I \times Y_{i+2}$;
	\item the metric is broken along $Y_{i+1,-}$  and the puncture lies in $W_i$;
	\item the metric is broken along $Y_{i+1,+}$ and the puncture lies in $W_{i+1}$.
\end{enumerate}

Let $I_i\colon \widecheck C_i\to \widecheck C_{i+2}$ be the operator counting
the zero-dimensional family moduli spaces over $\widetilde{Q}$
with $\mathfrak c_{-2}$ inserted at the puncture,
and with boundary-obstructed terms at edges~(3) and (4),
as in the upper-right entry of \eqref{eq:check-H-matrix}.
See also \cite[Lemma~8.2]{Bloom}.
The cone-version of the null homotopy is given by the matrix:
\begin{equation*}
	\widetilde H_i = \begin{bmatrix} \check H_i&0\\I_i&\check H_i\end{bmatrix}
	\colon\widetilde C_i\longrightarrow\widetilde C_{i+2}.
\end{equation*}

\begin{lem}\label{lem:null-homotopy-cone}
$\widetilde{\partial} \widetilde{H}_i+ \widetilde{H}_i \widetilde{\partial}_i=\widetilde F_{i+1}\widetilde F_i$.
\end{lem}

\begin{proof}
This follows from
\[\partial I_i+ I_i \check\partial = \upsilon_{i+2}\check H_i+ \check H_i\upsilon_i+\check K_{i+1} \check F_i+ \check F_{i+1} \check K_i.\]
The proof goes by counting the boundaries of the one-dimensional moduli spaces over $\widetilde{Q}$, as in \cite{Bloom} and \cite[Lemma~4.9]{Li24}.
To sketch the proof, the term $\check \partial I_i+I_i \check\partial_i$ arises from degeneration at the $Y_i,Y_{i+2}$ ends.
Edges~(1) and~(3) contribute $\check H_i\upsilon_i$ and $\upsilon_{i+2}\check H_i$.
Edges~(4) and~(5) contribute $\check F_{i+1} \check K_i$ and $\check K_{i+1}\check F_i$.
The auxiliary $S^3$-end at the puncture contributes zero as the Floer differential is zero mod two.
Also, edge~(2) contributes zero due to the same symmetry argument as in Lemma~\ref{lem:h-identity}.
Finally, all four entries of the matrix $\tilde{\partial} \tilde{H}_i+ \tilde{H}_i \tilde{\partial}_i-\tilde F_{i+1}\tilde F_i$ vanish, by the above and Lemma~\ref{lem:h-identity}.
\end{proof}

\subsection{Proof of exact triangles}
\label{ssec:detector}
Denote $\psi_i=\widecheck F_{i+2}\widecheck H_i+ \widecheck H_{i+1}\widecheck F_i$
and $\widetilde \psi_i =\widetilde F_{i+2}\widetilde H_i+ \widetilde H_{i+1}\widetilde F_i$.
The latter satisfies:
\begin{equation*}
	\widetilde\psi_i
	=\begin{bmatrix}\psi_i&0\\
	\check K_{i+2}\check H_i+\check F_{i+2}I_i+ I_{i+1}\check F_i+\check H_{i+1} \check K_i&\psi_i\end{bmatrix}.
\end{equation*}
Lemma~\ref{lem:null-homotopy-cone} implies $\tilde\psi_i$ is a chain map.
Moreover, since $\tilde\psi_i$ is lower
triangular, it induces a map of the short exact sequences
$0\to \check C_i\to\widetilde C_{i}\to \check C_i \to 0$ as in \eqref{eq:v-SES}, and thus a map $(\psi_i)_*$ on both $\check C_i$ terms.
Inside the triple composite cobordism $V_i = W_{i+2} \circ W_{i+1} \circ W_{i}$,
there are five separating hypersurfaces, ordered cyclically as:
\[S_i \sqcup \tau S_i, \quad Y_{i+2}, \quad Y_{i+1},\quad S_{i+1} \sqcup \tau S_{i+1}, \quad R \sqcup \tau R.\]
Here, each consecutive pair is disjoint, and has nontrivial intersection otherwise,
as in Figure~\ref{fig:twofivehyp}.
Following \cite[\S5.2]{KMOS}, there is a pentagon $P$ of Riemannian metrics, obtained from simultaneously stretching pairs of disjoint hypersurfaces.
At one of the edges, $Q(R)$, the metric is broken along $R \cup \tau R$,
where by \S\ref{ssec:top-or-triple} the cobordism $V_i$ decomposes as:
\begin{equation}
\label{eq:Vi-decomp-cylinder}
V_i = ([-1,+1] \times Y_i) \setminus \nu(\{0\} \times (K_i \sqcup \tau K_i)) \cup (N_i \sqcup \tau N_i).
\end{equation}
Here, $K_i\sqcup\tau K_i\subset Y_i$ are the cores of the two solid tori in the Dehn filling defining $Y_i$.
And $\nu(\{0\} \times (K_i \sqcup \tau K_i))$ denotes a tubular neighbourhood of $(K_i \sqcup \tau K_i)$,
diffeomorphic to $2(S^1 \times B^3)$.
Furthermore, $N_i \cong \overline{\mathbb{CP}}^{\,2}\setminus\nu(S^1)$.

Consider mixed operations on $R \cong S^1 \times S^2$.
Equip $R \sqcup \tau R$ with a $\tau$-invariant metric of positive scalar curvature and perturb the CSD functional with a Morse function as in \cite[\S5.3]{KMOS}.
The ordinary \textit{HM}-\emph{from} complex $\widehat C_{\bullet}(R)$ has trivial differential.
Moreover, the generators are two towers $\{e^0_j, e^1_j : j < 0\}$,
consisting of boundary-unstable critical points.
In particular, the ordinary $U$-map acts by $Ue^\mu_j=e^\mu_{j-1}$.
By Lemma~\ref{lem:exchanged-pairs}(ii) this gives a description of the real monopole Floer complex of $R\sqcup \tau R$.
Denote by
\[\mu_i = \mu(V_i \setminus (N_i \sqcup \tau N_i))\colon\widehat C_{\bullet}(R)\,\otimes\,\widecheck C_i\longrightarrow \widecheck C_i\]
the mixed operation of the complement of $N_i \sqcup \tau N_i$.
In particular, $N_i \sqcup \tau N_i$ contains no marking set.
Relating the $U$-map to $\upsilon$-map, combined with the calculations of \cite[\S5.3]{KMOS}, yields the following lemma.

\begin{lem}\label{lem:detector}
There is an endomorphism $q_i$ of the diagonal real monopole homology $\widecheck{\hmr}(Y_i)$ which commutes with $\upsilon_i$ and satisfies
	\begin{equation}\label{eq:detector-factor}
		(\psi_i)_*=1+\upsilon_i^2\,q_i .
	\end{equation}
\end{lem}

\begin{proof}
By ~\cite[Proposition~5.5]{KMOS}, counting 0-dimensional family moduli spaces over $P$ defines an operator $G_i\colon \check C_i\to \check C_i$ with
\[\check \partial_iG_i+G_i \check\partial_i=\psi_i+\mu_i(\bar n^u\otimes-),\]
	where $\bar n^u\in\widehat C_\bullet(R \sqcup \tau R)$ is the real relative
	invariant of $N_i \sqcup\tau N_i$ which,
	by Lemmas~\ref{lem:exchanged-pairs}(ii) and~\ref{lem:spinc-composites}(ii),
	is identified with the ordinary relative
	invariants of $(N,\mathfrak s_k)$.
The formula of $\bar n^u$ is identical to that of \cite[Lemmas~5.7 \&~5.10]{KMOS};
the contributions to $\bar n^u$ of $\mathfrak s_k$ and $\mathfrak s_{-1-k}$ add up to $e^1_{-1-k(k+1)/2}$ for $k\ge0$, so mod two
\begin{equation*}
	\bar n^u=\sum_{k=0}^{\infty}e^1_{-1-k(k+1)/2}
	=e^1_{-1} +Uf, \quad f :=\sum^{\infty}_{k=1}e^1_{-k(k+1)/2}.
\end{equation*}

The decomposition~\eqref{eq:Vi-decomp-cylinder} shows that gluing two copies of $(S^1 \times B^3)$ to $V_i \setminus (N_i \sqcup \tau N_i)$ gives the cylinder $I \times Y_i$.
By \cite[Lemma~5.11]{KMOS}, the relative invariant of $S^1\times B^3$ is $e^1_{-1}$.
From now on, we use the same notation $\mu_i$ for the operation on homology.
The real version of \cite[Proposition~5.6]{KMOS}, via Lemma~\ref{lem:exchanged-pairs}(iii), identifies $\mu_i(e^1_{-1}\otimes-)$ with the identity map by Proposition~\ref{prop:diagonal-functoriality}.
It follows that $\mu_{i}(e^1_{-1}\otimes x)=x$.
Moreover, moving the puncture through the connected exterior of $N_i \sqcup \tau N_i$ shows $\mu_{i}(z\otimes\upsilon_ix)=\upsilon_i\,\mu_{i}(z\otimes x)$ (\cite[Proposition~13.3]{Li22}).
Finally, we have:
	\begin{equation*}
		\mu_{i}(Uf \otimes x)=\upsilon_i^2\,\mu_{i}(f\otimes x).
	\end{equation*}
Indeed, $U$ on $\widehat C_\bullet(R)$ is the ordinary point operation $q\in R$ obtained by evaluating the first Chern class of the line bundle $\mathbb{L}_q$ defined by $g \mapsto g(q)$.
Moving $q$ continuously to a marked point $p\in\mathbf p$ gives a continuous family of complex line bundles.
At $p$, the homomorphism $g \mapsto g(p)$ has values in $\{\pm 1\} \subset S^1$, and thus coincides with $\epsilon_{\mathbf p}$.
See \cite[\S~13.4]{Li22}.
This complex bundle is $\mathbb L_{\mathbf p}\otimes_{\mathbb{R}}\mathbb{C}$, so the claim follows from
\[c_1(\mathbb{L}_{\mathbf{p}} \otimes_{\mathbb{R}} \mathbb{C})\bmod2=w_2(\mathbb L_{\mathbf p}\oplus\mathbb L_{\mathbf p})
	=w_1(\mathbb L_{\mathbf p})^2=\upsilon_i^2.\]
Identity \eqref{eq:detector-factor} follows by setting $q_i(x)=\mu_{i}(f\otimes x)$.
\end{proof}

\begin{proof}[Proof of Theorems~\ref{thm:exact-triangle-to}
and~\ref{thm:exact-triangle}]
The map $\upsilon$ on $\widecheck{\HMR}_{\Delta,\bullet}(Y,\tau,\mathbf p)$ is lcoally nilpotent, i.e. for every class $x$, there exists $N$ such that $\upsilon^N x=0$.
Indeed, the proof is identical to that of the local nilpotence of the $U$-map~\cite[Remark~2.4.1]{KMOS}, using the admissibility of perturbation.
By Lemma~\ref{lem:detector}, nilpotence implies, for every $x$ and for all $N$ sufficiently large,
$(\upsilon_i^2q_i)^N x=q_i^N \upsilon_i^{2N}x=0$.
Therefore the following series is finite:
\[
(\psi_i)_*^{-1}(x)
=
\sum_{n\geq0}(\upsilon_i^2q_i)^n x.
\]
We conclude that the \emph{to}-chain map
$\psi_i \colon \check C_i \to \check C_i$ is a quasi-isomorphism.
By Lemma~\ref{lem:h-identity} and the triangle detection lemma,
this proves Theorem~\ref{thm:exact-triangle-to}.
Furthermore,
$\tilde\psi_i$ is lower triangular (w.r.t.\eqref{eq:v-SES}), of the form:
\[
\widetilde\psi_i=
\begin{bmatrix}
\psi_i&0\\
*&\psi_i
\end{bmatrix}.
\]
Thus  $\widetilde\psi_i$ is also a quasi-isomorphism.
Together with the null-homotopy from Lemma~\ref{lem:null-homotopy-cone},
this proves Theorem~\ref{thm:exact-triangle}.
\end{proof}

\subsection{Absolute degrees of triangles}

Let  $L \subset S^3$ be an oriented link, equipped with an ordering of $\pi_0(C)$, and $F$ be a Seifert surface.
Let $\sigma(L)$ and $\nullity(L)$ denote the signature and nullity of the symmetrized Seifert form of $F$, respectively.
In this convention the positive trefoil $T_{2,3}$ has $\sigma=-2$ and the positive Hopf link has $\sigma=-1$.

Following \cite{LM}, we define the normalized mod-two grading $\widehat{\gr}=\gr^{(2)}+\epsilon(L)$ of $\thmr_{\Delta}(\dbcv(L))$ by the normalization term
	\begin{equation}\label{eq:alexander-normalization}
		\epsilon(L)\equiv\frac{-\sigma(L)-\nullity(L)-|L|+1}{2}\bmod 2.
	\end{equation}
In particular, $\epsilon(U)=0$, $\epsilon(H_\pm)=\tfrac12(\pm1-1)$, and $\epsilon(T_{2,p})=\tfrac12(p-1)$ for $p > 1$ odd.
These values match the computations in~\cite[\S5]{LM}.

Let $F_+,F_-$ be the Seifert surfaces for $L_{+},L_{-}$ in the oriented skein triple.
For $Y_{00}$, denote the induced Seifert surface of $L_+$ in $S^3_0(\gamma)$ by $F_{00}$.
In particular, $F_{00}$ lies in the exterior of $\nu(\gamma)$.
The surgeries carry the ordering of the components; in particular, the number of fixed components stays constant.
Let the outgoing cap of $S^3_0(\gamma)=S^1\times S^2$ be $M=D^2\times S^2$,
 which is admissible (\cite[Definition~2.8]{LM}) as $H_1(M)=0$.

\begin{thm}\label{thm:graded-triangle}
Each cobordism map $\widetilde F_i$ is homogeneous with respect to the absolute mod-two grading \eqref{eq:absolute-grading}, of degrees:
\begin{equation}\label{eq:triangle-degrees-intro}
	\begin{aligned}
		\delta_1&=\tfrac12\bigl(-\sigma_++\sigma_-+h_--h_+\bigr),\\
		\delta_2&=\tfrac12\bigl(1-\sigma_-+\sigma_0+h_0-h_-\bigr),\\
		\delta_3&=\tfrac12\bigl(1-\sigma_0+\sigma_++h_+-h_0\bigr).
	\end{aligned}
\end{equation}
Here, $h_+,h_-,h_0$ are the nullities of $L_+,L_-,L_0$, respectively.
And $\sigma_+,\sigma_-,\sigma_0$ refer to their signatures.
Moreover, the degrees $\delta_i$ (as integers) satisfy $\delta_1+\delta_2+\delta_3=1$.
Thus the triangle above is graded.
In particular, after shifting the grading of the second vertex by $\delta_1$ and the third by $\delta_1+\delta_2$,
the edges of the exact triangle have degrees $(0,0,1)$.
\end{thm}

\begin{cor}
	\label{cor:euler-char-equality}
The Euler characteristics $\chi_+,\chi_-,\chi_{00}$  in the mod-two grading~\eqref{eq:absolute-grading} of the three vertices in \eqref{eq:exact-triangle} satisfy
	\begin{equation}\label{eq:graded-euler}
		\chi_+-(-1)^{\delta_1}\chi_-
		+(-1)^{\delta_1+\delta_2}\chi_{00}=0.
	\end{equation}
\end{cor}

\begin{rem}
	The above result should be compared with the unoriented exact triangle~\cite{Li24}.
	The branched double covers of the unoriented skein triples are surgeries along circles that are fixed set-wise by the involution,
	as opposed to being swapped as in the oriented case.
	The index formula of the real Dirac operator shows that the cobordism maps in the unoriented case are not homogeneous modulo two, due to the additional factor of $1/2$.
	On the other hand, in the oriented cases the spin\textsuperscript{c} structures change in pairs over two disjoint submanifolds that are interchanged by the involution.
	The real index computation on this union reduces to the ordinary computation on a single copy, which does satisfy mod-two homogeneity.
\end{rem}

\begin{figure}[htbp]
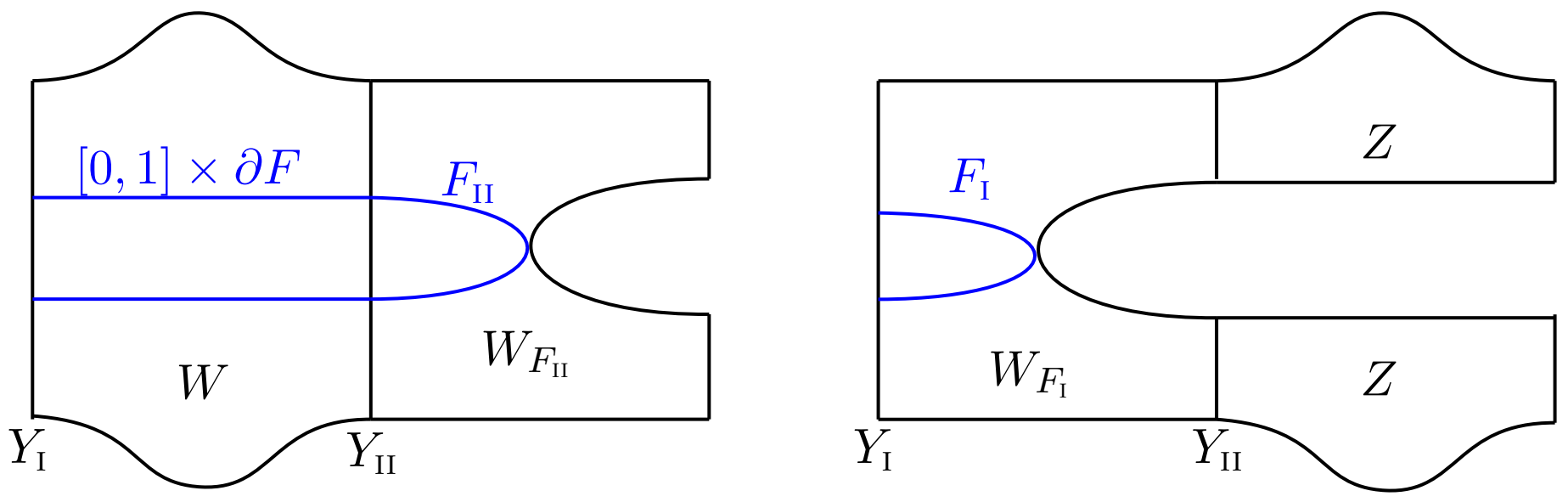

	\centering
	\figcapcommute
	\caption{Two diffeomorphic cobordisms, obtained by isotoping the branch set, which is away from the handle attaching regions.}
	\label{fig:cap-commute}
\end{figure}

The remainder of this section is devoted to the proof of Theorem~\ref{thm:graded-triangle}.
Let $W \colon Y_{\textsc{i}} \to Y_{\sctwo}$ denote one of $W_1$, $W_2$, $W_3$.
Denote  $Z=W/\tau$,
which is the trace of a surgery on $\gamma$ or on its dual circle, downstairs.
Denote the Seifert surfaces by $F_{\scone},F_{\sctwo}$,
and the corresponding Seifert cobordisms by $W_{F_{\scone}},W_{F_{\sctwo}}$.
For a real $4$-manifold $(X,\tau)$, let $\chi_R(X)$ and $\sigma_R(X)$ be $\tau$-skew-invariant Euler characteristic and the signature.
Both are additive under gluing along $\tau$-invariant 3-manifolds, using the same proof as the ordinary additivity, but with $(-\tau^*)$-invariant cohomologies.
To compute the mod two gradings, we relate the Seifert caps of $Y_{\scone}$ and $Y_{\sctwo}$ by the following lemma.
See Figure~\ref{fig:cap-commute} for a schematic picture. 

\begin{lem}\label{lem:cap-commutation}
There is an equivariant diffeomorphism rel boundary for which $\tau$ interchanges the outgoing $Z$-ends:
	\begin{equation}\label{eq:cap-commutation}
		W\cup_{Y_{\sctwo}}W_{F_{\sctwo}}\;\cong\;
		W_{F_{\scone}}\cup_{Y'_{\scone}\sqcup Y'_{\scone}}(Z\sqcup  Z).
	\end{equation}
Moreover, this diffeomorphism preserves the pushed-in Seifert surfaces $F_{\scone},F_{\sctwo}$ and the normal framings $\uptau$ (\cite[\S3.1]{LM}).
Thus for any admissible outgoing cap $M$ of $Y'_{\sctwo}$,
we have the following diffeomorphism:
	\begin{equation*}
		W\cup_{Y_{\sctwo}}W_{F_{\sctwo},M}\;\cong\;W_{F_{\scone},\,Z\cup M}.
	\end{equation*}
\end{lem}

\begin{proof}
The quotient, $W/\tau=Z$ is $I\times Y'_{\scone}$ with a $2$-handle attached along $\gamma$ (or the dual curve of $\gamma$) in the outgoing boundary of $Z$.
Recall that the Seifert surfaces $F_+,F_-,F_{00}$ are images of the same Seifert surface $F \subset S^3 \setminus \nu(\gamma)$ under the three fillings.
Thus $F_{\scone}$ and $F_{\sctwo}$ are identified inside the product region $I \times (Y_{\scone}' \setminus \nu(\gamma))$.
The 2-handle is attached inside $\nu(\gamma)$, and hence away from the Seifert surface $F$.

View both sides of \eqref{eq:cap-commutation} as branched covers of $Z$.
On the left side, the Seifert surface is $I \times \partial F$ inside $W/\tau$, joined by $F$ pushed in from the outgoing boundary of $W/\tau$.
On the right side, the branch surface is $F$ pushed in from the incoming boundary of $W_{F_{\scone}}$.
By isotopies supported in $I \times(Y'_{\scone}\setminus\nu(\gamma))$, the branch surface $F_{\sctwo}$ on the left can be isotoped into the $F_{\scone}$ by shrinking the product part $I \times \partial F$.
See Figure~\ref{fig:cap-commute}.
This isotopy can be chosen to preserve the normal framing $\uptau$ of $F$.
This proves the first statement.

The second statement follows by capping both sides of \eqref{eq:cap-commutation} with $M \sqcup M$ and viewing $W_{F_{\scone}}\cup(Z\sqcup Z)\cup(M\sqcup M)$ as the Seifert cap of $Y_{\scone}$, using outgoing cap $Z \cup M$.
\end{proof}

Next, we show the degree of the cobordism map of a surgery trace is determined by the de Rham part of~\eqref{eq:absolute-grading}, and hence expressible in terms of Betti numbers.
This is due to the fact that all surgeries take place away from the Seifert surfaces and hence do not affect the Arf terms.
Compare following express~\eqref{eq:general-degree} with its real analogue \cite[Definition~25.4.1]{KM}.

\begin{lem}\label{lem:general-degree}
	For all \RRSC\ structure $\widetilde{\sr}_W$ on $W$, the map
	$\widecheck{\HMR}_\Delta(W;\widetilde{\sr}_W)$ has the same mod-two degree:
	\begin{equation}\label{eq:general-degree}
		\frac{\chi_R(W)+\sigma_R(W)+b^1_{R}(Y_{\sctwo})-b^1_{R}(Y_{\scone})}{2}
		\pmod 2.
	\end{equation}
Moreover, this holds also for the operator $\check K_W$ in the definition of the cone map $\thmr_{\Delta}(W)$ as in~\eqref{eq:checkKdefn}.
\end{lem}

\begin{proof}
 Fix $\widetilde{\sr}_W$ with restrictions
	$\widetilde{\sr}_{\scone}$, $\widetilde{\sr}_{\sctwo}$.
Choose an admissible outgoing cap $M$.
Denote $W_{\sctwo} = W_{F_{\sctwo},M}$, and let $\hat{\mathfrak t}_{\sctwo}$ be a real spin$^c$ extension of $\widetilde{\sr}_{\sctwo}$.
By Lemma~\ref{lem:cap-commutation}, the cap
	$W_{\scone}:=W\cup_{Y_{\sctwo}}W_{\sctwo}$ is an admissible Seifert cap of $Y_{\scone}$.
Denote $\hat{\mathfrak t}_{\scone}=\mathfrak{t}_W\cup\hat{\mathfrak t}_{\sctwo}$.
Fix nondegenerate critical points $\mathfrak a$ on $Y_{\scone}$ and $\mathfrak b$ on $Y_{\sctwo}$.
We analyze pair-by-pair the terms in the mod-two grading difference~\eqref{eq:absolute-grading} induced by the caps $W_{\scone},W_{\sctwo}$.
That is,
\begin{equation*}
\gr^{(2)}_\Delta(\mathfrak b)-\gr^{(2)}_\Delta(\mathfrak a)=\gr^\sigma_{SW}(\mathfrak b;W_{\sctwo},\widehat{\mathfrak t}_{\sctwo})- \gr^\sigma_{SW}(\mathfrak a;W_{\scone},\widehat{\mathfrak t}_{\scone})+\iota_R(W_{\sctwo})-\iota_R(W_{\scone})+\Arf(\widehat{\mathfrak s}_{\mathfrak t_{\sctwo}}) -\Arf(\widehat{\mathfrak s}_{\mathfrak t_{\scone}}).
\end{equation*}

Suppose that the coefficient of $\mathfrak b$ in $\widecheck m(W)(\mathfrak a)$ is nonzero.
Every contribution to this coefficient in \eqref{eq:diagonal-cobordism-map} has total Fredholm index zero.
Indeed, $m_o^o$,$m_s^o$,$m_o^u$,$m_s^u$ and $\bar m_s^s$ count index zero solutions.
Moreover,
a trajectory counted by $\bar\partial_{\scone,u}^s$ has index zero,
a solution counted by $\bar m_u^s$ has index $-1$,
and a trajectory counted by $\partial_{\sctwo,o}^u$ or $\partial_{\sctwo,s}^u$ has index $1$.
Thus the terms
$m_o^u\bar\partial_{\scone,u}^s$ and
$m_s^u\bar\partial_{\scone,u}^s$ have total index $0+0$,
and the terms $\partial_{\sctwo,o}^u\bar m_u^s$ and
$\partial_{\sctwo,s}^u\bar m_u^s$ have total index $1-1$.
By additivity of APS index,
we have the cancellation of $\gr^{\sigma}_{SW}$ terms:
\[
\gr^\sigma_{SW}
(\mathfrak a;W_{\scone},\widehat{\mathfrak t}_{\scone})
-
\gr^\sigma_{SW}
(\mathfrak b;W_{\sctwo},\widehat{\mathfrak t}_{\sctwo})
\equiv 0 \bmod 2.
\]

Next, we claim $\Arf(\widehat{\mathfrak s}_{\mathfrak t_{\scone}})=\Arf(\widehat{\mathfrak s}_{\mathfrak t_{\sctwo}})$.
The construction~\cite[\S3.1--3.3]{LM} of the spin structure $\mathfrak s_{\mathfrak t}$ on $F$ is local.
Indeed, the spin bundle is defined as the eigenbundle of the complex structure $\rho(n_1)\rho(n_2)$ along the pushed-in Seifert surface, using any local frame $(n_1,n_2)$ of the normal bundle.
The spin marking data on $\mathfrak s_{\mathfrak t}$ is given by the \rrsc\ structures on the boundary 3-manifold.
Lastly, the marking on the planar capping surface $G$ is determined by the chosen ordering of $\pi_0(C)$  \cite[\S3.6]{LM}.
All these data are preserved by the diffeomorphism~\eqref{eq:cap-commutation},
which is the identity on $F$ and the marked arcs,
so $\mathfrak{s}_{\mathfrak t_{\scone}} \cong \mathfrak s_{\mathfrak t_{\sctwo}}$.
Since the marking set on $W$ consists of $I \times \mathbf{p}$, lying in the complement of the surgery region,
the Clifford trivialization $T_{I \times \mathbf{p}}$ transports the markings.
It follows that the $\mathfrak{s}_{\mathfrak t_{\scone}} \cong \mathfrak s_{\mathfrak t_{\sctwo}}$ as marked spin structures.
Hence $\widehat{\mathfrak s}_{\mathfrak t_{\scone}} \cong \widehat{\mathfrak s}_{\mathfrak t_{\sctwo}}$.

Lastly, if $X$ is any connected cap with boundary $Y$ also connected, the abelian ASD index $\iota_R(X)$ with APS boundary condition can be expressed as
$(\chi_R(X)+\sigma_R(X)-b^1_R(Y))/2$.
Indeed, this follows from the proof of its ordinary counterpart~\cite[Lemma~25.4.2]{KM}, using the $-\tau^*$-invariant version of Poincar\'{e}--Lefschetz duality.
The additivity of $\chi_R$ and $\sigma_R$ implies
\[\iota_R(W_{\sctwo})-\iota_R(W_{\scone})
=-\frac{1}{2}\left(\chi_R(W)+\sigma_R(W)+b_R^1(Y_{\sctwo})-b_R^1(Y_{\scone})\right).\]

Since the APS indices and Arf terms both cancel out,
we conclude that a nonzero entry of $\widecheck m(W)$ from $\mathfrak a$ to $\mathfrak b$ gives rise to the grading shift~\eqref{eq:general-degree}.
Finally, the family moduli space defining $\check K_W$ has one more dimension than the moduli spaces defining $\widecheck{m}(W)$,
but with dimension cut down by one by inserting $\mathfrak{c}_{-2}$ at the puncture.
The same argument as above can be applied to family moduli spaces as in \cite[\S2.2]{Li24}.
We conclude that $\check K_W$ and $\widecheck m(W)$, and hence $\widetilde{m}(W)$ have the same degree with respect to the absolute mod-two grading.
\end{proof}

By Lemma~\ref{lem:general-degree}, to finish the proof of Theorem~\ref{thm:graded-triangle} it suffices to compute the grading shifts~\eqref{eq:general-degree} of the cobordisms for $W_1,W_2,W_3$.

\begin{proof}[Proof of Theorem~\ref{thm:graded-triangle}]
Since $W_i$ is obtained by attaching a pair of $2$-handles,
$H_2(W_i,Y_i;\mathbb{R}) = \mathbb{R}^2$.
Since $\tau_*$ anti-fixes a one-dimensional subspace of this $\mathbb{R}^2$,
we have $\chi_R(W_i)=1$.
Take $\sigma_R$ of both sides of \eqref{eq:cap-commutation}:
\begin{equation*}
	\sigma_R(W) + \sigma_R(W_{F_{\sctwo}})=
	\sigma_R(W_{F_{\scone}}) + \sigma_R(Z\sqcup Z) = \sigma_R(W_{F_{\scone}}) + \sigma(Z).
\end{equation*}
Here, we used Novikov additivity and $\sigma_R(Z\sqcup Z)=\sigma(Z)$.
Moreover, since the quotient of a Seifert cobordism $W_F$ is a cylinder, we have $\sigma_R(W_F)=\sigma(W_F)-\sigma(W_F/\tau)=\sigma(W_F)$.
Specializing to the Seifert surfaces $F_+,F_-,F_{00}$, we have
	\[
	\sigma_R(W_i)=-\sigma_{i}+\sigma_{i+1}+\sigma(Z_i).
	\]
Here, $\sigma_i$ denotes the signature of the corresponding link $L_+$, $L_-$, or $L_0$, with label $i$ as in the vertices of the triangle.
By our orientation convention of the Seifert cobordisms and \cite[Theorem~3]{GordonLitherland}, $\sigma(W_{F_i}) = -\sigma_i$.
Since the downstairs traces $Z_1,Z_2,Z_3$ relate the $(\infty,-1,0)$-surgeries of $S^3$, they have intersection forms $(-1)$, $(0)$, $(0)$, respectively.
Thus $\sigma(Z_1)=-1$ and $\sigma(Z_2)=\sigma(Z_3)=0$.

Recall $Y_{00}$ is the outgoing boundary of $I \times Y_0$ with a 4-dimensional equivariant 1-handle attached.
The quotient of this 1-handle is again a 1-handle,
so resulting class in $H^1(Y_{00};\mathbb{R})$ is $\tau$-invariant\footnote{This class is dual to the non-separating sphere $\Sigma$ as in \S\ref{ssec:top-or-triple}, whose orientation is preserved by $\tau$.}.
In particular, $b^1_R(Y_{00}) = \nullity(L_0)$.
To compare the signatures of $W_{F_{00}}$ and $W_{F_0}$, note that the Seifert cap $W_{F_{00},S^1\times B^3}$ is obtained from $W_{F_0,B^4}\sqcup B^4$ by attaching two equivariant 4-dimensional 1-handles,
as in Figure~\ref{fig:Fgamma}.
Attaching 1-handles preserves signature, and the outgoing caps $S^1 \times B^3$ and $B^4$ have zero signature.
Therefore, $\sigma (W_{F_{00}})=-\sigma(L_0)$.
Plugging the calculations above for $(Y_+,Y_-),(Y_-,Y_{00}),(Y_{00},Y_+)$ into \eqref{eq:general-degree} yields \eqref{eq:triangle-degrees-intro}.

Finally, the operators $\check H_i$ and $I_i$ are defined using families of one dimension higher than that of $\check F_{i+1}\check F_{i}$.
Let $\deg(\cdot)$ denote the mod-two grading shift.
Then
\[\deg(\check H_i)=\deg(I_i)=\delta_i+\delta_{i+1}+1, \quad
\deg(\psi_i)=\deg(\widetilde \psi_i)=0.\]
This shows the null-homotopies and isomorphisms in the triangle detection lemma are homogeneous.
Hence the triangle is exact also with the mod-two grading.
\end{proof}

\section{Connected sum theorem}
\label{sec:connected-sum}
This section proves an equivariant connected sum formula in a more general setting.
Unlike the preceding discussions, the Floer homology groups in this section use the framework of \cite{Li26}.
In particular, all configurations are downstairs,
and the real gauge groups are framed rather than diagonal.
\subsection{Double connected sum}
\label{ssec:double-connected-sum}
Let $(Y, \tau, \mathbf{p})$ be a closed marked real 3-manifold,
and $(Z, \tau_Z, \mathbf{q})$ be another real 3-manifold,
where the basepoint set $\mathbf{q}$ is allowed to empty.
Choose $x_1, x_2 \in \fix(\tau)$ and $z_1, z_2 \in \fix(\tau_Z)$, away from the basepoint sets.
Let $B_{x1}^3,B_{x2}^3,B_{z1}^3,B_{z2}^3$ be invariant 3-balls around $x_1,x_2,z_1,z_2$, respectively.
Fix equivariant identifications $\partial B^3_{xi} \cong \partial B^3_{zi}$,
where the involution acts by $180^{\circ}$-rotation around an axis.
Define
\begin{equation*}
	Y \#_2 Z := (Y \setminus (B^3_{x1} \sqcup B^3_{x2})) \, \cup_{\partial B^3_{xi} = S^2 \times \{-1\}}
	2(S^2 \times [-1,+1])\, \cup_{\partial B^3_{zi} = S^2 \times \{+1\}}
	\, (Z \setminus (B^3_{z1} \sqcup B^3_{z2})).
\end{equation*}
Denote $\Sigma^2_i = S^2 \times \{0\}$.

\begin{figure}[htbp]
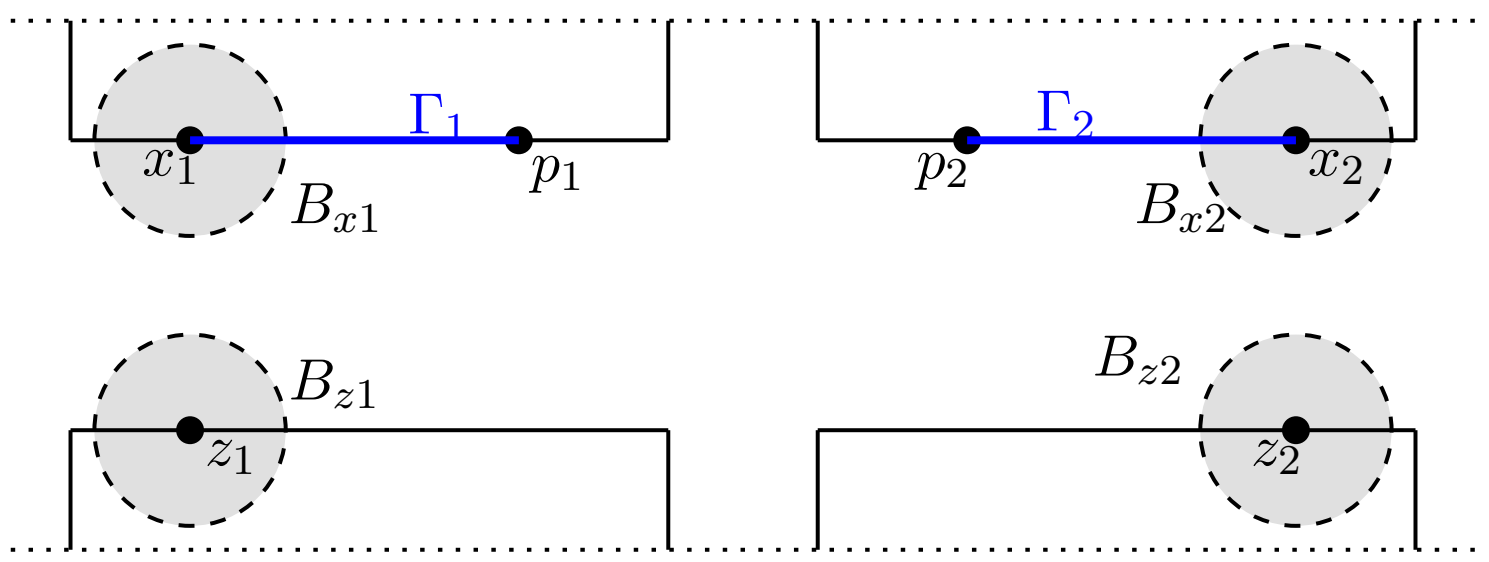

	\centering
	\figfourballs
	\caption{Local balls around the points $x_1,x_2,z_1,z_2$, shaded grey. The arcs $\Gamma_1,\Gamma_2$ from basepoints $p_i$ to $x_i$ are coloured in blue. The upper half of the figure lives in $Y$ and the lower half lives in $Z$.}
	\label{fig:fourballs}
\end{figure}

\begin{figure}[htbp]
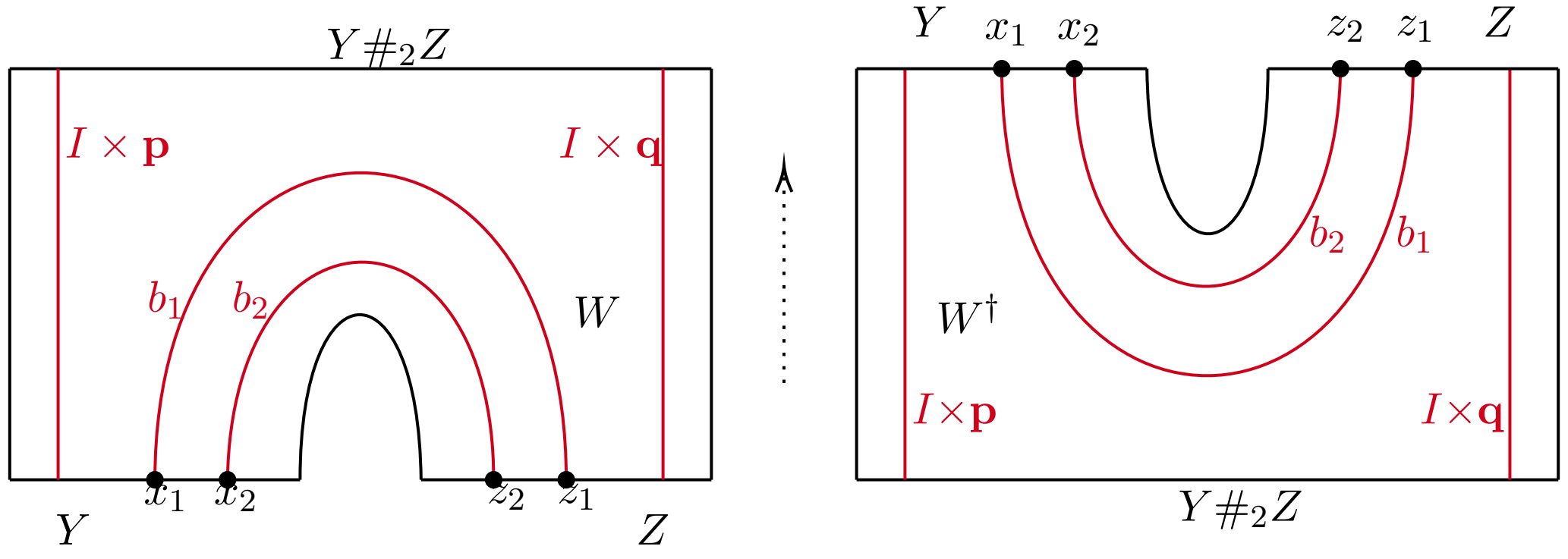

	\centering
	\figWWdag
	\caption{The 1-handle cobordism $W \colon (Y \sqcup Z) \to Y \#_2 Z$ (left) and its reversal $\overline{W}$ (right).
	Marking sets are coloured red.}
	\label{fig:WWdag}
\end{figure}

Let $W$ be the real cobordism obtained by attaching a pair of 1-handles to $I \times (Y \sqcup Z)$.
The marking set $\mathbf{a}$ is the union of $I \times (\mathbf{p} \sqcup  \mathbf{q})$
with arcs $b_i$ that go through the cores of the 1-handles and join $x_i$ to $z_i$:
\[
	(W,\mathbf{a}) \colon
(Y,\mathbf{p}\cup\{x_1,x_2\}) \sqcup (Z, \mathbf{q} \cup \{z_1,z_2\})
\longrightarrow
(Y \#_2 Z,\mathbf{p} \sqcup \mathbf{q}).
\]
Let $(\overline{W},\bar{\mathbf{a}})$ be the reverse cobordism.
See Figure~\ref{fig:WWdag}.

For convenience, denote $\mathbf{p}^{\sharp} = \mathbf{p} \cup \{x_1, x_2\}$ and $\mathbf{q}^{\sharp} = \mathbf{q} \cup \{z_1, z_2\}$.
Assume $x_1,x_2$ lie in components of $\fix(\tau)$ of $Y$ that contain some points $p_1,p_2$ in $\mathbf{p}$.
It follows that the forgetful map
from $\RRSC(Y,\mathbf{p}^{\sharp})$ to $\RRSC(Y,\mathbf{p})$ is 4-to-1.
Let $\Gamma = (\Gamma_1,\Gamma_2)$ be a pair of arcs in $\fix(\tau)$ going from any pair of basepoints $p_1,p_2$ in $\mathbf{p}$ to $x_1$ and $x_2$, respectively.
See Figure~\ref{fig:fourballs}.
Then the Clifford trivializations at the chosen basepoints determine, via trivialization along $\Gamma$, Clifford trivializations $T_{x_1},T_{x_2}$.

In particular, every fibre of $\RRSC(Y,\mathbf{p}^{\sharp}) \to \RRSC(Y,\mathbf{p})$ of $(\mathfrak{s},T, \sr)$ has Clifford trivializations $(T,\epsilon_1 T_{x_1}, \epsilon_2 T_{x_2})$ for $(\epsilon_1, \epsilon_2) \in (\pm 1)^2$.
The corresponding configuration spaces for the $(\epsilon_1,\epsilon_2)$-fibres are identical.
Including the $\RRSC$ structure $(\mathfrak{s}, T, \sr)$ into $(\mathfrak{s}, T \cup \{ T_{x_1},T_{x_2}\}, \sr)$ defines the following inclusion of chain complexes:
\begin{equation*}
	i_{\Gamma} \colon \widetilde{\cmr}(Y, \mathbf{p}) \to  \widetilde{\cmr}(Y, \mathbf{p}^{\sharp}).
\end{equation*}
Further, consider the projection onto $(\mathfrak{s}, T, \sr)$ that sends every $\RRSC$-summand, except for $(\mathfrak{s}, T \cup \{T_{x_1}, T_{x_2}\}, \sr)$, to zero:
\begin{equation*}
	p_{\Gamma} \colon \widetilde{\cmr}(Y, \mathbf{p}^{\sharp})
	\to \widetilde{\cmr}(Y, \mathbf{p}).
\end{equation*}

Let
\begin{align*}
I_{\Gamma} = i_{\Gamma} &\otimes 1 \colon  \widetilde{\cmr}(Y, \mathbf{p}) \otimes \widetilde{\cmr}(Z, \mathbf{q}^{\sharp}) \
\longrightarrow \widetilde{\cmr}(Y, \mathbf{p}^{\sharp}) \otimes \widetilde{\cmr}(Z, \mathbf{q}^{\sharp}),\\
\Pi_{\Gamma} = p_{\Gamma} &\otimes 1 \colon \widetilde{\cmr}(Y, \mathbf{p}^{\sharp}) \otimes \widetilde{\cmr}(Z, \mathbf{q}^{\sharp})
\longrightarrow
\widetilde{\cmr}(Y, \mathbf{p}) \otimes \widetilde{\cmr}(Z, \mathbf{q}^{\sharp})
\end{align*}
Both are chain maps.
Define
\begin{align*}
	\mathsf{F}_{\Gamma} = \Pi_\Gamma \circ \widetilde{\cmr}(\overline{W},\overline{\mathbf{a}})
&\colon \widetilde{\cmr}(Y \#_2 Z, \mathbf{p} \sqcup \mathbf{q})\longrightarrow \widetilde{\cmr}(Y, \mathbf{p}) \otimes \widetilde{\cmr}(Z, \mathbf{q}^{\sharp}),\\
	\mathsf{G}_{\Gamma} = \widetilde{\cmr}(W,\mathbf{a}) \circ I_\Gamma &\colon
\widetilde{\cmr}(Y, \mathbf{p}) \otimes \widetilde{\cmr}(Z, \mathbf{q}^{\sharp}) \longrightarrow \widetilde{\cmr}(Y \#_2 Z,\mathbf{p} \sqcup \mathbf{q}).
\end{align*}

\begin{thm}
	\label{thm:connected-sum}
	We have chain homotopy equivalences
	$
		\mathsf{F}_{\Gamma} \circ \mathsf{G}_{\Gamma} \simeq \id$ and $
		\mathsf{G}_{\Gamma} \circ \mathsf{F}_{\Gamma} \simeq \id,
	$
	and thus the isomorphism
	\begin{equation*}
		\widetilde{\HMR}(Y \#_2 Z, \mathbf{p} \cup \mathbf{q}) \cong
		\widetilde{\HMR}(Y, \mathbf{p}) \otimes \widetilde{\HMR}(Z, \mathbf{q} \cup \{z_1,z_2\}).
	\end{equation*}
\end{thm}

Before commencing the proof, we analyze the marking data and relative real \text{spin\textsuperscript{c}} structures on the cobordisms.
\begin{lem}
Restriction of the set of \RRSC\ structures on $W$ defines a 4-to-1 map and an isomorphism as follows:
\begin{equation}
\label{eq:rrsc-W-Y-Z-YZ}
	\RRSC(Y \#_2 Z, \mathbf{p} \sqcup \mathbf{q}) \xleftarrow{1:4}
	\RRSC(W, \mathbf{a}) \xrightarrow{\cong}
	\RRSC(Y, \mathbf{p}^{\sharp}) \times
	\RRSC(Z, \mathbf{q}^{\sharp}).
\end{equation}
Furthermore, the four fibres over an \RRSC\ structure of $\RRSC(W, \mathbf{a}) \to \RRSC(Y\#_2 Z, \mathbf{p} \cup \mathbf{q})$ differ by the signs of the markings on the arcs $b_i$ from $x_i$ to $z_i$.
Finally, exactly one of these four fibres, when restricted to $\rrsc(Y,\mathbf{p}^{\sharp})$ on the right, lies in the image of $i_{\Gamma}$.
\end{lem}

\begin{proof}
From $Y\#_2 Z$ to $Y \sqcup Z$, the cobordism $W$ is obtained by attaching two $3$-handles to $I \times Y\#_2 Z$.
Since $\tau$ preserves the orientation of $\Sigma_i$
and $\tau^*c_1(\mathfrak s)=-c_1(\mathfrak s)$,
we have $\langle c_1(\mathfrak s),[\Sigma_i]\rangle=0$.
This implies the ordinary spin\textsuperscript{c} structure on $Y \#_2 Z$ extends over the 3-handle.
Since $H^2(W,Y\#_2 Z;\mathbb{Z})=0$, the extension is unique.
Furthermore, the real structure on the attaching sphere also extends over the standard invariant ball.
Since the basepoints $\mathbf{p} \cup \mathbf{q}$ are away from the connected sum region, attaching $3$-handles preserves $H^1(W,Y\#_2 Z \cup (I \times (\mathbf{p} \cup \mathbf{q}));\mathbb{Z})$.
Thus the classification of \RRSC\ structures~\cite[\S3]{Li26} shows the relative real extension is unique.

Next, consider the marking arcs $b_1,b_2$.
Each $b_i$ admits two compatible Clifford trivializations $\pm T_{b_i}$, differing by a sign.
To show that these $4=2 \times 2$ choices give distinct \RRSC\ structures,
suppose $g$ is a real gauge transformation such that $u|_{\mathbf{p}\cup \mathbf{q}} = 1$.
Since $x_i$ is connected to $p_i$ through $\Gamma_i$ in the fixed set, so $g(x_i) = 1$.
Furthermore, $z_i$ is connected to $x_i$ through $b_i$, so $g|_{z_i} = g|_{b_i} =1$.
Thus no two choices are equivalent via $g$.
This also proves the last statement, as the marking changes on $(b_1,b_2)$ induce marking changes on $(x_1,x_2)$.

As for the second arrow of~\eqref{eq:rrsc-W-Y-Z-YZ}, regard $W$ as the result of adding two $1$-handles to $I \times (Y\sqcup Z)$.
Collapse each 1-handle to its core arc $b_i$.
The first and second cohomology rel $(\mathbf p\sqcup\mathbf q \cup b_1 \cup b_2)$ of the resulting space agree with those of
$(Y\sqcup Z,\mathbf p\sqcup\mathbf q\sqcup\{x_1,x_2,z_1,z_2\})$.
To define an \RRSC\ structure on $(W,\mathbf{a})$,
we use the markings at $x_i$ and $z_i$ to identify the spinor bundles over the two attaching balls with the standard local models.
Then the marking data along the $b_i$'s allow us to extend the real spin$^{c}$ structure over to the $1$-handle.
The classification of \RRSC\ structures~\cite[\S3]{Li26} proves the bijection on the right hand side.
\end{proof}

\begin{figure}[htbp]
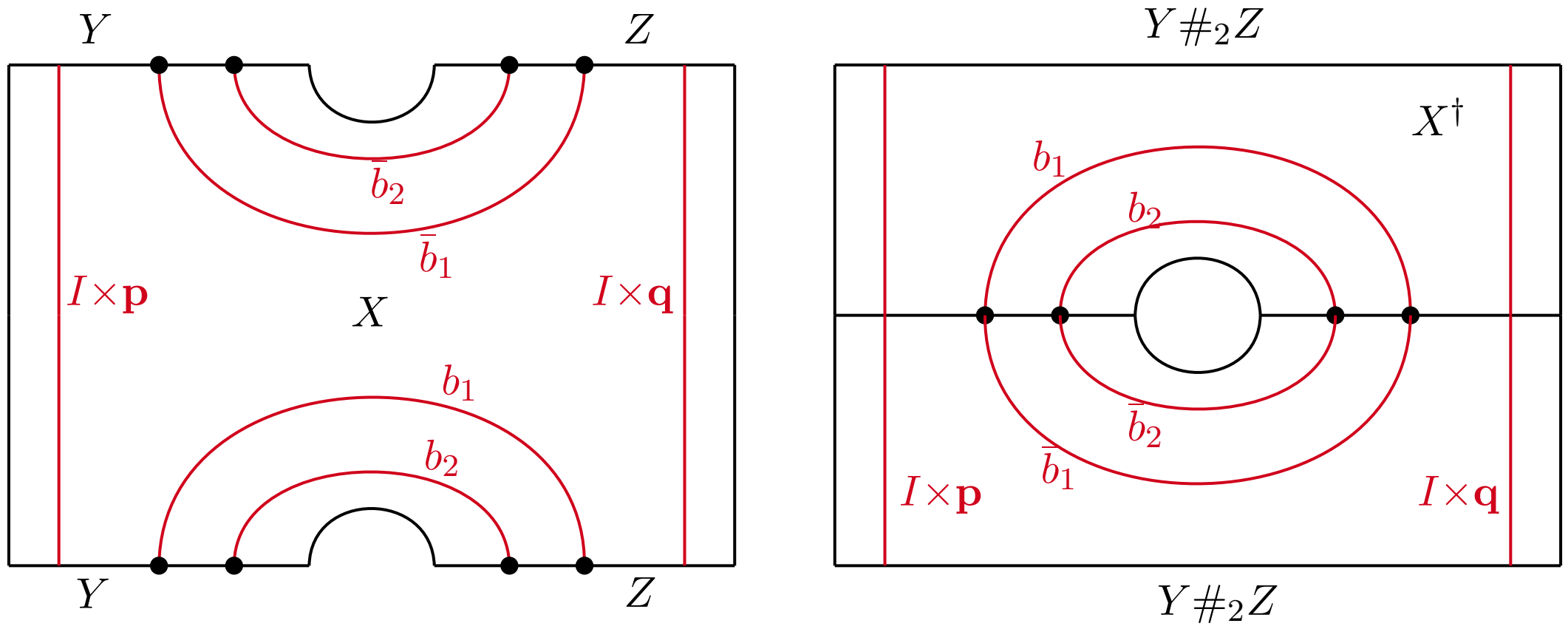

	\centering
	\figcomposites
	\caption{The composite cobordisms $X$ and $X^\dag$, together with their marking set in red.}
	\label{fig:composites}
\end{figure}

Define
\[X := W \cup_{Y \#_2 Z} \overline{W}, \quad X^{\dag} := \overline{W} \cup_{Y \sqcup Z} W.\]
See Figure~\ref{fig:composites}.
There is no ambiguity from gluing $\RRSC$ structures along $Y \#_2 Z$ or $Y \sqcup Z$,
thanks to the following lemma.
\begin{lem}
The restrictions of $\RRSC$ structures to the cobordisms $W$ and $\overline{W}$ induce fibre product descriptions:
\begin{align}
	\RRSC(X,\mathbf{a} \cup
	\overline{\mathbf{a}})&\xrightarrow{\ \cong\ }
	\RRSC(W,\mathbf{a}) \mathop{\times}_{\RRSC(Y \#_2 Z, \mathbf{p}
	\cup \mathbf{q})}\RRSC(\overline{W},\overline{\mathbf{a}}),\\
	\RRSC(X^{\dag}, \overline{\mathbf{a}}\cup \mathbf{a})
	&\xrightarrow{\ \cong\ } \RRSC(\overline{W},\overline{\mathbf{a}})
	\mathop{\times}_{\RRSC(Y \sqcup Z, \mathbf{p}
		\cup \mathbf{q} \cup \{x_1,x_2,z_1, z_2\})}\RRSC(W,\mathbf{a}).
\end{align}
\end{lem}
\begin{proof}
It suffices to show that every framed gauge transformation $g$ on the  3-manifolds $Y \#_2 Z$ and $Y \sqcup Z$, in the middle of $X$ and $X^{\dag}$, extends over $\overline{W}$ and $W$, respectively.
At the middle $Y\sqcup Z$ of $X^{\dag}$, the cobordism $W$ is obtained by attaching two $1$-handles to $I \times (Y \sqcup Z)$.
On each attaching 3-ball, $g=e^{ih}$ for some $\tau^*h=-h$.
For each such 1-handle, extend the function $h$ to a function $f$, and let $\widetilde h=(f-\tau^*f)/2$.
Since $\widetilde h=0$ at every fixed point, $\exp(i\widetilde h) = 1$
along the marking arc $b_i$.
Next, the cobordism $\overline{W}$ attaches two $3$-handles to $I \times (Y\#_2 Z)$.
Each attaching region $S^2\times[-1,1]$ is simply connected.
Choose a fixed point $r_0$ in this region,
and set $g = g(r_0) e^{ih}$ with $h(r_0)=0$.
The same argument gives $h\circ\tau=-h$,
and the same extension and averaging argument extends $g$ over the 3-handles.
\end{proof}

\subsection{Cutting the composite cobordism}
The proof of Theorem~\ref{thm:connected-sum} involves four elementary 4-manifolds:
\[\mathsf{B} = B^4 \sqcup B^4,\quad \mathsf{T} = S^3 \times I,
\quad \mathsf{A} = D^2 \times S^2, \quad \mathsf{S}=S^1 \times B^3.\]

Suppose $I = [-1,+1]$.
In the cylinder $I \times (Y \sqcup Z)$, choose small invariant 4-balls $B^4_{xi},B^4_{zi}$ centred at $x_i,z_i$, respectively.
The involution on the cylinder acts on $B^4_{xi},B^4_{zi}$ and fixes 2-discs.
The marking sets are the interior arcs $(-\epsilon, +\epsilon) \times \{x_i\}$ and $(-\epsilon, +\epsilon) \times \{ z_i\}$.
The involutions on the boundary 3-spheres fix unknots, each marked with two points $x_i^{\pm}$ or $z_i^{\pm}$.
Denote $\mathsf{B}_i = B^4_{xi} \sqcup B^4_{zi}$.
Removing the two pairs of 4-balls yields a cobordism:
\[V := (I \times (Y \sqcup Z)) \setminus (\mathsf{B}_1 \sqcup \mathsf{B}_2).\]
Then
\begin{equation*}
	V \cup (\mathsf{T}_1 \sqcup \mathsf{T}_2) = X, \quad
V \cup (\mathsf{B}_1 \sqcup \mathsf{B}_2) = I \times (Y \sqcup Z),
\end{equation*}
where $\partial \mathsf{T}_i \cong \partial (B^4_{xi} \sqcup B^4_{zi})$.
In particular, the involution on the tube $S^3 \times I$ acts by $180^{\circ}$-rotation around an axis on $S^3$ that fixes the unknot, and acts trivially on $I$.
The marking set on each $\mathsf{T}_i$ consists of two arcs,
which join $x_i^{\pm}$ to $z_i^{\pm}$.

\begin{figure}[htbp]
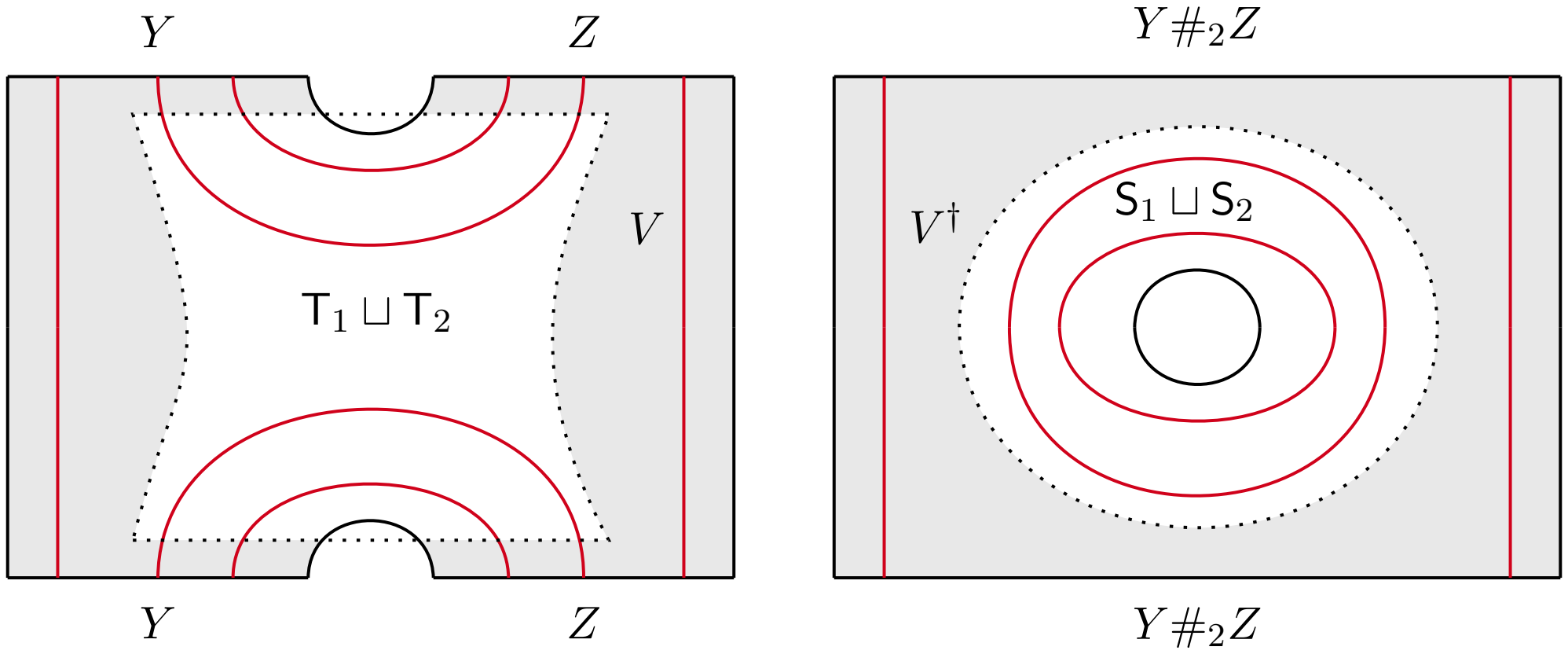

	\centering
	\figVVdag
	\caption{The cobordisms $V,V^{\dag}$, shaded grey inside $X,X^{\dag}$.}
	\label{fig:VVdag}
\end{figure}

Similarly, consider $I \times (Y \#_2 Z)$ and choose neighbourhoods $\mathsf{A}_i \cong D^2 \times \Sigma^2_i$ of $\{0\} \times \Sigma^2_i$.
Removing this pair of $\mathsf{A}$'s defines
\begin{equation*}
	V^{\dag} := I \times (Y \#_2 Z) \setminus (\mathsf{A}_1 \sqcup \mathsf{A}_2).
\end{equation*}
Then \[V^{\dag} \cup (\mathsf{S}_1\sqcup \mathsf{S}_2) = X^{\dag}, \quad
V^{\dag} \cup (\mathsf{A}_1 \sqcup \mathsf{A}_2) = I \times (Y \#_2 Z).\]
Here $\mathsf{S}_1\cup\mathsf{S}_2$ denotes a pair of $S^1 \times B^3$.
By a slight abuse of notation, let $\rot$ be the $180^{\circ}$-rotation around an axis in $B^3$ and $S^2$.
The involution $\tau$ acts on $\mathsf{A}$ by $\id \times \rot$, fixing two 2-discs,
and $\tau$ acts by $\id \times \rot$ on $\mathsf{S}$, fixing one annulus.
The marking sets are empty for $\mathsf{A}$ and $S^1 \times \{0\}$ for $\mathsf{S}$, respectively.
In particular, $S^1 \times \{0\}$ is precisely the union $\bar{b}_i \cup b_i$  of the marking arcs obtained from the cores of the 1-handles and its reverse.
In other words, the marking set $\bar{\mathbf{a}} \cup \mathbf{a}$ is a self-cobordism of $\mathbf{p} \cup \mathbf{q}$,
comprising $I \times (\mathbf{p} \cup \mathbf{q})$ and two circles $\bar{b}_1 \cup b_1$ and $\bar{b}_2 \cup b_2$ in the interior of $X^{\dag}$.
See Figure~\ref{fig:VVdag}.

Assume the $\tau$-invariant Riemannian metrics on the boundaries of $\mathsf{A}_i,\mathsf{B}_i,\mathsf{S}_i,\mathsf{T}_i$ have positive scalar curvature.
In particular, each boundary component of $\partial \mathsf{B}$ and $\partial \mathsf{T}$ is diffeomorphic to $S^3$ where the involution fixes the unknot.
For each \RRSC\ structure on $S^3$ with two basepoints, there is a unique unperturbed framed critical point, which is reducible and nondegenerate by the PSC assumption.
The other two boundaries
\[(\partial \mathsf{A}_i, \tau) = (\partial \mathsf{S}_i, \tau) \cong (S^1 \times S^2, \id \times \rot).\]
The real Picard torus of $S^1 \times S^2$ is trivial for the involution $\id \times \rot$,
as the induced action on cohomology is trivial.
By the same PSC argument as before,
there is a unique unperturbed nondegenerate reducible critical point, for every absolute or relative real spin\textsuperscript{c} structure on $S^1 \times S^2$.
Thus there are no nonconstant trajectories on $\mathbb{R} \times \partial \mathsf{C}_i$, for any $\mathsf{C} \in \{\mathsf{A},\mathsf{B},\mathsf{S},\mathsf{T}\}$.

\begin{rem}
	In contrast to the analysis on $R_1 \cong S^1 \times S^2$ in~\cite[\S5.3]{KMOS}, this real $S^1 \times S^2$ has a particularly simple Floer chain complex, as there is no $S^1$ family of flat connections, modulo gauge.
	Compare with~\cite[\S14.3--4]{Li22}.
\end{rem}

Now, on $\mathsf{A}_i,\mathsf{B}_i,\mathsf{S}_i,\mathsf{T}_i$ we choose $\tau$-invariant metrics of positive scalar curvature that are product near their boundaries.
Attach cylindrical ends to $\mathsf{A}_i,\mathsf{B}_i,\mathsf{S}_i,\mathsf{T}_i$  and to $V, V^{\dag}$.
Since $\tau$ acts as identity on the cohomology on every $\mathbf{C}_i \in \{\mathsf{A},\mathsf{B},\mathsf{S},\mathsf{T}\}$,
only the torsion spin\textsuperscript{c} structures admit real structures.

\begin{lemma}
\label{lem:solution-on-C}
Fix $\mathsf{C} \in \{\mathsf{A},\mathsf{B},\mathsf{S},\mathsf{T}\}$, a (relative) real spin\textsuperscript{c} structure on $\mathsf{C}$,
and a PSC cylindrical-end metric as above.
Then the unperturbed Seiberg--Witten equation on $\mathsf{C}$ admits a unique finite-energy solution modulo the real gauge group $\mathcal{G}(\mathsf{C})$.
Moreover, this solution is regular and of index zero.
\end{lemma}

\begin{proof}
The standard positive scalar curvature argument on $\mathsf{C}$ implies that a finite energy solution is reducible and satisfies $F_A^+=0$.
By Hodge theory~\cite[Proposition~4.9]{APS} on $\mathsf{C}$ with cylindrical end, the curvature $F_A$ lies in
\[
\operatorname{im}\!\left(
H^2(\mathsf C,\partial\mathsf C;\mathbb R)
\longrightarrow H^2(\mathsf C;\mathbb R)
\right)^{-\tau^*},
\]
which vanishes as $\tau$ acts trivially on cohomology.
Thus $F_A=0$. But $b^1_R(\mathsf{C})=0$, so up to real gauge transformations there is only one such reducible solution.
For the same reason, $(d^+ + d^*)_R$ is invertible.
Moreover, PSC implies that $D^+_{A,R}$ is also invertible,
so the linearization of the Seiberg--Witten operator is invertible.
\end{proof}

By Lemma~\ref{lem:solution-on-C}, what remains is the problem of matching \RRSC\ structures on the different pieces.

\begin{lem}\label{lem:match-rrsc-count}
Fix \RRSC\ structures $\widetilde \sr_V$ on $V$, and $\widetilde \sr_{V^{\dag}}$ on $V^{\dag}$.
\begin{enumerate}[label=(\roman*)]
\item
On $V$, suppose the markings at $x_1,x_2$ on both $Y$ ends are specified by $\Gamma$.
Then $\widetilde\sr_V$ extends over $\mathsf T_i$ if and only if it extends over $\mathsf B_i$.

\item
On $V^{\dag}$, the \RRSC\ structure $\widetilde \sr_{V^\dag}$ extends over $\mathsf A_i$
if and only if it admits an extension over $\mathsf S_i$
whose induced marking at $x_i$ on the middle $Y$ is the one
specified by $\Gamma$.
\end{enumerate}
Any such extension of $\widetilde \sr_V$ or $\widetilde \sr_{V^{\dag}}$, if it exists, is unique.
\end{lem}

\begin{proof}
For (i), define $\delta_x=+1$ if the two markings on the two points of $\partial B_i$ extend over the corresponding marked 4-balls in $\mathsf B_i$, and let $\delta_x=-1$ otherwise.
Define $\delta_z \in \{\pm 1\}$ similarly.
Assume the marking at $x_i$ is specified by $\Gamma_i$.
Then $\Gamma_i\times I$ lies in the fixed set once $\mathsf{B}_i$ is filled into $V$, so $\delta_x=1$.
A marking on $\partial \mathsf T_i$ extends over $\mathsf T_i$ if and only if $\delta_x=\delta_z$.
On the other hand, a marking on $\partial \mathsf B_i$ extends into $\mathsf B_i$ if and only if $\delta_x=\delta_z=1$.
Since the markings on the $Y$-ends are given by transporting along $\Gamma$,
the extension conditions on $\mathsf T_i$ and $\mathsf B_i$ are both equivalent to $\delta_z=1$. This proves (i).

For (ii), there is no marking set on $S^1 \times S^2$.
Let $\sr_0$ be the absolute real spin\textsuperscript{c} structure on $\partial\mathsf A_i \cong S^1 \times S^2$ that extends over $\mathsf A_i$.
Then $\sr_0$ extends also to $\mathsf{S}_i$.
The marking set of $\mathsf S_i$ is its core $S^1 \times \{0\}$,
and the two markings along the core differ by $\pm 1$.
Exactly one such marking agrees with the marking specified by $\Gamma$.
Thus the extension of $\widetilde{\sr}_{V^{\dag}}$ from the exterior to $\mathsf S_i$ exists if and only if it extends over $\mathsf A_i$.
This proves (ii).
Finally, since the real gauge groups on the $S^3$ and $S^2 \times S^1$ are homotopic equivalent to $\pm 1$, the extensions are unique.
\end{proof}

\begin{proof}[Proof of Theorem~\ref{thm:connected-sum}]
First, consider $X=V \cup (\mathsf{T}_1 \sqcup \mathsf{T}_2)$ equipped with a one-parameter family $Q \cong [0,\infty]$ of metrics that are of product form near $\partial \mathsf{T}_i$, have positive scalar curvature on $\mathsf{T}_i$, and stretch the neck of $\partial \mathsf{T}_i$ as $q \in Q$ goes to $+\infty$.
Then we have
\begin{align*}
	\mathsf{F}_{\Gamma} \mathsf{G}_{\Gamma}\
	&= \Pi_\Gamma \circ \widetilde{\HMR}(\overline{W},\overline{\mathbf{a}}) \circ \widetilde{\HMR}(W,\mathbf{a}) \circ I_\Gamma \\
	&\simeq\Pi_\Gamma\thmr(X,\mathbf a\cup\bar{\mathbf a})I_\Gamma\\
	&\simeq\Pi_\Gamma\thmr\bigl(I\times(Y\sqcup Z),
	I\times(\mathbf p\sqcup\mathbf q\sqcup\{x_1,x_2,z_1,z_2\})\bigr)I_\Gamma\\
	&\simeq\Pi_\Gamma I_\Gamma=\id.
\end{align*}
The first and third chain homotopies follow from the composition law.
The second chain homotopy follows from neck stretching.
Indeed, counting family moduli spaces over $Q$ provides a homotopy from the cobordism map of $X$ to the cobordism map induced by $V$,
where we insert the relative real Seiberg--Witten invariant of $\mathsf{T}$ at the additional end.
Lemmas~\ref{lem:solution-on-C} and~\ref{lem:match-rrsc-count} show that the relative invariants of $\mathsf{T}$ and $\mathsf{B}$ coincide.
Here, we used the fact that there are no nonconstant trajectories on $\mathbb{R} \times \partial \mathsf{T}_i$.
Hence the resulting cobordism map is induced by the cylinder with product marking set.

Next, let $Q$ be a neck-stretching family of metrics on $X^{\dag} = V^{\dag} \cup (\mathsf{S}_1\sqcup \mathsf{S}_2)$ for the cut $\partial \mathsf{S}_1 \sqcup \partial \mathsf{S}_2$.
The composition $I_\Gamma \circ \Pi_\Gamma$ selects only the \RRSC\ structures on $Y$ whose markings on $x_1,x_2$ are transported from $\Gamma$.
By Lemma~\ref{lem:match-rrsc-count}, the \RRSC\ structures on $V^{\dag}$ that contribute to the image of $I_\Gamma \circ \Pi_\Gamma$ are precisely the ones that extend over $\mathsf{A}_1 \sqcup \mathsf{A}_2$.

Therefore,
\begin{align*}
	\mathsf{G}_\Gamma \mathsf{F}_\Gamma
	&= \widetilde{\HMR}(W,\mathbf{a}) \circ I_\Gamma \circ \Pi_\Gamma \circ \widetilde{\HMR}(\overline{W},\overline{\mathbf{a}})\\
	&\simeq
 \thmr\left(I\times(Y\#_2 Z),I\times(\mathbf p\sqcup\mathbf q)\right)
 \simeq\id.
\end{align*}
The first homotopy follows from the same neck stretching argument as above and Lemma~\ref{lem:solution-on-C}.
The second follows from naturality.
In the first step, recall $\mathsf{A}_i$ contains an empty marking set,
whereas $\mathsf{S}_i$ is marked by a core circle.
Therefore the real gauge groups on each $\mathsf{A}_i$ or $\partial \mathsf{A}_i$ contains the constant gauge transformation $\pm 1$.
Moreover, the gauge transformation $-1$ that extends from $\partial \mathsf{A}_i$ to  $\mathsf{A}_i$ identifies pairs of glued-up Seiberg--Witten solutions on $V^{\dag}\cup (\mathsf{A}_1 \sqcup \mathsf{A}_2)$.
On the other hand, $-1$ does not extend over $\mathsf{S}_i$ as a framed gauge transformation due to the interior marking set.
Since changing the gluing identification on $\partial \mathsf{S}_i$ by $-1$ negates the induced marking at $x_i$,
only one of the two corresponding \RRSC\ summands lies in the image of $I_\Gamma\Pi_\Gamma$.
Thus the count of Seiberg--Witten solutions on the cobordism $V^{\dag}\cup (\mathsf{A}_1 \sqcup \mathsf{A}_2)$ coincides with and those on $V^{\dag}\cup (\mathsf{S}_1 \sqcup \mathsf{S}_2)$ that factor through $I_{\Gamma} \circ \Pi_{\Gamma}$.
\end{proof}

\section{Alexander polynomial at $i$}
\label{sec:alex}

Let $L\subset S^3$ be an oriented link, and let $F$ be a connected Seifert surface.
Fix an ordering of the components of $L$.
Let $\mathbf{p}$ be a set of basepoints with exactly one basepoint on each component.
Denote the total Euler characteristic of the framed (or equivalently by~Proposition~\ref{prop:euler}, the cone) real monopole Floer homology and its normalization by:
\begin{equation*}
 \chi(L)=\chi\left(\thmr(\dbcv_2(S^3,L),\tau,\mathbf p)\right),
 \qquad
 \widehat\chi(L)=(-1)^{\epsilon(L)}\chi(L).
\end{equation*}
where $\epsilon(L)$ is the normalization term~\eqref{eq:alexander-normalization}.
That is, $\chi(L)$ uses the mod-two grading $\grtwo$, and $\widehat\chi(L)$ uses the normalized grading $\widehat{\gr}$.

\subsection{The third vertex}

\begin{figure}[htbp]
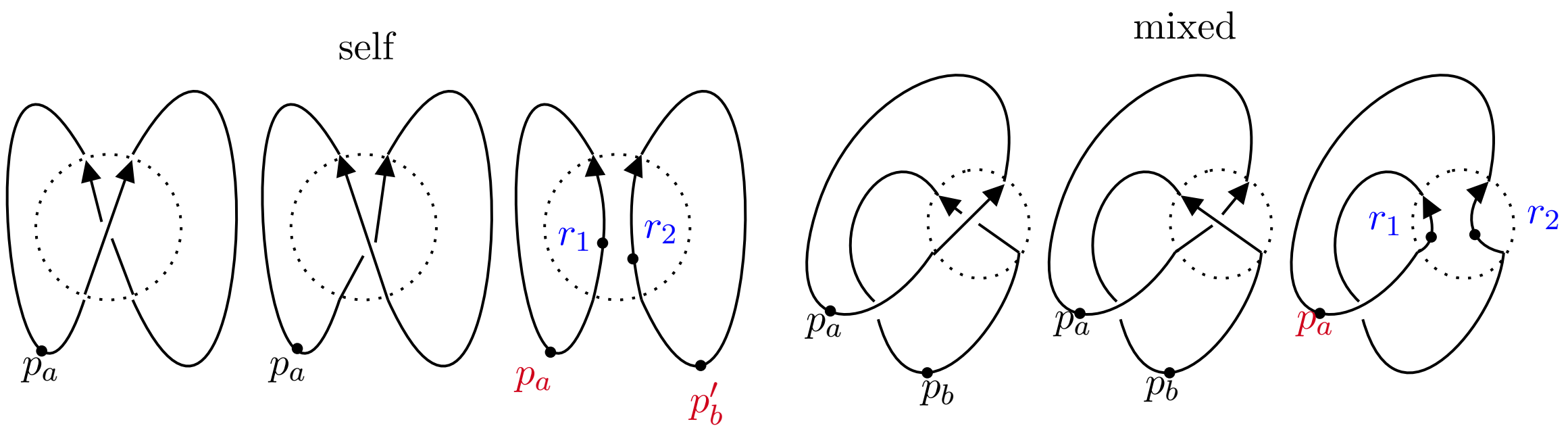

	\centering
	\figselfmixed
	\caption{Skein triples associated to a self crossing (left) and a mixed crossing (right).
	The basepoints $\mathbf{p}_0$ are placed on the $L_0$'s.}
	\label{fig:selfmixed}
\end{figure}

Let $\mathbf p$ be the inherited marking set on $Y_{00}$.
In the case of mixed crossing, let $p_a,p_b$ be the points on the two components meeting the crossing,
and set $\mathbf p_0=\mathbf p\setminus\{p_b\}$.
In the case of self crossing,
let $p_b'$ be a point away from the two balls in \eqref{eq:tube},
which lies on the component of $L_0$ containing no basepoint from $\mathbf{p}$.
Set $\mathbf p_0=\mathbf p\cup\{p_b'\}$.
In either case, $\mathbf p_0$ has one basepoint on every component of $L_0$. See Figure~\ref{fig:selfmixed}.

By \eqref{eq:tube}, we have $Y_{00}\cong Y_0\#_2 S^3$,
where the involution on $S^3$ fixes an unknot.
Denote by $r_1,r_2$ the points on $L_0$ used for the double connected sum, and by $z_1,z_2$ the corresponding points on $S^3$.
In the mixed crossing case, both $r_1,r_2$ lie on the component of $L_0$ containing $p_a$.
Inside $Y_{00}$, move $p_b$ into the $S^3$ side through the fixed point set.
Theorem~\ref{thm:connected-sum}, with $\mathbf q=\{p_b\}$ and a choice $\Gamma$ of arcs, gives a quasi-isomorphism
\begin{equation}\label{eq:mixed-connected-sum}
 \mathsf{F}_{\Gamma} \colon \tcmr(Y_{00},\mathbf p)
 \xrightarrow{\cong}
 \tcmr(Y_0,\mathbf p_0)\otimes
 \tcmr(S^3,\{z_1,z_2,p_b\}).
\end{equation}

In the self crossing case, $r_1,r_2$ lie on different components of $L_0$.
Since both components contain basepoints from $\mathbf p_0$,
Theorem~\ref{thm:connected-sum} with $\mathbf q= \emptyset$ gives a quasi-isomorphism
\begin{equation}\label{eq:self-connected-sum}
 \mathsf{F}_{\Gamma} \colon \tcmr(Y_{00},\mathbf p\cup\{p_b'\})
 \xrightarrow{\cong}
 \tcmr(Y_0,\mathbf p_0)\otimes
 \tcmr(S^3,\{z_1,z_2\}).
\end{equation}
On $Y_{00}$ the extra point $p_b'$ lies on a component already marked by $\mathbf p$.
In this case, we use the grading on $(Y_{00},\mathbf p)$ for the group on the left-hand side.
Forgetting the extra marking $p_b'$ does not change the framed gauge groups, configuration spaces, nor the perturbed Seiberg--Witten operators.
The Euler characteristic doubles:
\begin{equation}
\label{eq:factor-2-marking}
 \chi\left(\thmr(Y_{00},\mathbf p\cup\{p_b'\})\right)
 =2\chi\left(\thmr(Y_{00},\mathbf p)\right).
\end{equation}
Indeed, for each \RRSC\ structure on $(Y_{00},\mathbf p)$,
the two choices for an additional marking $T_{p_b'}$ at $p_b'$ are inequivalent.
This is because a real gauge transformation preserving the markings on $\mathbf{p}$
must evaluate to $+1$ at $p_b'$.
For the standard real $S^3$ with $k$ marked points on the fixed unknot, there are $2^{k-1}$ \RRSC\ structures.

\subsection{The unnormalized skein relation}
Let us compute the absolute mod-two degree of the splitting map $\mathsf{F}_{\Gamma}$.
Let $\mathfrak a$ and $\mathfrak b\otimes\vartheta$ be homogeneous representatives of nonzero elements of the left and right sides of \eqref{eq:mixed-connected-sum} or \eqref{eq:self-connected-sum}, respectively.
Suppose $\mathsf F_\Gamma([\mathfrak a]) = [\mathfrak b]\otimes[\vartheta]$.
Denote by $\widehat{\mathfrak s}_{\mathfrak t_{00}}$ and $\widehat{\mathfrak s}_{\mathfrak t_0}$ the closed surface spin
structures used in the mod-two grading of $\mathfrak a$ and $\mathfrak b$, on the closures $\hat{F}_{00}$ and $\hat{F}_0$, respectively.
\begin{lem}\label{lem:splitting-degree}
The degree of the splitting map $\mathsf{F}_{\Gamma}$ is given in terms of the Arf invariants on $F_0$ and $F_{00}$; that is,
\begin{equation}\label{eq:splitting-degree}
 \gr^{(2)}(\mathfrak a)-\gr^{(2)}(\mathfrak b)
 \equiv
 \Arf(\widehat{\mathfrak s}_{\mathfrak t_{00}})
 -\Arf(\widehat{\mathfrak s}_{\mathfrak t_0})\pmod2.
\end{equation}
\end{lem}
\begin{figure}[htbp]
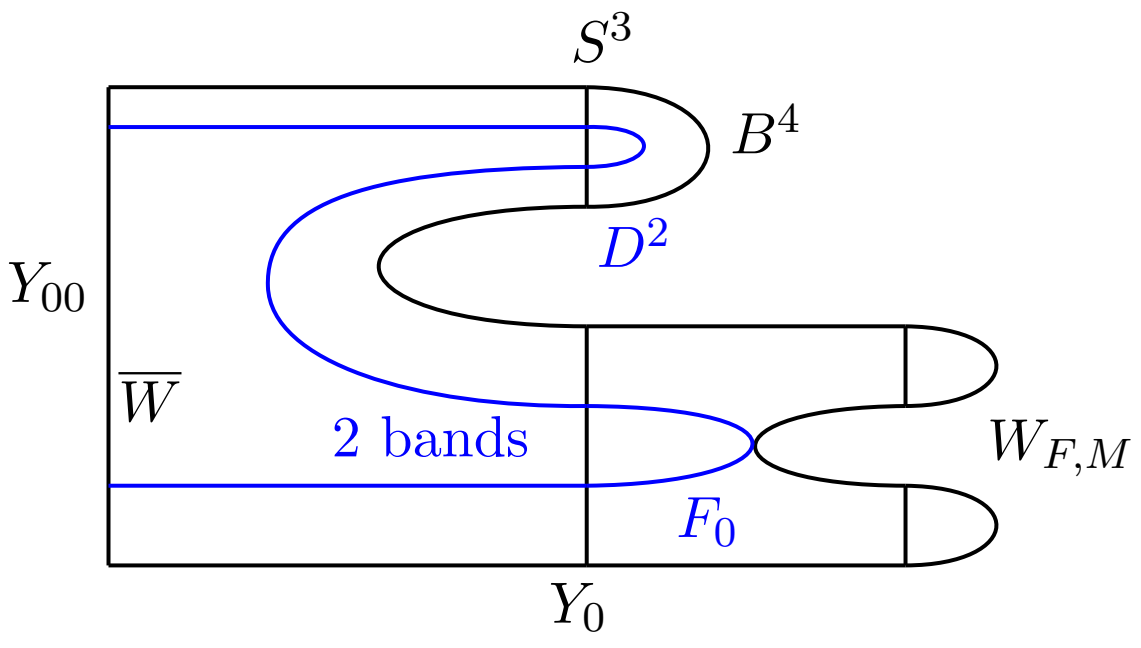

	\centering
	\figFgamma
	\caption{Seifert cap of $Y_{00}$ from Seifert cap of $Y_0$.
	The (single) curved blue arc on the left schematically represents the two bands that join $F^2$ and $D^2$.}
	\label{fig:Fgamma}
\end{figure}

\begin{figure}[htbp]
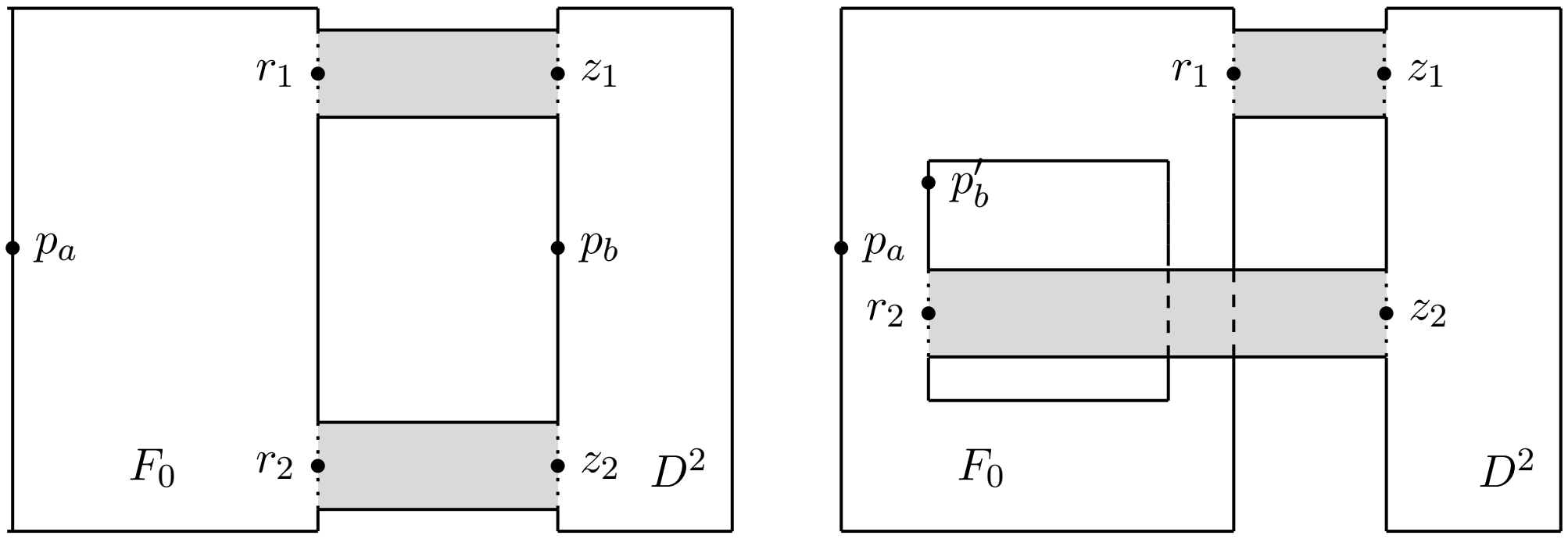

	\centering
	\figbands
	\caption{Attach a pair of bands (shaded) to obtain $F_{00}$ from $F_0$,
	in the mixed crossing (left) and self crossing (right) cases. The $D^2$'s refer to the discs bounding the unknot in $S^3$.
	Equivalently, we regard $F_{00}$ as obtained from attaching a single band to $F_{0}$, by viewing the whole right hand side $D^2 \sqcup \{\text{2 bands}\}$ as a band.}
	\label{fig:bandstoF}
\end{figure}

\begin{proof}
Let $W_{F_0,B^4}$ and the standard real $B^4$ be the Seifert caps of $Y_0$ and $S^3$, respectively.
Also, let $\overline{W} \colon Y_{00} \to (Y_0 \sqcup S^3)$ be the splitting cobordism (\S\ref{sec:connected-sum}).
Capping $\overline{W}$ by $W_{F_0,B^4}$ and $B^4$ gives rise to a cap $X$ of $Y_{00}$:
\[
X=\overline{W} \cup\bigl(W_{F_0,B^4}\sqcup B^4\bigr).
\]
Equivalently, $X$ is obtained from $W_{F_0,B^4}\sqcup B^4$ by attaching, for each $i=1,2$,
a $1$-handle joining the boundary ball at $r_i$ to the boundary ball at $z_i$.
In the quotient, the two 4-balls $W_{F_0,B^4}/\tau$ and $B^4/\tau$ joined by these two 1-handles form an
$S^1\times B^3$.
The branch locus is $F_0\sqcup D^2$ with two bands attached,
as in Figure~\ref{fig:Fgamma}.
Here, $D^2$ is the fixed set in the standard $B^4$.
The union of the two bands with $D^2$ can be viewed as a single band attached to $F_0$ along the two original
attaching intervals,
as in Figure~\ref{fig:bandstoF}.
This identifies the resulting surface and its normal framing, with $F_{00}$ and the normal framing\footnote{See~\cite[\S3.1]{LM} for the definition of the normal framing $\uptau$ of a pushed-in Seifert surface.} $\uptau$ of $F_{00}$.
Thus $X$ can be regarded as the Seifert cap of $Y_{00}$ built using the outgoing cap $S^1\times B^3$.
This capping manifold does not have trivial $H_1$,
but it suffices for this index comparison (cf.~Remark~\ref{rem:inadmissible}).
Indeed, replacing the $S^1 \times B^3$ with $D^2\times S^2$ yields the same index of real Seiberg--Witten operator modulo two.

Furthermore,
since each of the two 1-handles in $X$ that are attached to $W_{F_0,B^4} \sqcup B^4$ above is $\tau$-invariant, we have
\[
 H^*(X,W_{F_0,B^4}\sqcup B^4;\mathbb{R})^{-\tau^*}=0,
 \qquad
 \iota_R(X)=\iota_R(W_{F_0,B^4})+\iota_R(B^4).
\]
For the reducible $\vartheta$ from the spin connection and the spin structure over the standard $B^4$, both $\gr_{SW}(\vartheta;B^4)$ and $\iota_R(B^4)$ are zero.
Capping the disc by another disc gives a sphere, which has zero Arf invariant.
Since $F_\Gamma$ counts a zero-dimensional moduli space,
additivity of the APS Seiberg--Witten index\footnote{Here the extending spin$^c$ structures are suppressed from the notation.} implies
\[\gr_{SW}(\mathfrak a;X)
 =\gr_{SW}(\mathfrak b;W_{F_0,B^4})
  +\gr_{SW}(\vartheta;B^4).\]
Thus both the Seiberg--Witten formal dimension terms and the $\iota_R$ terms in $\grtwo(\mathfrak a)-\grtwo(\mathfrak b)$ cancel out.
The remaining terms form Equation~\eqref{eq:splitting-degree}.
\end{proof}

With the degree of the splitting map understood,
we can now compute the mod-two gradings of the framed real monopole Floer homology of $Y_0$ from that on $Y_{00}$.
The computation of the Arf invariants in the proof of the following statement is closely related to the computation of the gradings on the framed real monopole Floer homology of $S^1 \times S^2$ in \cite[\S5]{LM}.
This is natural, as $Y_{00} \cong Y_0 \# (S^1 \times S^2)$.

\begin{prop}\label{prop:euler}
With the inherited basepoint set $\mathbf p$ and Seifert surface $F_{00}$, the total Euler characteristic for $Y_{00}$ satisfies
	\begin{equation}\label{eq:euler-result}
		\chi\left(\thmr(Y_{00},\mathbf p)\right)=
		\begin{cases}
			2\chi(L_0),&\text{if the crossing is mixed};\\
			\chi(L_0),&\text{if the crossing is self}.
		\end{cases}
	\end{equation}
\end{prop}

\begin{proof}
The Seifert surface $F_{00}$ for $Y_{00}$ is obtained from $F_0$ by attaching a band.
Since the spin structures $\widehat{\mathfrak s}_{\mathfrak t_{00}}$ and $\widehat{\mathfrak s}_{\mathfrak t_0}$ agree on $F_0$ and $F_{00}$ away from the bands,
the resulting quadratic forms agree on the complement of the band.
Denote by $q$ the quadratic form associated to $\widehat{\mathfrak s}_{\mathfrak{t}_{00}}$ on the closure of $F_{00}$.

Suppose the crossing is \emph{mixed}.
The two attaching arcs of the band both lie on the component $L_{0a}$ of $L_0 = \partial F_0$ that contains $p_a$,
obtained by smoothing the crossing.
Band attachement splits $L_{0a}$ into two components $L_{a},L_{b}$ of $\partial F_{00}$, containing basepoints $p_a$ and $p_b$ respectively.
Fix an \RRSC\ structure on $(Y_0,\tau,\mathbf p_0)$,
extend its marking $T$ to $r_1,r_2$ using the arcs and Clifford trivializations along $\Gamma_1,\Gamma_2$ as in \S\ref{sec:connected-sum}.
Modulo overall sign,
the four marking choices in \eqref{eq:mixed-connected-sum} are obtained by changing the markings at $z_2$ and $p_b$,
while fixing the marking at $z_1$.
The effects are given by marking changes $\e_{z_2}$ and $\e_{p_b}$ as in Definition~\ref{defn:marking-change}.

Changing the marking at $z_2$ multiplies the spin identification used to attach the band at that end of the second 1-handle by $-1$.
Since each resulting boundary circle meets the attaching arc containing $z_2$ at exactly one endpoint,
its spin monodromy is multiplied by $-1$.
Thus the marking change $\e_{z_2}$ shifts the boundary spin types $\lambda_a,\lambda_b$ of $\widehat{\mathfrak s}_{\mathfrak{t}_{00}}$ on the components $L_a,L_b$ of $\partial F_{00}$ via simultaneous translation:
\begin{equation}
\label{eq:lambda-ab-shift}
	(\lambda_a,\lambda_b)\longmapsto(\lambda_a+1,\lambda_b+1).
\end{equation}
The marking change $\e_{p_b}$ at $p_b$ preserves the spin identifications used to attach the bands,
and preserves the boundary spin types of $F_{00}$.
Since the cobordism from the circle $L_{0a}$ in $Y_0$
to the two circles $L_a,L_b$ in $Y_{00}$ is a pair of pants,
we have $\lambda_a^0=\lambda_a+\lambda_b$.

In sum, the four choices of markings $(\pm T_{z_1}, \pm T_{z_2})$ form two pairs of choices, interchanged by $\e_{z_2}$.
Within each pair, $\e_{p_b}$ swaps the two choices that have the same boundary type.
Assume the boundary types of $\widehat{\mathfrak s}_{\mathfrak t_{00}}$ and $\widehat{\mathfrak s}_{\mathfrak t_{0}}$ are zero except for possibly $\lambda_a,\lambda_b,\lambda_a^0$.
Then by Corollary~\ref{cor:marking-change-euler},
only the choices with $(\lambda_a,\lambda_b) = (0,0)$ contribute nontrivially to the total
Euler characteristic for $Y_{00}$.
By \eqref{eq:lambda-ab-shift}, exactly two choices satisfy this when $\lambda_a^0=0$, and none do when $\lambda_a^0=1$.

Since the band splits one boundary component into two, the inclusion $F_0\subset F_{00}$ induces an isomorphism on $H_1$ modulo boundary classes.
Thus the quadratic forms of $\mathfrak s_{\mathfrak{t}_{00}}$ and $\mathfrak s_{\mathfrak{t}_{0}}$ agree modulo boundary classes.
For either contributing choice of markings on $(z_1,z_2)$ above, all boundary types of $\mathfrak s_{\mathfrak{t}_{00}}$ vanish, so $\Arf(\widehat{\mathfrak s}_{\mathfrak t_{00}})
 =\Arf(\widehat{\mathfrak s}_{\mathfrak t_0})$.
By \eqref{eq:splitting-degree}, the degree of $\mathsf{F}_{\Gamma}$ restricted to the two contributing summands is zero. Therefore, in the mixed case,
\[\chi\left(\thmr(Y_{00},\mathbf p)\right)= 2\chi(L_0).\]

Suppose the crossing is \emph{self}.
The band joins two boundary components $L_{0a}, L_{0b}$ of $F_0$.
Hence their boundary spin types $\lambda_a^0,\lambda_b^0$ satisfy
$\lambda_a=\lambda_a^0+\lambda_b^0$,
where $\lambda_a$ is the spin type on the joined component of $\partial F_{00}$ containing $p_a$.
By Corollary~\ref{cor:marking-change-euler} again,
in order for spin structure $\mathfrak s_{\mathfrak t_{00}}$ to contribute to the Euler characteristic, $\lambda_a = 0$.
Thus
the boundary spin structure types on $F_0$ must satisfy $\lambda_a^0=\lambda_b^0 =\varepsilon$, and zero on every other boundary component.

Let $c \in H_1(\widehat{F}_{00})$ represent one of the two boundary components $L_{0a},L_{0b}$ of $F_0$.
Represent $d$ by the core of the added band together with an arc in $F_0$ joining its ends.
Then $c \cdot d = 1$, and $q(c)=\varepsilon$.
The two \RRSC\ structures of $(S^3, \{z_1,z_2\})$ differ by changing the marking at $z_2$ while fixing the marking at $z_1$.
This multiplies the spin identification of band attachement by $-1$.
Hence $q(d)$ changes by one.
The curve representing $c$ and the remaining symplectic basis can be chosen away from the band, where the spin structures agree.
Thus the Arf invariant of $\widehat{s}_{\mathfrak{t}_{00}}$ changes by
$q(c) \cdot 1 =\varepsilon$.
This implies that,
when $\varepsilon=1$, the two \RRSC\ structures on $(S^3,\{z_1,z_2\})$ induce opposite degrees of $\mathsf F_\Gamma$ by Lemma~\ref{lem:splitting-degree}, whose contributions cancel.

Suppose $\varepsilon=0$.
In $\hat F_0$, a curve $d$ dual to $c$ consists of an arc in $F_0$ joining the boundary components $L_{0a},L_{0b}$, and closed by an arc in the planar cap.
In $\hat F_{00}$, define $d$ using the same arc in $F_0$, closed by the core of the added band.
Choose the planar caps so that the quadratic forms $q_0$ and $q$ agree on the symplectic complements of $(c,d)$.
Since $q(c)=q_0(c) = 0$, again $\Arf(\widehat{\mathfrak s}_{\mathfrak t_{00}})
 =\Arf(\widehat{\mathfrak s}_{\mathfrak t_0})$.
It follows that
\[\chi\left(\thmr(Y_{00},\mathbf p \sqcup \{p_b'\})\right)= 2\chi(L_0).\]
Dividing this by two via~\eqref{eq:factor-2-marking} implies the self case of \eqref{eq:euler-result}.
\end{proof}

\subsection{The normalized skein relation}

\begin{prop}
The normalized Euler characteristic satisfies\footnote{Compare with \cite[Proof of Proposition 1.6]{Srivastava}.}
\begin{equation}\label{eq:normalized-skein}
 \widehat\chi(L_+)-\widehat\chi(L_-)=
 \begin{cases}
  2\widehat\chi(L_0),&\text{if the crossing is mixed};\\
  -\widehat\chi(L_0),&\text{if the crossing is self}.
 \end{cases}
\end{equation}
\end{prop}

\begin{proof}
Let $n=|L_\pm|$ and $n_0 = |L_0|$.
By the normalization formula~\eqref{eq:alexander-normalization} and
Theorem~\ref{thm:graded-triangle} we have, modulo two,
\[
 \delta_1=\epsilon(L_+)-\epsilon(L_-),\qquad
 \delta_2=\epsilon(L_-) -\epsilon(L_0)+\frac{1+n-n_0}{2}.
\]
Multiplying \eqref{eq:graded-euler} from Corollary~\ref{cor:euler-char-equality} by
$(-1)^{\epsilon(L_+)}$ yields
\[
 \widehat\chi(L_+)-\widehat\chi(L_-)
 =-(-1)^{(1+n-n_0)/2} (-1)^{\epsilon(L_0)}
       \chi\left(\thmr(Y_{00},\mathbf p)\right).
\]
A mixed crossing satisfies $n_0=n-1$, and $n_0=n+1$ otherwise.
Proposition~\ref{prop:euler} now implies~\eqref{eq:normalized-skein}.
\end{proof}

For the unknot with one basepoint, \cite[\S5.1]{LM} shows $\chi(U)=\widehat\chi(U)=1$.
Applying the normalized skein relation~\eqref{eq:normalized-skein} to the triple $(U_n,U_n,U_{n+1})$ iteratively yields the following.
\begin{cor}
	Let $U_n$ be the $n$-component unlink.
	Then $\widehat\chi(U_1) = 1$, and $\widehat\chi(U_n) = 0$ for all $n \ge 2$.
\end{cor}
\begin{rem}
Alternatively, one can show that $\thmr_*(\dbcv_2(U_n)) \cong H_*(T^{n-1})$, which has zero Euler characteristic for all $n \ge 2$, using a PSC metric on $\dbcv_2(U_n) \cong \#^{n-1}(S^1 \times S^2)$.
See \cite[\S14.1]{Li22}.
\end{rem}

\subsection{Alexander polynomial evaluation}
Let $i=\sqrt{-1}$.
Let $\nabla_L(z)$ be the Conway polynomial with the convention that
\begin{equation}\label{eq:conway-convention}
	\nabla_U(z)=1,\qquad
	\nabla_{L_+}(z)-\nabla_{L_-}(z)=z\nabla_{L_0}(z).
\end{equation}
The Alexander polynomial is
\begin{equation*}
	\Delta_L(t)=\nabla_L(t^{1/2}-t^{-1/2}).
\end{equation*}
If $L$ is a knot, this is the symmetric Alexander polynomial normalized by $\Delta_L(1)=1$.
If $L$ is a link, then $\Delta_{L}(t)$ may have half-integral powers.
At $t=i$ we take $t^{1/2}=e^{\pi i/4}$, so that $t^{1/2}-t^{-1/2}=i\sqrt2$.
Moreover, $\Delta_{L}(t)$ vanishes for unlinks with more than one component.

\begin{thm}
For every oriented link $L\subset S^3$,
\begin{equation}\label{eq:alexander-result}
	\widehat\chi(L)=(-i\sqrt2)^{|L|-1}\Delta_L(i).
\end{equation}
\end{thm}

\begin{proof}
The link invariant
\begin{equation}\label{eq:alex-value-skein}
	(-i\sqrt2)^{|L|-1}\nabla_L(i\sqrt2)
\end{equation}
has the same unlink values as $\widehat\chi$.
At a mixed crossing the component factor in its oriented skein relation is $(-i\sqrt2)$,
	and at a self crossing it is $(-i\sqrt2)^{-1}$.
By \eqref{eq:conway-convention}, its skein coefficients for the mixed and self crossings are
\[(-i\sqrt2)(i\sqrt2)=2,\qquad (-i\sqrt2)^{-1}(i\sqrt2)=-1,\]
respectively.
Since \eqref{eq:alex-value-skein} and $\widehat{\chi}(L)$ satisfy the same skein relation \eqref{eq:normalized-skein} and have the same evaluations on unlinks,
Equation~\eqref{eq:alexander-result} follows by induction.
\end{proof}

\bibliographystyle{plain}
\bibliography{yat.bib}
\end{document}